\documentclass[a4paper]{article}
\usepackage{latexsym, amsfonts, amsmath, amssymb, amsthm}
\usepackage{hyperref}
\hypersetup{
	colorlinks=true,
	linkcolor=red,
	citecolor=blue,
}
\usepackage{xcolor}
\usepackage[
margin=1in,
]{geometry}
\newcommand{\be}{\begin{equation}}
\newcommand{\ee}{\end{equation}}
\newcommand{\bea}{\begin{eqnarray}}
\newcommand{\eea}{\end{eqnarray}}
\newcommand{\bean}{\begin{eqnarray*}}
\newcommand{\eean}{\end{eqnarray*}}

\newcommand{\brray}{\begin{array}}
\newcommand{\erray}{\end{array}}
\newcommand{\ben}{\begin{equation}{nonumber}}
\newcommand{\een}{\end{equation}{nonumber}}

\newtheorem{dfn}{Definition}[section]
\newtheorem{thm}[dfn]{Theorem}
\newtheorem{lmma}[dfn]{Lemma}
\newtheorem{ppsn}[dfn]{Proposition}
\newtheorem{crlre}[dfn]{Corollary}
\newtheorem{xmpl}[dfn]{Example}
\newtheorem{rmrk}[dfn]{Remark}
\newcommand{\bdfn}{\begin{dfn}}
\newcommand{\bthm}{\begin{thm}}
\newcommand{\blmma}{\begin{lmma}}
\newcommand{\bppsn}{\begin{ppsn}}
\newcommand{\bcrlre}{\begin{crlre}}
\newcommand{\bxmpl}{\begin{xmpl}}
\newcommand{\brmrk}{\begin{rmrk}}
\newcommand{\edfn}{\end{dfn}}
\newcommand{\ethm}{\end{thm}}
\newcommand{\elmma}{\end{lmma}}
\newcommand{\eppsn}{\end{ppsn}}
\newcommand{\ecrlre}{\end{crlre}}
\newcommand{\exmpl}{\end{xmpl}}
\newcommand{\ermrk}{\end{rmrk}}

\newcommand{\IC}{\mathbb{C}}

\newcommand{\IR}{\mathbb{R}}

\newcommand{\IT}{\mathbb{T}}

\newcommand{\IZ}{\mathbb{Z}}

\newcommand{\oneformclassical}{\Omega^1 ( M )}
\newcommand{\oneform}{{\Omega}^1_D ( \mathcal{A} )}

\newcommand{\twoformdeformed}{{\Omega}^2_D ( \mathcal{A}_\theta )}
\newcommand{\twoform}{{{\Omega}^2_D}( \mathcal{A} )}
\newcommand{\tensora}{\otimes_{\mathcal{A}}}
\newcommand{\tensorsym}{\otimes^{{\rm sym}}_{\mathcal{A}}}
\newcommand{\tensorc}{\otimes_{\mathbb{C}}}
\newcommand{\A}{\mathcal{A}}
\newcommand{\Atheta}{\mathcal{A}_{\theta}}
\newcommand{\B}{\mathcal{B}}

\newcommand{\C}{\mathcal{C}}

\newcommand{\D}{\mathcal{D}}

\newcommand{\E}{\mathcal{E}}
\newcommand{\Etheta}{\mathcal{E}_{\theta}}
\newcommand{\F}{\mathcal{F}}

\newcommand{\G}{\mathcal{G}}

\newcommand{\Acenter}{\mathcal{Z}( \mathcal{A} )}

\newcommand{\Ecenter}{\mathcal{Z}( \mathcal{E} )}

\newcommand{\Aprime}{\mathcal{A}^{\prime}}

\newcommand{\Vgtwo}{V_{g^{(2)}}}
\newcommand{\Psym}{P_{\rm sym}}

\newcommand{\cla}{{\cal A}}
\newcommand{\clb}{{\cal B}}
\newcommand{\clc}{{\cal C}}

\newcommand{\clf}{{\cal F}}

\newcommand{\clh}{{\cal H}}

\newcommand{\clj}{{\cal J}}

\newcommand{\cll}{{\cal L}}

\newcommand{\cls}{{\cal S}}

\newcommand{\clz}{{\cal Z}}

\def\a*{{\cal A}_{h,*}}
\def\B{{\cal B}(h)}
\def\B1{{\cal B}_1(h)}
\def\b{{\cal B}^{\rm s.a.}(h)}
\def\b1{{\cal B}^{\rm s.a.}_1(h)}

\newcommand{\raro}{\rightarrow}

\def \qed {$\Box$}

\def\a*{{\cal A}_{h,*}}
\def\B{{\cal B}(h)}
\def\B1{{\cal B}_1(h)}
\def\b{{\cal B}^{\rm s.a.}(h)}
\def\b1{{\cal B}^{\rm s.a.}_1(h)}

\newcommand{\RNum}[1]{\uppercase\expandafter{\romannumeral #1\relax}}

\begin{document}
\begin{center}
{\Large{\bf Levi-Civita connections for a class of spectral triples}}\\
\vspace{0.2in}
{\large {Jyotishman Bhowmick, Debashish Goswami and Sugato Mukhopadhyay}}\\
Indian Statistical Institute\\
203, B. T. Road, Kolkata 700108\\
Emails: jyotishmanb$@$gmail.com, goswamid$@$isical.ac.in, m.xugato@gmail.com \\
\end{center}
\begin{abstract}
We give a new definition of Levi-Civita connection for a noncommutative pseudo-Riemannian metric on a noncommutative manifold given by a spectral triple. We prove 
the existence-uniqueness result for a class of modules of one-forms over a large class of noncommutative manifolds, including the matrix geometry of the fuzzy 3-sphere, the quantum Heisenberg manifolds and Connes-Landi deformations of spectral triples on the Connes-Dubois Violette-Rieffel-deformation of a compact manifold equipped with a free toral action. It is interesting to note that in the example of the quantum Heisenberg manifold, the definition of metric compatibility given in \cite{frolich} failed to ensure the existence of a unique Levi-Civita connection.  In the case of the matrix geometry, the Levi-Civita connection that we get coincides with the unique real torsion-less unitary connection obtained by Frolich et al in \cite{frolich}. 
\end{abstract}

\section{Introduction} 
The concepts of connection and curvature occupy a central place in any form of geometry, classical or noncommutative. In noncommutative geometry (NCG for short) \`{a} la Connes, there have been several attempts over the last few years to formulate and study analogues of curvature. There seem to be mainly two different approaches to this problem so far:\\
(a) formulating an analogue of Levi-Civita connection and computing the corresponding curvature operator, in particular scalar and Ricci curvatures
(see, e.g. \cite{frolich}, \cite{connes_1} );

or\\
(b) defining Ricci and scalar curvature through an asymptotic expansion of the noncommutative Laplacian
(see \cite{Connes_moscovici}, \cite{dabrowski_sitarz}, \cite{dabrowski_sitarz2}, \cite{scalar_1}, \cite{scalar_2}, \cite{fourtorus}, \cite{scalar_3}, \cite{scalar_4}, \cite{scalar_5}, \cite{fatizadeh_connes}, \cite{khalkhaliricci} etc.).

In the classical case, at least for a compact Riemannian manifold, these two approaches turn out to be equivalent.
In the first approach, one gets the Levi-Civita connection and the full curvature operator, which are important in their own right. 

In NCG, the definition of Levi-Civita connection given in \cite{frolich} seemed to face an obstacle because it was not possible to get a unique Levi-Civita (i.e. both torsion-less and metric-compatible in a suitable sense) connection in some standard examples such as the fuzzy 3-sphere and the quantum Heisenberg manifolds. We refer to \cite{frolich}, \cite{chak_sinha} as well as the appendix B of \cite{heckenberger_etal} for results regarding such non-existence/non-uniqueness.

In a recent article, J. Rosenberg (\cite{Rosenberg}) proposed an alternative definition giving an existence and uniqueness theorem for some noncommutative manifolds including the noncommutative tori. This has been followed by existence-uniqueness results for Levi-Civita connections and  computations of scalar curvature by a number of authors (e.g. \cite{pseudo,arnlind2,arnlind3,sheu}). In \cite{pseudo}, the definition of Rosenberg has been extended to the case of (noncommutative) pseudo-Riemannian metrics as well.

The aim of this article is to give a new definition of Levi-Civita connections for the space of one-forms of a spectral triple. The underlying pseudo-Riemannian metric is a bilinear version of the sesquilinear form constructed in Theorem 2.9 of \cite{frolich}. In fact, there are two points of departure of this article from the existing literature on Levi-Civita connections on spectral triples.

Firstly, as opposed to the approach taken in \cite{pseudo}, \cite{sheu}, \cite{Rosenberg} etc., we use one-forms instead of derivations. Indeed, spaces of vector fields in NCG are not as well-behaved as in the case of classical geometry. In fact, they do not form a module over the underlying (noncommutative) smooth algebra. The simple description of vector fields, i.e. derivations of the smooth algebra of the noncommutative tori, 
which played a crucial role in the success of Rosenberg's approach for such noncommutative manifolds, hold only for spectral triples equivariant w.r.t. a 
toral action which is also ergodic on the underlying $C^*$ algebra. On the other hand, the space of one-forms is quite well-behaved in NCG and it does have a natural
module structure over the noncommutative algebra of smooth functions. In algebra and algebraic geometry, including noncommutative algebraic geometry, the notion of connections on the module of one-forms of a differential calculus is quite familiar and standard (See \cite{frolich}, \cite{heckenberger_etal}, \cite{landi} etc.).

The second point of departure is the definition of metric compatibility of a connection. For a spectral triple $(\A,\mathcal{H},D)$ and $\oneform$ as the associated space of one-forms (see Subsection \ref{sptripleforms}), we work with an $\A$-bilinear map from $\oneform \tensora \oneform$ to $\A$ as the candidate for a pseudo-Riemannian metric. This is sharp contrast to the approach taken in \cite{frolich} where the authors had worked with a sesquilinear form. As a result, we work with a different definition of metric compatibility of a connection (see Subsection \ref{metric_compatibility}) compared to that in \cite{frolich}.

In some sense, we have combined the approaches of \cite{frolich} and \cite{Rosenberg} as we work in the setting of one-forms instead of vector fields but define a pseudo-Riemannian-metric to be a symmetric, bilinear non-degenerate form instead of a sesquilinear inner product. Under some assumptions on the space of one-forms of a spectral triple, we prove the existence and uniqueness of a Levi-Civita connection. 

Our assumptions are satisfied by a large class of noncommutative manifolds,
 which do include all the Rieffel-deformations of classical compact Riemannian manifolds obtained by isometric and free toral actions.
 Moreover, our existence-uniqueness result covers the case of quantum Heisenberg manifold for which the approach of \cite{frolich} did not succeed. Interestingly, it turns out that the Levi-Civita connection for the Heisenberg manifold as obtained from this article has constant negative scalar curvature (see Section \ref{heisenberguniqueness}). For the matrix geometry on the fuzzy 3-sphere, the authors of \cite{frolich} prove that there exists a family of torsion-less and unitary connections. Uniqueness can be obtained if in addition, one assumes that the connection is real. With our definition of metric compatibility, we get a unique torsion-less and metric compatible connection in this example which coincides with the unique real torsion-less unitary connection obtained in \cite{frolich}. 
 
  The results of this article apply to a class of pseudo-Riemannian bilinear metrics on the space of one-forms of a spectral triple $(\A,\mathcal{H},D)$ satisfying certain assumptions discussed in this article. Our results do not cover the examples of the conformal perturbations of a Riemannian bilinear metric which is the subject of study of a number of recent works (\cite{scalar_3}, \cite{fatizadeh_connes}, \cite{dabrowski_sitarz}, \cite{dabrowski_sitarz2}, \cite{scalar_1}, \cite{scalar_2}, \cite{fourtorus}, \cite{khalkhaliricci}, \cite{scalar_4}, \cite{scalar_5}). However, in a companion article (\cite{article2}), it is proven that under one additional technical assumption on the spectral triple as in this paper, there exists a unique Levi-Civita connection for any pseudo-Riemannian metric (which is only right $\A$-linear as opposed to being both left and right $\A$-linear). Thus one can compute the Ricci and scalar curvature for the conformal perturbation of the canonical metric on the noncommutative 2-torus (see \cite{article2}).
  
  Although the existence and uniqueness theorem of \cite{article2} hold under weaker assumptions, the novelty of this article based only on the spectral triple framework is twofold. Firstly, given a spectral triple, we provide a natural candidate for the Riemannian metric on the space of one-forms and then check that our candidate indeed satisfies the required properties (see Definition \ref{11thjuly20183}). The second novelty is in the proof of Theorem \ref{Phigiso} which nicely adapts the classical proof of existence and uniqueness of Levi-Civita connections by exploiting the bilinearity of $g$ and the isomorphism
$$ ( \Psym )_{12}: \oneform \tensora ( \oneform \tensorsym \oneform ) \rightarrow ( \oneform \tensorsym \oneform ) \tensora \oneform. $$
Here, the symbols $ \Psym $ and $ \tensorsym $ are as in Subsection \ref{splittingsubsection}. Let us mention that in \cite{article3}, an alternative proof of existence and uniqueness of Levi-Civita connections has been derived in the set up of this article. The proof in \cite{article3} mirrors exactly the one in the classical case by deriving a Koszul-formula on the space of forms. 
    
   It should also be mentioned that there is a very different approach by S. Majid, E. Beggs and their co-authors (\cite{majid_1}, \cite{majid_2}, \cite{majid_3}, \cite{majid_4} etc. and references therein). In \cite{article3}, it has been shown that our Levi-Civita connection is a bimodule connection with respect to a very natural braiding-like operator. 
	
	Let us discuss the plan of the article. We begin with some generalities on bimodules and spectral triples. In particular, we focus on a certain class of bimodules called `centered bimodules' and the flip map on the tensor product of two copies of such bimodules. In Section \ref{4thjune2018}, we start by discussing pseudo-Riemannian metrics. Then we recall the $A''$-valued sesquilinear form on the bimodule of forms for a spectral triple $(\A,\mathcal{H},D)$ constructed in \cite{frolich} and show that under some regularity assumptions on the spectral triple, one can construct an $\A$-valued pseudo-Riemannian bilinear metric. In Section \ref{mainsection}, we define the metric compatibility of a connection on $\oneform$ for a class of spectral triples and then state and prove the main result giving the existence and uniqueness of the Levi-Civita connection.
	
	The rest of this article deals with examples. In Section \ref{fuzzysection} and Section \ref{heisenberguniqueness}, we apply our results to spectral triples on the fuzzy 3-sphere and quantum Heisenberg manifold defined and studied in \cite{frolich} and \cite{chak_sinha} respectively. Finally, Section \ref{prelimdef} is devoted to the verification of our assumptions for the one-forms of the Connes-Landi spectral triples on a large class of noncommutative manifolds obtained by Rieffel deformation of classical Riemannian manifolds.

We fix some notations which we will follow. Throughout the article, $ \mathcal{A} $ will denote a complex unital $ \ast $-subalgebra of a $ C^* $ algebra. $ \mathcal{Z} ( A ) $ will denote the center of $ \cla. $ $ {\rm Span}(V) $ will denote the complex linear span of a subset $V$ of a complex vector space. Moreover, for a subset $ S $ of a right $ \A $ module $ \E, ~ S \A $ will denote $ {\rm Span} \{ sa: s \in S, ~ a \in \A \}. $ For right $ \A $ modules $ \mathcal{E} $ and $ \mathcal{F}, ~ {\rm Hom}_\A ( \mathcal{E}, \mathcal{F} ) $ will denote the set of all right $ \A $-linear maps from $ \E $ to $ \F. $ The tensor product over the complex numbers $ \IC $ is denoted by $ \tensorc $ while the notation $\tensora$ will denote the tensor product over the algebra $ \A. $  

If $ \mathcal{E} $ and $ \mathcal{F} $ are bimodules, $ {\rm Hom}_\A ( \mathcal{E}, \mathcal{F} ) $ has a left $ \A $ module structures given by left multiplication by elements of $ \A, $ i.e, for elements $ a $ in $ \A, ~ e $ in $ \E $ and $ T $ in $ {\rm Hom}_\A ( \mathcal{E}, \mathcal{F} ) , ~ ( a T ) ( e ) := a T ( e ) \in \F. $ The right $ \A $ module structure on $ {\rm Hom}_\A ( \mathcal{E}, \mathcal{F} ) $ is given by $ T a ( e ) = T ( a e ). $    Lastly, for a linear map $ T $ between suitable modules, $ {\rm Ran} ( T ) $ will denote the Range of $ T. $

\section{Preliminaries}

In the first subsection, we recall some concepts about centered bimodules. Then we discuss some generalities of the flip map on vector spaces (and bimodules) and the space of forms $\Omega^k(\A)$ on spectral triples $(\A,\mathcal{H},D)$. Finally, we discuss a certain bilinear splitting of $\oneform \tensora \oneform$ and some of its consequences.

\subsection{Centered bimodules}

We recall the concept of centered bimodules in this subsection. Unless otherwise mentioned, we will assume that the left and right $\A$-actions  on an $ \A$-$\A $-bimodule $ \mathcal{E} $ are not the same. 

\bdfn \label{skeide}
The center of an $ \A$-$\A $-bimodule $ \mathcal{E} $ is defined to be the set
$$ \mathcal{Z} ( \mathcal{E} ) = \{ e \in \mathcal{E}: e a = a e ~ \forall ~ a ~ \in \A \}. $$
The bimodule $ \mathcal{E} $ is called centered if $ \mathcal{Z} ( \mathcal{E} ) $ is right $ \A $-total in $ \E $, i.e, the right $\A$-linear span of $ \mathcal{Z} ( \mathcal{E} ) $ equals $ \mathcal{E}. $
\edfn

It is easy to see that $ \mathcal{Z} ( \E ) $ is a $ \mathcal{Z} ( \A ) $-bimodule. Indeed, if $ e $ is an element of $ \Ecenter,~ a  $ belongs to $ \Acenter $ and $ b  $ belongs to $ \A, $ then 
$$ b ( e a ) = e  b a  = ( e a ) b. $$

For a related notion of central bimodules, we refer to the paper \cite{dubois} of Dubois-Violette and Michor. It is easy to see that a centered bimodule in the sense of Skeide ( \cite{Skd} ) is a central bimodule in the sense of Dubois-Violette and Michor ( \cite{dubois} ), i.e, if $ \E $ is a centered module in the sense of Definition \ref{skeide}, then $ e a = ae $ for all $ e $ in $ \E $ and for all $ a $ in $ \mathcal{Z} ( \A ).$ We refer to \cite{article3} for a proof of this fact.

\bxmpl If $ \A = C^\infty ( M ) $ for some compact manifold $ M, $ and $ \Gamma ( E ) $ is the $ \A$-$\A $-bimodule of sections of some
smooth vector bundle $ E $ on $ M, $ then since $ \A $ is commutative, the right $\A$-action on $ \E $ can be defined to be the left $\A$-action and so  $ \Gamma (E) $ is centered. In particular, the $ \A$-$\A $-bimodule $ \Omega^k ( M ) $ of $ k $-forms on $ M $ is centered.

  Now suppose $ \A $ is a (possibly) noncommutative algebra and $ \E $ is an $\A$-$\A$-bimodule. If $ \mathcal{E} $ is free as a right $\A$ module of the form $ \mathcal{A} \tensorc \mathbb{C}^n $ such that $a(b \tensorc v) = ab \tensorc v$ and $(a \tensorc v)b= ab \tensorc v$ for all $a,b \in \A$ and $v \in \IC^n$, then $\E$ is centered with $ \mathcal{Z} ( \mathcal{E} ) = {\rm Span}\{ a \tensorc e_i: i = 1,2, \cdot, \cdot \cdot n, ~ a \in \mathcal{Z} ( \cla ) \},$ where $\{e_i:i=1,2,...,n\}$ is the canonical basis of $\IC^n$.
	\exmpl
	
We record the following well-known facts which are known to experts (see \cite{pinzari}). We provide a proof for the sake of completeness.
\bppsn \label{28thmay2018}
\begin{enumerate}
	\item Let $S$ be a right $\A$-total subset of a right $\A$ module $\E$. If $T_1$ and $T_2$ are two right $A$ linear maps from $\E$ to another right $\A$ module $F$ such that they agree on $S$, then they agree everywhere on $\E$.
	\item Let $\E$ and $\F$ be $\A$-$\A$ bimodules which are finitely generated and projective as both left and right $\A$ modules. Then for elements $e_i \in \E$, $f \in \F$ and $\phi_i \in \F^{*}$, the map $\zeta_{\E,\F}:\E\tensora \F^{*} \rightarrow {\rm Hom}_{\A}(\F,\E)$ defined by $\zeta_{\E,\F}(\sum_i e_i\tensora \phi_i)(f)=\sum_i e_i\phi_i(f)$ defines an isomorphism of $\A$ bimodules.
\end{enumerate}
\eppsn
{\bf Proof:} Let $S$ be a right $\A$-total subset as in 1. If $e$ is an element of $\E$, there exist elements $s_i$ in $S$ and $a_i$ in $\A$ such that $e = \sum_i s_i a_i$. Then, $$T_1(e) = \sum_i T_1(s_i) a_i = \sum_i T_2 (s_i) a_i = T_2(e).$$
Now we prove the second observation. Since $e$ is finitely generated projective, there exists $n$ and an idempotent $p$ in ${\rm Hom}_\A(\A^n, \A^n)$ such that $\E = p(\A^n)$. Let $\{ \omega_1, \omega_2, ..., \omega_n \}$ be a basis of $\A^n$ so that $\E$ is generated by $\{ p(\omega_i) \}_i$. If $T$ be an element of ${\rm Hom}_\A(\F, \E)$, then there exists elements $\phi_j$ in $\F^*$ such that 
$$ T(f) = \sum_j p(\omega_j) \phi_j (f) \ {\rm for \ all} \ e. $$ Clearly, $T = \zeta_{\E , \F}(\sum_j p(\omega_j) \tensora \phi_j)$, proving that $\zeta_{\E , \F}$ is onto. For proving that $\zeta_{\E , \F}$ is one-one, we observe that it can be easily verfied that the map $\zeta_{\E , \F}$ is a restriction of $\zeta_{\A^n , \F}$ and therefore also one-one. This completes the proof. \qed

In this article, we will be dealing with a class of centered bimodules. The next proposition deals with such bimodules.

\bppsn \label{assumption3implycentered}
Suppose $\E$ is an $\A$-$\A$-bimodule such that the map	$ u^\E: \mathcal{Z} ( \E ) \otimes_{\mathcal{Z} ( \A )} \A \rightarrow \E $ defined by
		$$  u^\E \left( \sum_{i} e^{\prime}_i  \otimes_{\mathcal{Z} ( \A )} a_i \right) =\sum_{i} e^{\prime}_i a_i $$
		is an isomorphism of vector spaces. Then we have the following isomorphism of $\A$ bimodules: 
		$$ \E \cong \A \otimes_{\mathcal{Z} ( \A )} \mathcal{Z} ( \E ) \cong \mathcal{Z} ( \E ) \otimes_{\mathcal{Z} ( \A )}\A. $$
		In particular, the right $ \A $-linear span of $ \mathcal{Z} ( \E ) $ is total in $ \E, $ i.e, $ \E $ is centered.
\eppsn
{\bf Proof:} It is clear from the definitions that the map $ u^\E $ is left $\mathcal{Z} ( \A )$, right $ \cla $-linear. Let us
	define a left $\A$, right $\mathcal{Z} ( \A )$-linear map
	$ v^\E: \A \otimes_{\mathcal{Z} ( \A )} \mathcal{Z} ( \E ) \rightarrow \E $ by
	 $$ v^\E \left(
	\sum_i a_i \otimes_{\mathcal{Z} ( \A )} e^{\prime}_i \right) = \sum_i a_i e^{\prime}_i.$$
	Consider the map
	$p: \mathcal{Z} ( \E ) \times \A \raro \A \otimes_{\mathcal{Z} ( \A )} \mathcal{Z} ( \E )$
	 given by $(e,a) \mapsto (a \otimes_{\mathcal{Z} ( \A )}e)$. It is clear that $p(ea^\prime, a)=p(e,a^\prime a)$, so that
	 we get a well-defined map $$\overline{p}: \mathcal{Z} ( \E ) \otimes_{\mathcal{Z} ( \A )} \A \raro \A \otimes_{\mathcal{Z} ( \A )} \mathcal{Z} ( \E ),
	 {\rm \ given \ by \ } (e \otimes_{\mathcal{Z} ( \A )}a) \mapsto (a \otimes_{\mathcal{Z} ( \A )} e).$$ It is in fact an isomorphism, with the inverse map, say $q$, given by
	 $$ q(a \otimes_{\mathcal{Z} ( \A )} e) = (e \otimes_{\mathcal{Z} ( \A )} a).$$ Observe that $v^\E=u^\E \circ q$, hence $v^\E$ is an isomorphism as well. Thus, the map $v^\E$ is also a vector space isomorphism as well.
		
	Next, we endow $ \mathcal{Z} ( \E ) \otimes_{\mathcal{Z} ( \A )} \A $ with an $ \A$-$\A $ bimodule structure defined by 
		$$ b ( e \otimes_{\mathcal{Z} ( \A )} a ) c = e \otimes_{\mathcal{Z} ( \A )} b a c, $$
		where $ e \in \mathcal{Z} ( \E ), a \in \A, b,c \in \A.$	 Then	it is easy to see that $ u^\E $ defines an $ \A$-$\A $ bimodule isomorphism. The other isomorphism follows by using the map $ v^\E. $	\qed		
	
The following theorem is of crucial importance in the sequel.

\bthm \label{skeide2}
({\rm Theorem 6.10, \cite{Skd}}) Let $\E$ be an $\A$-$\A$ bimodule which is centered. Then there exists a unique $\A$-$\A$ bimodule isomorphism $\sigma^{{\rm can}}: \E \tensora \E \rightarrow \E \tensora \E$ such that $\sigma^{{\rm can}}(\omega \tensora \eta)=\eta \tensora \omega$ for all $\omega,\eta \in \Ecenter$. Moreover, $(\sigma^{{\rm can}})^2= {\rm id}$ so that
$P_{{\rm sym}}^{{\rm can}}:=\frac{1}{2}(1 + \sigma^{{\rm can}}): \E \tensora \E \rightarrow \E \tensora \E$ is an $\A$-$\A$ bilinear idempotent map.
\ethm	

{\bf Proof:} We only need to remark that the equation $(\sigma^{\rm can})^2={\rm id}$ is derived in the proof of {\rm Theorem 6.10, \cite{Skd}}. Indeed, since $\E$ is centered, $\E \tensora \E = {\rm Span}\{e\tensora f a: e,f \in \Ecenter, a \in \A\}$ and $(\sigma^{\rm can})^2(e \tensora f a)=\sigma^{\rm can}(\sigma^{\rm can}(e \tensora f a))= \sigma^{\rm can}(f \tensora e a)=e \tensora f a$. \qed

 Let us make the following observation:
\blmma
\label{sigmaoncenteroneside}
For a centered $\A$-$\A$ bimodule $\E$, we have $ \sigma^{\rm can} ( \omega \tensora e ) = e \tensora \omega $ and $ \sigma^{\rm can} ( e \tensora \omega ) = \omega \tensora e $ for all
$ \omega \in \mathcal{Z} ( \mathcal{E} ) $ and $ e \in \mathcal{E}. $
\elmma

{\bf Proof:} Since $ \mathcal{E} $ is centered and $ \sigma $ is right $\A$-linear, it is enough to prove the lemma for elements $ e  $ of the
form $ \eta b $ where $ \eta \in \mathcal{Z} ( \mathcal{E} ) $ and $ b \in \A. $

We compute $ \sigma^{\rm can} ( \omega \tensora \eta b ) = \sigma^{\rm can} ( \omega \tensora \eta ) b = ( \eta \tensora \omega ) b = \eta \tensora \omega b =
\eta \tensora b \omega = \eta b \tensora \omega = e \tensora \omega. $

The other equality follows similarly. \qed
 
 We will end this subsection with Lemma \ref{module_symm_iso}. But before that, we want to state and prove Proposition \ref{braidingargument} whose proof is basically a reformulation of the proof of the existence and uniqueness of Levi-Civita connections for pseudo-Riemannian manifolds.

	Let $ V $ be a complex vector space and $ \sigma^{\mathbb{C}} $ denotes the map from $ V \tensorc V \rightarrow V \tensorc V $ defined on simple tensors by the formula $ \sigma^{\mathbb{C}} ( v \tensorc w ) = w \tensorc v. $ We will use the maps $ \sigma^{\mathbb{C}}_{12}:= \sigma^{\mathbb{C}} \tensorc {\rm id}_V, ~ \sigma^{\mathbb{C}}_{23}: = {\rm id}_V \tensorc \sigma^{\mathbb{C}} $ and $ \sigma^{\mathbb{C}}_{13}:= \sigma^{\mathbb{C}}_{12} \sigma^{\mathbb{C}}_{23} \sigma^{\mathbb{C}}_{12} .$
		
	Then the map $ P^{\mathbb{C}}: = \frac{\sigma^{\mathbb{C}} + 1}{2} $ is an idempotent. We will denote $ {\rm Ran} ( P^{\mathbb{C}} ) $ by $ V \otimes^{{\rm sym}}_{\mathbb{C}} V. $ We will need the maps $ P^{\mathbb{C}}_{12}:= P^{\mathbb{C}} \tensorc{\rm id}_V $ and $ P^{\mathbb{C}}_{23}:= {\rm id}_V \tensorc P^{\mathbb{C}} .$
	Thus, for elements $ v_1, v_2, v_3 $ in $ V, ~ P^{\IC}_{12} ( v_1 \tensorc v_2 \tensorc v_3 ) = \frac{1}{2} ( v_1 \tensorc v_2 + v_2 \tensorc v_1 ) \tensorc v_3 $ and $ P^{\IC}_{23} ( v_1 \tensorc v_2 \tensorc v_3 ) = v_1 \tensorc \frac{1}{2} ( v_2 \tensorc v_3 + v_3 \tensorc v_2 ). $
	
	\bppsn
  \label{braidingargument}
		If $ V $ is a vector space, then each of the following maps is an isomorphism of vector spaces.
		$$ P^{\mathbb{C}}_{12}|_{{\rm Ran} ( P^{\mathbb{C}}_{23} ) }: {\rm Ran} ( P^{\mathbb{C}}_{23} ) = V \tensorc ( V  \otimes^{{\rm sym}}_{\mathbb{C}} V  ) \rightarrow {\rm Ran} ( P^{\mathbb{C}}_{12} ) = ( V \otimes^{{\rm sym}}_{\mathbb{C}} V ) \tensorc V $$
		$$ P^{\mathbb{C}}_{23}|_{{\rm Ran} ( P^{\mathbb{C}}_{12} ) }: {\rm Ran} ( P^{\mathbb{C}}_{12} ) =
	(	V \tensorc^{\rm sym} V )  \otimes_{\mathbb{C}} V  \rightarrow {\rm Ran} ( P^{\mathbb{C}}_{23} ) =  V \tensorc ( V \otimes^{{\rm sym}}_{\mathbb{C}} V ) $$
		\eppsn

	{\bf Proof:} We prove the statement about the first of the two maps since the proof for the other map is similar.
	Let us begin by proving that the first map is one-one. Let $ X \in {\rm Ran} ( P^{\mathbb{C}}_{23} ) $ such that $P^{\mathbb{C}}_{12}(X)=0$. That is, $\sigma^{\mathbb{C}}_{23}(X)=X$
 and $\sigma^{\mathbb{C}}_{12}(X)=-X.$ Now, it is easy to verify the following braid relations:
 \begin{equation} \label{braid_eq}
 \sigma^{\mathbb{C}}_{12}\sigma^{\mathbb{C}}_{23}\sigma^{\mathbb{C}}_{12}=\sigma^{\mathbb{C}}_{23}\sigma^{\mathbb{C}}_{12}\sigma^{\mathbb{C}}_{23}.
 \end{equation}
But we have $\sigma^{\mathbb{C}}_{12}\sigma^{\mathbb{C}}_{23}\sigma^{\mathbb{C}}_{12}(X) = - \sigma^{\mathbb{C}}_{12}\sigma^{\mathbb{C}}_{23}(X)=
-\sigma^{\mathbb{C}}_{12}(X)= X$. On the other hand, $ \sigma^{\mathbb{C}}_{23}\sigma^{\mathbb{C}}_{12}\sigma^{\mathbb{C}}_{23}(X)= \sigma^{\mathbb{C}}_{23}\sigma^{\mathbb{C}}_{12}(X)
=-\sigma^{\mathbb{C}}_{23}(X)= - X.$ This implies, $X=-X, $ i.e. $X=0$. Thus, the map $ P^{\mathbb{C}}_{12}|_{{\rm Ran} ( P^{\mathbb{C}}_{23} ) } $ is injective.

			Now we come to surjectivity.
			 If $ V $ is finite dimensional, surjectivity follows since $ {\rm Ran} ( P^{\mathbb{C}}_{23} ) $ and $ {\rm Ran} ( P^{\mathbb{C}}_{12} ) $ are of the same dimension. In the general case, given any $ \xi \in ( V \otimes^{{\rm sym}}_{\mathbb{C}} V ) \tensorc V $ such that $ \sigma^{\mathbb{C}}_{23} ( \xi ) = \xi, $ there exists a natural number $ n $ and linearly independent elements $ e_1, e_2,..., e_n  $ of $ V $ such that $ \xi $ belongs to $ ( K \otimes^{{\rm sym}}_{\mathbb{C}} K ) \tensorc K,$ where $ K:= {\rm span} \{ e_1, e_2,..., e_n  \} .$ If $ P^{\mathbb{C}}_{K,12} $ denotes the map $ P^{\mathbb{C}}_{12}|_{K \tensorc K \tensorc K }, $ then by the surjectivity of $ P^{\mathbb{C}}_{12}|_{{\rm Ran}(P^{\mathbb{C}}_{23} )} $ for finite dimensional vector spaces, there exists $ \eta \in K \tensorc ( K  \otimes^{{\rm sym}}_{\mathbb{C}} K  ) $ such that $ P^{\mathbb{C}}_{K, 12} ( \eta ) = \xi. $ Since $\xi $ is arbitrary, the proof of surjectivity is complete. 	\qed

\blmma
 \label{module_symm_iso}
Let $\E$ be a centered $\A$-$\A$ bimodule and $\sigma^{\rm can}: \E \tensora \E \to \E \tensora \E$ be as in Theorem \ref{skeide2}. Define  $P^{\rm can}_{ij}:= \frac{1}{2}(1+\sigma^{\rm can}_{ij}): \E \tensora \E \tensora \E \to \E \tensora \E \tensora \E$, $(i,j)=(12), (13), (23)$. Then the following maps are bimodule isomorphisms:
$$ P^{\rm can}_{12}|_{{\rm Ran} ( P^{\rm can}_{23} ) }: {\rm Ran} ( P^{\rm can}_{23} )
\rightarrow {\rm Ran} ( P^{\rm can}_{12} ),~~~~
	 P^{\rm can}_{23}|_{{\rm Ran} ( P^{\rm can}_{12} ) }: {\rm Ran} ( P^{\rm can}_{12} )
		 \rightarrow {\rm Ran} ( P^{\rm can}_{23} ). $$
		\elmma

{\bf Proof :}We begin by noting that since $ \sigma^{\rm can} $ is a bimodule map, the maps $ \sigma^{\rm can}_{12}, \sigma^{\rm can}_{23}, \sigma^{\rm can}_{13} :
\E \tensora \E \tensora \E \rightarrow \E \tensora \E \tensora \E $ defined as $ \sigma^{\rm can} \tensora {\rm id}_\E, ~ {\rm id}_\E \tensora \sigma^{\rm can}, ~ \sigma^{\rm can}_{12} \sigma^{\rm can}_{23} \sigma^{\rm can}_{12} $ respectively are well defined bimodule morphisms, the proof of injectivity can be given by verbatim adaptation of the arguments of Proposition \ref{braidingargument}, as the braid relations (\ref{braid_eq}) do hold for
 the maps $\sigma^{\rm can}_{ij}$ as well. For surjectivity, we also use Proposition \ref{braidingargument}. Indeed, the ${\mathbb{C}}$-vector space
 ${\mathcal W}:={\mathcal Z}({\mathcal E}) \otimes_{\mathbb{C}} {\mathcal Z}({\mathcal E}) \otimes_{\mathbb{C}} {\mathcal Z}({\mathcal E})$ is right ${\mathcal A}$-total in ${\mathcal E} \otimes_{\mathcal A} {\mathcal E} \otimes_{\mathcal A} {\mathcal E}$. As each of the maps $P^{\rm can}_{ij}$ is 
 right ${\mathcal A}$-linear and leaves ${\mathcal W}$ invariant, ${\rm Ran}(P^{\rm can}_{ij})$ is the right ${\mathcal A}$-linear span of
 the image of ${\mathcal W}$ under $P^{\rm can}_{ij}$ for all $i,j.$ Now, the restriction of $P^{\rm can}_{ij}$ to ${\mathcal W}$ is nothing but the map $P^{\mathbb{C}}_{ij}$ as in the statement of Proposition
  \ref{braidingargument}, with $V={\mathcal Z}({\mathcal E})$ as vector spaces. By Proposition \ref{braidingargument}, we get $P^{\rm can}_{12}P^{\rm can}_{23}({\mathcal W})=P^{\rm can}_{12}({\mathcal W})$.
  Taking the right ${\mathcal A}$ linear spans on both sides, we get the surjectivity of $P^{\rm can}_{12}|_{{\rm Ran}(P^{\rm can}_{23})}.$ The surjectivity of the other map follows in a similar way. \qed
			
			\subsection{Generalities on spectral triples} \label{sptripleforms}
		In this subsection, we quickly recall the definitions of spectral triple (see \cite{connes}) and the  construction of the space of forms on it.\\
	A spectral triple on a unital *-algebra $\A$ is given by $(\A,\mathcal{H},D)$ where $\mathcal{H}$ is a separable Hilbert space admitting a representation of $\A$ and $D$ is a (possibly unbounded) self adjoint operator on $\mathcal{H}$ such that $[D,a]$ has bounded expansion for all $a \in \A$.
	
	Let $ ( \mathcal{A}, \mathcal{H}, D ) $ be a spectral triple.  
Let $ d$ denote the canonical derivation on $ \A $ defined by 
$$ d ( \cdot ) = \sqrt{- 1} [ D, \cdot ]. $$
We fix the spectral triple and henceforth denote $d$ simply by $d$. We refer to \cite{connes}, \cite{frolich}, \cite{landi} for definition and detailed discussion on the bimodule
 of noncommutative differential forms. However, we only need to consider the spaces of one and two-forms, to be denoted by $\oneform$ and $\twoform$ respectively. They are defined as follows.
 The space $\oneform$ is the linear span of elements of the form $a [D, b]$, $a,b \in \A$, in ${\mathcal B}({\mathcal H})$. It is clearly a right $\A$ module.
 As any elements $\omega, \eta$ are elements of ${\mathcal B}({\mathcal H})$, we have a natural multiplication map $m_0: \oneform \tensora \oneform \ni
 (\omega \tensora \eta) \mapsto \omega \eta \in {\mathcal B}({\mathcal H}).$
 Let $\clj$ be the module of junk forms, defined to be the
 right $\A$-submodule of the ${\rm Im}(m_0)$ spanned by elements of the form $\sum_i [D,a_i] [D,b_i]$ (finite sum) such that $\sum_i a_i [D,b_i]=0$ $(a_i, b_i \in \A$).
 We define $\twoform={\rm Im}(m_0)/\clj$ and let $ \wedge: \oneform \tensora \oneform \rightarrow \twoform $ be the composition of $m_0$ and the quotient map from  ${\rm Im}(m_0)$ to $\twoform$. From now on, we adopt the notation $\omega \wedge \eta := \wedge(\omega \tensora \eta)$.

Now we are in a position to recall the extension of the map $d$ to higher order forms. For an element $T$ in ${\rm Im}(m_0)$, let us denote the equivalence class of $T$ in $\twoform$ by the symbol $[T]$. Then it can be checked that the map which send $[D,a]b$ to $[[D,a][D,b]]$ extends to a well-defined map (to be denoted again by d, by abuse of notation) from $\oneform$ to $\twoform$. Moreover, for an element $\omega$  in $\oneform$ and $a$ in $\A$, we have $d(wa) = (d\omega)a - \omega \wedge da$. The map $d$ can be extended to  higher order forms of the spectral triple is a similar way. We refer to Chapter 7 of \cite{landi} for the details.
	
\blmma
The map $\wedge: \oneform \tensora \oneform \rightarrow \twoform$ is an $\A$-$\A$ bilinear map. \label{misbilinear}
\elmma	
{\bf Proof:} Clearly, the map $m_0$ is a bilinear map, so we have to prove that the quotient map from ${\rm Im}(m_0)$ to $\twoform$ is bilinear. Again, for this it is enought to prove that $\mathcal{J}$ is closed under left $\A$ module multiplication since by definition, $\mathcal{J}$ is closed under right $\A$ multiplication. To this end, let $\sum_i a_i[D,b_i]=0$ be a finite sum, where $a_i,b_i \in \A$. For $c \in \A$, $c\sum_i[D,a_i][D,b_i]=\sum_i[D,ca_i][D,b_i]-[D,c]\sum_ia_i[D,b_i]=\sum_i[D,ca_i][D,b_i]$, as $\sum_ia_i[D,b_i]=0$. But since $\sum_ica_i[D,b_i]=0$, we have $\sum_i[D,ca_i][D,b_i] \in \mathcal{J}$ which implies that $c\sum_i[D,a_i][D,b_i] \in \mathcal{J}$. This proves the claim.
 \qed

	\bdfn
Let $ ( \A, \mathcal{H}, D ) $ be a spectral triple. A (right) connection on $ \oneform $ is a $ \IC $-linear map $ \nabla: \oneform \rightarrow \oneform \tensora \oneform $ such that $ \nabla ( \omega a ) = \nabla ( \omega ) a + \omega \tensora da $ for all $ \omega $ in $ \oneform $ and $ a $ in $\A$.
The torsion of a connection $\nabla$ is defined to be the map $T_{\nabla}:=\wedge \circ \nabla + d : \oneform \to \twoform. $
\edfn

Let us recall ( \cite{connes} ) that a spectral triple $ ( \mathcal{A}, \mathcal{H}, D ) $  is said to be of compact type if the operator $ ( 1 + D^2 )^{-1} $ is compact. 	However, in this subsection, we have not used this assumption  although we will need it in Subsection \ref{assumptionsonsptriple} for constructing the candidate of the pseudo-Riemannian bilinear metric for a spectral triple.	

	\subsection{A natural splitting of $ \oneform \tensora \oneform $} \label{splittingsubsection}
		
	For a Riemannian manifold $M$ with $\oneformclassical$ as the space of one-forms, we have the following decomposition of $C^\infty(M)$ bimodules:
	$$\oneformclassical \otimes_{C^\infty(M)} \oneformclassical = {\rm Sym}^2 (\oneformclassical)\oplus \F,$$	
	where ${\rm Sym}^2(\oneformclassical):={\rm Span}(\omega \otimes_{C^\infty(M)}\eta + \eta \otimes_{C^\infty(M)}\omega: \omega,\eta \in \oneformclassical)$ is the space of all symmetric 2-tensors and $\F$ is isomorphic to $\Omega^2(M)$. If $\wedge:\oneformclassical\otimes_{C^\infty(M)}\oneformclassical \rightarrow \Omega^2(M)$ denotes the canonical map, then ${\rm Sym}^2(\oneformclassical)$ is nothing but ${\rm Ker}(\wedge)$. The aim of this subsection is to discuss an analogous decomposition of $\oneform\tensora\oneform$ and some of its consequences.
 
  \bdfn \label{defn24thmay2018}
  We say that a spectral triple $(\A,\mathcal{H},D)$ is quasi-tame if the following conditions hold:
  \begin{itemize}
  	\item[i.] The bimodule $\oneform$ is finitely generated and projective as a right $\A$ module.
  	\item[ii.] There exists a right $\A$-module $\F$ such that:
  \begin{equation}
  \oneform \tensora \oneform = {\rm Ker}(\wedge) \oplus \F \label{splitting25thmay2018}
  \end{equation}
    \item[iii.] The idempotent $P_{sym} \in {\rm Hom}_{\A}(\oneform \tensora \oneform, \oneform \tensora \oneform)$ mapping onto ${\rm Ker}(\wedge)$ and with kernel $\F$ is an $\A$-$\A$ bimodule map.
\end{itemize}
Then we will denote ${\rm Ker}(\wedge)$ by the symbol $\oneform \tensorsym \oneform$. \\ Moreover, $\sigma$ will denote the map $2P_{sym}-1$. 
	\edfn

\blmma \label{28thmay20182}
Let $(\A,\mathcal{H},D)$ be a quasi-tame spectral triple. Then we have the following
\begin{enumerate}
\item[i.] $ \oneform \tensorsym \oneform =: {\rm Ker} ( \wedge ) $ and $ {\rm Im} ( \wedge ) $ are $\A$ bimodules.
\item[ii.]  $ \sigma $ is an $\A$-$\A$ bimodule map.
\item[iii.] $ P^2_{{\rm sym }} = P_{{\rm sym }} $ and $ \sigma^2 = {\rm id}$.
\end{enumerate}
\elmma

{\bf Proof:} By Lemma \ref{misbilinear} the map $\wedge$ is bilinear, hence ${\rm Ker}(\wedge)$ and ${\rm Im}(\wedge)$ are $\A$ bimodules. This gives us the first claim. The second claim, i.e the $ \A$-$\A $-bilinearity of $ \sigma $ follows from the $ \A $-$ \A$-bilinearity of $ P_{{\rm sym}}. $ The third claim follows from the fact that $\Psym$ is an idempotent.
 \qed

\vspace{4mm}
 
 Let us recall (\cite{frolich}) that a connection $\nabla$ on $\oneform$ is said to be torsion-less if $T_{\nabla}=0$. We are using the definition the same definition of torsion of a connection as used by previous authors ( \cite{frolich}, \cite{heckenberger_etal}  ) with the only difference being that we are using right connections as opposed to left connections used in \cite{frolich} and \cite{heckenberger_etal}.  We have the following result as a consequence of the assumptions made in Definition \ref{defn24thmay2018}.

	\bthm
\label{torsionless}
For a quasi-tame spectral triple, there exists a torsion-less connection on $\oneform$.
\ethm
{\bf Proof:} We have a sub-bimodule ${\mathcal F}:= {\rm Im}(1-\Psym)$ of $\oneform \tensora \oneform$ and a bimodule isomorphism, say $Q,$ from ${\mathcal F}$ to ${\rm Im}(\wedge) =\twoform$, satisfying $Q((1-P_{\rm sym})(\beta) ) = \wedge(\beta )$ for $\beta \in \oneform \tensora \oneform$. Moreover, as $\oneform$ is finitely generated and projective, 
 we can find a free rank $n$ ( right ) module $\A \tensorc {\mathbb C}^n$ containing $\oneform$ as a complemented right submodule. Let $ p $ be an idempotent in $M_n(\A) \cong {\rm Hom}_{\A}(\A \tensorc {\mathbb C}^n, \A \tensorc {\mathbb{C}}^n )$ such that $\oneform =p(\A \tensorc {\mathbb C}^n).$ 
 . 
Let $e_i, i=1, \ldots n,$ be the standard basis of $ {\mathbb C}^n$ and define $\tilde{\nabla}_0 : \A \tensorc {\mathbb C}^n \rightarrow \oneform \tensora \oneform$ by 
$$ \label{20thjan3} \tilde{\nabla}_0(e_i a):= - Q^{-1}(d(p(e_i)))a+p(e_i) \tensora da, ~~~i=1, \ldots, n, a \in \A.$$ 
\be \label{20thjan4} {\rm Then} ~ \nabla_0=\tilde{\nabla}_0|_{\oneform}\ee
defines a connection on $\oneform$ Note that, as $ \wedge \circ \Psym = 0 $ by definition, we have
	$$ \wedge \circ Q^{ - 1 } ( \wedge ( \beta ) ) = \wedge \circ Q^{ - 1 } ( Q ( 1 - \Psym ) \beta ) = \wedge ( ( 1 - \Psym ) \beta ) = \wedge (\beta) ~ \forall ~ \beta \in ~ \oneform \tensora \oneform, $$
	i,e, $ \wedge Q^{ - 1 } : {\rm Im} ( \wedge ) \rightarrow {\rm Im} ( \wedge ) $ is the identity map. Since $ d ( p ( e_i ) a ) $ belongs to the image of the map $\wedge,$ we can write
	
	\begin{eqnarray*}
	 \wedge \circ \nabla_0 ( p ( e_i ) a ) &=& - \wedge ( Q^{-1} (  d ( p ( e_i ) ) a  ) ) + p ( e_i ) \wedge d a \\
	                 &=& - d ( p ( e_i ) ) a + p ( e_i ) \wedge d a \\
																		&=& - d ( p ( e_i ) a  ).
	\end{eqnarray*}
	Therefore, $ \nabla_0 $ is a torsion-less connection on $\oneform$. \qed	
	
	\section{Pseudo-Riemannian metrics on a spectral triple} \label{4thjune2018}

Next, we want to introduce a noncommutative analogue of pseudo-Riemannian metrics. In classical differential geometry, a Riemannian metric on a manifold $M$ 
is a smooth, positive definite, symmetric 
bilinear form on the tangent (or, equivalently, co-tangent)
 bundle. One can extend it to the complexification of the tangent/cotangent spaces in two ways: either as a sesquilinear pairing (inner product) on the module of one-forms, 
 which is conjugate-linear in one variable and linear in the other,
 or, as a complex bilinear form, i.e. a $C^\infty(M)$-linear map on $\Omega^1(M) \otimes_{C^\infty(M)} \Omega^1(M)$.
 
 Somehow, the first of these two alternatives seems to be more popular in the literature so far  to formulate a noncommutative analogue of metric, with some exceptions like the formulation in \cite{heckenberger_etal} in the framework of bicovariant differential calculi on quantum groups. One advantage of defining a (Riemannian) metric for noncommutative manifold as a non-degenerate sesquilinear pairing (i.e. inner product) taking value in the underlying $C^*$ algebra is that one can use the rich and popular theory of Hilbert modules.
 
 However, when one wants to deal with pseudo-Riemannian metrics, there is no assumption of positive definiteness and the relative advantage of sesquilinear extension over the bilinear extension no longer exists. Moreover, the existence and uniqueness of classical Levi-Civita connection for a classical manifold do not need any positive definiteness and hold for an arbitrary pseudo-Riemannian metric. For this reason, it makes sense to consider bilinear non-degenerate forms as pseudo-Riemannian metrics in the noncommutative set-up. In classical case there is no difference between right module maps or bimodule maps, as the left and right $C^\infty(M)$-actions on the module of forms coincide. This is no longer true in the noncommutative framework. In fact, as we will see, requiring a pseudo-metric to be a bimodule map restricts the choice of metrics. It is reasonable to require one-sided (right/left) $\A$-linearity only. For this reason, we give the following definition  which is actually valid for a class of bimodules ( see \cite{article2} ) and not only spectral triples:

\bdfn \label{11thjuly20183}
Let $(\A,\mathcal{H},D)$ be a spectral triple and let $\E:= \oneform$ be the $\A$-$\A$ bimodule of one-forms. Let us assume that $\E$ satisfies the conditions of Definition \ref{defn24thmay2018} and let $\sigma$ be as in that definition. A pseudo-Riemannian metric $ g $ on $ \E$ is an element of $ {\rm Hom}_{\A} ( \E \tensora \E, \A ) $ such that

(i) $g$ is symmetric, i.e. $ g \circ \sigma = g, $

 (ii) $g$ is non-degenerate, i.e, the right $ \A$-linear map $ V_g: \E \rightarrow {\E}^* $ defined by $ V_g ( \omega ) ( \eta ) = g ( \omega \otimes_\A \eta ) $ is
 an isomorphism of right $ \A$ modules.

We will say that a pseudo-Riemannian metric $g$ is a pseudo-Riemannian bilinear metric if $ g $ is an $ \A $-$\A$ bimodule map. It is called a Riemannian metric if for all $\omega_1, \omega_2, \dots, \omega_n \in \E$, the matrix $((g(\omega_i^*\omega_j)))_{i,j}$ is a positive element of $M_n(\A)$ for all $n.$
\edfn

As an immediate consequence of the definition, we have the following important proposition.

\bppsn \label{17thoct20191}
Suppose $ g $ is a pseudo-Riemannian bilinear metric on the space of one-forms $\E:= \Omega^1_D ( \A )$ of a quasi-tame spectral triple. Then 
  \begin{equation} \label{eq1} g ( \omega \tensora \eta ) \in \Acenter ~ {\rm if} ~ {\rm both} ~  \omega ~ {\rm and} ~  \eta  ~ {\rm belong} ~ {\rm to} ~  \Ecenter.  \end{equation}
 In particular, if $ \E $ is a free right $ \A $-module of rank $n$  admitting a central basis $ \{ \omega_i \}_i \subseteq \Ecenter, $ then the components of the metric $ g_{ij}:= g ( \omega_i \tensora \omega_j ) $ belong to $\Acenter.$  
\eppsn
{\bf Proof:} The proof is a trivial consequence of the fact that  $ g $ is an $\A$-$\A$ bimodule map. Indeed, since $ \omega, \eta $ are in $ \mathcal{Z} ( \E ),$
	$$ g ( \omega \tensora \eta ) a = g ( \omega \tensora \eta a ) = g ( \omega \tensora a \eta ) = g ( \omega a \tensora \eta ) = a g ( \omega \tensora \eta ). $$ 
\qed

However, if $\A$ is noncommutative, the metric need not take values in the center on the whole of $\E \tensora \E.$ For example, if $ \omega, \eta \in \Ecenter $ and $a \in \A,$ then $ g ( \omega \tensora \eta a ) $ typically does not belong to $\Acenter$ unless $ a \in \Acenter. $ Pseudo-Riemannian bilinear metrics have been also used by Beggs and Majid. We refer to \cite{beggsmajidbook} for the details.

\brmrk
Our definition of nondegeneracy of $ g $ is stronger than the definition given by most authors who require only the injectivity of $ V_g. $ However, in the classical situation, i.e, when $ \A = C^\infty ( M ), $ these two definitions are equivalent as $ V_g $ is a bundle map from $ T^* M $ to $ ( T^* M )^* \cong TM $ in that case and the fibers are finite dimensional. 
\ermrk

To compare our definition of a pseudo-Riemannian metric with that of \cite{frolich}, \cite{Rosenberg} and \cite{pseudo}, let us consider the case when $\E$ is free (of rank $ n $) as a right $ \A $ module, i.e, $ \E $ is isomorphic to $ \IC^n \tensorc \A $ as a right $ \A $ module. Let $e_i, i=1, \ldots, n$ be the standard basis of ${\mathbb C}^n$. A pseudo-Riemannian metric in our sense is determined by an invertible element
 $A:=(( g_{ij} ))$ of $M_n(\A)$, where $ g_{ij} = g ( ( e_i \tensorc 1 ) \tensora ( e_j \tensorc 1 ) ) $ and $ g ( ( e_i \tensorc a ) \tensora ( e_j \tensorc b  ) ) = g_{ij} a b $ for all $ a, b \in \A. $  On the other hand, a pseudo-metric in the sense of \cite{pseudo}
  corresponding to $A$ will be given by the sesquilinear pairing $ \left\langle \left\langle  e_i \otimes_{\IC} a, e_j \otimes_{\IC} b \right\rangle\right\rangle =a^*g_{ij}b$. Thus, there is a one-to-one correspondence between these
  two notions of pseudo-metric at least for the case when $\E$ is free as a right $ \A $ module. In fact, they do agree in a sense on the basis elements.
  But their extensions are quite different as maps.

	Throughout this subsection, we will assume that $ ( \A, \mathcal{H}, D ) $ is a spectral triple satisfying the conditions of Definition \ref{defn24thmay2018} so that we can freely use the notation $ \sigma $ introduced in that definition and the results in Lemma \ref{28thmay20182}. 
	
	The rest of this subsection is devoted to proving some results on pseudo-Riemannian metrics on $\E:=\oneform$ which will be used in Subsection \ref{subsectionmainproof}. In the next subsection, we will discuss a candidate for a canonical pseudo-Riemannian bilinear metric on $\E$.
	
	\bppsn \label{28thmay20183}
	Let $ ( \A, \mathcal{H}, D ) $ be a quasi-tame spectral triple such that $\E$ is centered as an $\A$-$\A$-bimodule and $\sigma = \sigma^{{\rm can}}, $ where $ \sigma^{{\rm can}} $ is as in Theorem \ref{skeide2}. If $g$ is a pseudo-Riemannian metric on $\E,$ then we have 
 $$	g ( \omega \tensora \eta ) = g ( \eta \tensora \omega ) $$
if either $\omega$ or $ \eta $ belongs to $ \Ecenter. $
	\eppsn
	{\bf Proof:} Let $ \omega \in \Ecenter $ and $ \eta \in \E. $ 
		Then by by Lemma \ref{sigmaoncenteroneside}, we compute
		$$ g ( \omega \tensora \eta ) = g \circ \sigma ( \omega \tensora \eta ) = g (\sigma^{{\rm can}} ( \omega \tensora \eta )) = g ( \eta \tensora \omega ).$$
		\qed

\blmma \label{lemma4}
 Let $ g $ be a pseudo-Riemannian metric on $\E$. Then for all $ T \in {\rm Hom}_\A ( \E, \E ) $ which is also left $\A$-linear, there exists a unique element $ T^* \in {\rm Hom}_\A ( \E, \E ) $ such that for all $ e, f$ in $\E, $
 $$ g ( T^* ( e ) \tensora f ) = g ( e \tensora T ( f ) ). $$
\elmma
 {\bf Proof:} Suppose $e \in \E$. We define an element $z(e) \in \E$ by the equation $V_g(z(e))(f) = g(e \tensora T(f))$. The above definition is well-defined by the non-degeneracy of $V_g$. Clearly, the element $z(e)$ is the unique choice for $T^{*}(e)$.\\
 For proving that the map $e \mapsto T^{*}(e):= z(e)$ is right $\A$-linear, we compute 
 $$V_g(T^{*}(ea))(f) = V_g(ea)(T(f)) = g(ea \tensora T(f))$$
 $$= g(e \tensora aT(f))=g(e \tensora T(af)) \text{ [Since $T$ is left $\A$-linear]}$$ $$=V_g(T^{*}(e))(af)=V_g(T^{*}(e)a)(f)$$
  Since $V_g$ is an isomorphism, we have $T^{*}(ea)=T^{*}(e)a$.
	  \qed

\blmma \label{19thfeb2018}
 Suppose $\E:=\oneform$ is centered. Then for all $a$ in $\Acenter$ and $\omega$ in $\E$, \label{19thfeb20182} $a \omega = \omega a$.
Moreover, if $ g $ is a pseudo-Riemannian bilinear metric on $\E$ and $ \omega, \eta \in \Ecenter, $ then 
$$ g ( \omega \tensora \eta ) e = e g ( \omega \tensora \eta ) ~ {\rm for} ~ {\rm all} ~ e ~ \in ~ \E. $$ 
\elmma
 {\bf Proof:} Let $e $ be an element of $\E.$ Since $\E$ is centered, we can write $ e = \sum_i e_i a_i $ for some elements $ e_i \in \Ecenter $ and $ a_i $ in $ \A. $ Then for all $ a $ in $ \Acenter, $ we have
	$$ a e = \sum_i a e^\prime_i a_i = \sum_i e^\prime_i a a_i = \sum_i e^\prime_i a_i a = ea. $$
This proves the first assertion. The second assertion then follows from a combination of the first assertion and Proposition \ref{17thoct20191}.
 \qed

\bdfn
Suppose $g$ is a pseudo-Riemannian bilinear metric on $\E$. We define $g^{(2)}:(\E\tensora \E)\tensora (\E \tensora \E) \rightarrow \A$ by
$$g^{(2)}((e\tensora f)\tensora (e'\tensora f')) = g(e \tensora g(f \tensora e') f')$$
\edfn

	We spell out the relationship between $g^{(2)}$ and the inner product on the internal tensor product of Hilbert modules. Suppose $(\A, \mathcal{H}, D)$ is a spectral triple, $\E$ the bimodule of one-forms. We will need  to make explicit use of the ${}^*$-structure on $\A$ and $\E:=\oneform$ inherited from $ B ( \mathcal{H} ).$  
	Let us recall the conjugate bimodule $ \overline{\E} $ ( see \cite{pinzari}, \cite{majid_2} and references therein ) which is equal to $ \E $ as a set but with the $\A$-bimodule structures defined by the following equations:
	$$ a \overline{e} = \overline{e a^*}, ~ \overline{e} a = \overline{a^* e}. $$
	Here, $ \overline{e} $ is an element of $ \E $ viewed in $ \overline{\E}. $
	
 We have a well-defined map $S: \E \tensora \E \to \E \tensora \E$, defined by 
	\[ S(e \tensora f) = \overline{f} \tensora \overline{e} .\]
	Now suppose $g$ is a pseudo-Riemannian bilinear metric on $\E$.  
	Then the following map makes $\E$ into a right $\A$-pre-Hilbert module:
	\[ \langle \langle e, f \rangle \rangle_g = g(\overline{e} \tensora {f}). \]
	On the right hand side of this equation, we have used the obvious identification between $ \E $ and $ \overline{\E}. $  
	
	Consequently, the $\A$-valued inner product on the internal tensor product $\E \tensora \E$ is given by
	\[\langle \langle e \tensora f, e^\prime \tensora f^\prime \rangle \rangle_{g^{(2)}} = \langle \langle f, \langle \langle e, e^\prime \rangle \rangle_g f^\prime \rangle \rangle_g.
\]
We refer to \cite{pinzari} for the details. We claim that $\langle \langle e \tensora f, e^\prime \tensora f^\prime \rangle \rangle_{g^{(2)}} = g^{(2)}(S(e \tensora f) \tensora (e^\prime \tensora f^\prime))$. Indeed,
\begin{eqnarray*}
\langle \langle e \tensora f, e^\prime \tensora f^\prime \rangle \rangle_{g^{(2)}} &=& \langle \langle f, g( \overline{e} \tensora e^\prime) f^\prime \rangle \rangle_g\\
&=& g( \overline{f} \tensora g( \overline{e} \tensora e^\prime) f^\prime)\\
&=& g^{(2)}((\overline{f} \tensora \overline{e}) \tensora (e^\prime \tensora f^\prime))\\
&=& g^{(2)}(S(e \tensora f) \tensora (e^\prime \tensora f^\prime)).
\end{eqnarray*}

We end this subsection by showing that the map $ g^{(2)} $ is nondegenerate in a suitable sense.

\bppsn \label{vg2nondegenerate}
Suppose $\E$ is centered as an $\A$-$\A$-bimodule and also that $\E$ is finitely generated and projective as a right $\A$ module. Let $ g $ be a pseudo-Riemannian bilinear metric on $\E.$ Then the map $ \Vgtwo : \E \tensora \E \rightarrow ( \E \tensora \E )^* $ defined by
$$V_{g^{(2)}}(e\tensora f)(e'\tensora f')=g^{(2)}((e \tensora f) \tensora (e' \tensora f'))$$
 is an isomorphism of right $ \A $ modules. Moreover, $ g^{(2)} $ is a pseudo-Riemannian bilinear metric on $\E \tensora \E.$
\eppsn
{\bf Proof:} Let us start by proving that the map $ \Vgtwo $ is onto. Since $ \E $ is a finitely generated projective module over $\A,$ we can use the isomorphism of $ ( \E \tensora \E )^* $ with $ \E^* \tensora \E^*. $Thus, it is enough to show that $ V_g ( e ) \tensora V_g ( f ) $ belongs to the range of $ \Vgtwo $ for arbitrary elements $ e,f $ of $ \Ecenter $. Indeed, if $x_{ij} \in \E\tensora\E$ is such that $V_{g^{(2)}}(x_{ij})=V_g(e_i)\tensora V_g(e_j)$, then for $\omega = \sum e_ia_i$, $\eta=\sum f_jb_j$, where $e_i, f_j \in \Ecenter$ and $a_i, b_j \in \A$, we have
  \begin{eqnarray*}
	  V_g ( \omega ) \tensora V_g ( \eta ) &=& \sum_{i,j} V_g ( e_i ) a_i \tensora V_g ( f_j ) b_j\\
		                   &=& \sum_{i,j} V_g ( e_i ) \tensora V_g ( a_i f_j ) b_j\\
																					&=& \sum_{i,j} V_g ( e_i ) \tensora V_g ( f_j ) a_i b_j\\
																					&=& \sum_{i,j} \Vgtwo ( x_{ij} ) a_i b_j\\
																					&=& \sum_{i,j} \Vgtwo ( x_{ij} a_i b_j ).	 
	 \end{eqnarray*}
	Now, for $e,f$ in $ \Ecenter $ and $\omega,\eta$ in $\E,$ we compute
	 \begin{eqnarray*}
	 \Vgtwo ( f \tensora e ) ( \omega \tensora \eta ) &=& g^{(2)} ( ( f \tensora e ) \tensora ( \omega \tensora \eta ) )\\
		                           &=& g ( f \tensora g ( e \tensora \omega ) \eta )\\
																												 &=& g ( g ( e \tensora \omega ) f \tensora \eta  )\\
																												 &=& g ( e \tensora \omega ) g ( f \tensora \eta )\\
																												 &=& (V_g ( e ) \tensora V_g ( f ) ) ( \omega \tensora \eta ).	
	 \end{eqnarray*}
	$${\rm Hence,} ~ {\rm we} ~ {\rm have} ~ V_g ( e ) \tensora V_g ( f ) = \Vgtwo ( f \tensora e ).$$
	For proving that $\Vgtwo$ is one-one, let us suppose that for $ i = 1,2, \cdots n, $ there exist $ \omega_i, \eta_i $ in $ \E $ such that for all $ \omega^\prime, \eta^\prime$ in $\E, $ 
		$$ g^{(2)} ( ( \sum_i \omega_i \tensora \eta_i ) \tensora ( \omega^\prime \tensora \eta^\prime )  ) = 0. $$
		Then by the definition of $ g^{(2)}, $ we see that 
		$$ V_{g} ( \sum_i \omega_i g ( \eta_i \tensora \omega^\prime ) ) = 0. $$
		By nondegeneracy of $ g, $ we conclude that
		 $$ \sum_i \omega_i g ( \eta_i \tensora \omega^\prime ) = 0. $$
			Thus, if $\zeta_{\E,\E}$ is the map introduced in Proposition \ref{28thmay2018}, then we have:
			$$ \zeta_{\E,\E} ( \sum_i \omega_i \tensora \eta_i ) ( \omega^\prime ) = 0 ~ {\rm for} ~ {\rm all} ~ \omega^\prime \in \E $$
			implying that $ \sum_i \omega_i \tensora \eta_i = 0. $ 
			
			The bilinearity of $\Vgtwo$ is easy to check and hence we omit its proof. This completes the proof of the proposition.
  \qed

	\subsection{The canonical Riemannian (bilinear) metric for a spectral triple} \label{assumptionsonsptriple}
	
	Let $ ( \cla, \clh, D ) $ be a $p$-summable spectral triple (\cite{connes}) of compact type. The goal of this section is to derive some sufficient conditions for obtaining a canonical bilinear form (candidate of a pseudo-Riemannian bilinear metric) on the module $ \E := \oneform $ of one-forms.

 Consider the positive linear functional $ \tau $ on $ \clb ( \clh ) $ given by 
$$ \tau ( X ) = {\rm Lim}_\omega \frac{{\rm Tr} ( X \left| D \right|^{- p}  ) }{{\rm Tr} ( \left| D \right|^{- p} ) }, $$
where $ {\rm Lim}_{\omega} $ is the Dixmier trace as in \cite{connes} and the spectral triple is $p$-summable. We will denote the $\ast$-subalgebra generated by $\cla$ and $ [ D, \cla ] $ in $\clb ( \clh ) $ by $ \cls_0 .$ We will assume that $\tau$ is a faithful normal trace on the von Neumann algebra generated by $ \cls_0. $

 Let us recall from \cite{frolich} the construction of an $ \cla^{\prime \prime}$-valued inner product $ \left\langle \left\langle \cdot ~, ~ \cdot \right\rangle \right\rangle $ on $ \E = \oneform $ defined by the following equation:
 $$ \tau ( \left\langle \left\langle \omega, \eta \right\rangle\right\rangle a ) = \tau ( \omega^* \eta a  ) ~ \forall ~ a ~ \in \cla^{\prime \prime} ~ {\rm and} ~ \omega, \eta \in \E \subseteq \clb ( \clh ). $$
	Here, $ \omega^* $ denotes the usual adjoint of $ \omega $ in $ \clb ( \clh ). $ 
	
	As seen in Theorem 2.9 of \cite{frolich}, it can be proved that $ \left\langle \left\langle \omega, \eta \right\rangle\right\rangle $ takes values in $ \cla^{\prime \prime} \subseteq L^2 ( \cla^{\prime \prime}, \tau  ). $ We denote by $ \left\langle ~ \cdot ~, ~ \cdot ~ \right\rangle $ the positive functional $ \tau \circ \left\langle \left\langle ~ \cdot ~, ~ \cdot ~ \right\rangle \right\rangle. $
	
	Now define a natural $ \cla^{\prime \prime} $-valued bilinear form $ g $ as follows:
	
	\blmma \label{22ndjune2018}
	Let $ g: \E \tensorc \E \rightarrow \cla^{\prime \prime} $ be given by 
	$$ g ( \omega \tensorc \eta ) = \left\langle \left\langle \omega^*, \eta \right\rangle\right\rangle. $$
	Then for all $ \omega, \eta \in \E $ and $ a \in \cla, $ we have:
	$$ g ( \omega a \tensorc \eta ) = g ( \omega \tensorc a \eta ),~ g ( a \omega \tensorc \eta ) = a g ( \omega \tensorc \eta ),~ g ( \omega \tensorc \eta a ) = g ( \omega \tensorc \eta ) a.$$
	\elmma
	{\bf Proof:} The proof of the above statements are straightforward consequences of the properties of an inner product and the fact that $ ( X a )^* = a^* X^*  $ for all $ a, X \in \clb ( \clh ). $ \qed
	
	Thus, $ g $ descends to an $ \cla $-bilinear, $ \cla^{\prime \prime} $-valued map, to be denoted by $ g $ again. The restriction of $ g $ to $ \oneform \tensora \oneform $ is the candidate of a Riemannian bilinear metric in our sense, provided $ g ( \omega \tensora \eta ) \in \A $ for all $ \omega, \eta $ in $ \oneform. $
	
	Let us recall the definition of a quasi-tame spectral triple as well as the notations $P_{{\rm sym}}$ and $\sigma$ from Definition \ref{defn24thmay2018}. Then we have the following definition:
	
	\bdfn \label{defn25thmay2018}
	Let $(\A,\mathcal{H},D)$ be a quasi-tame spectral triple. Suppose the $\A$-$\A$ bilinear map $g$ as in Lemma \ref{22ndjune2018} is $\A$-valued, $V_g:\E \rightarrow \E^*$ is nondegenerate and $g\circ \sigma = g$, i.e. it gives a bilinear metric. Then we call $g$ a canonical Riemannian bilinear metric for the spectral triple $(\A,\mathcal{H},D)$.
	\edfn
	
	When $ A = C^\infty ( M ) $ for a compact Riemannian manifold $ M, $ then this construction recovers the usual Riemannian metric ( see page 128-129 of \cite{frolich} and Subsection 2.1.3  of \cite{frolich2} ). However, in the general noncommutative set-up, one usually needs additional regularity assumptions to ensure that  $ g $ takes values in $\A$ ( as opposed to $ \A^{\prime \prime} $ ). This is the content of the next proposition. 
	
	Since we have assumed that $\tau$ is faithful on the von Neumann algebra generated by $S_0$, let us oberve that the vector space $\oneform$ can be equipped with a semi-inner product defined by the equation:
	$$\langle \eta, \eta^\prime \rangle = \tau(\eta^*\eta^\prime).$$ 
	\bdfn ( \cite{frolich} ) \label{18thoct20191}
	The Hilbert space completion of $\oneform$ with respect to the inner product $ \langle \cdot, \cdot \rangle  $ is called the Hilbert space of one-forms and is denoted by the symbol $ \mathcal{H}^1_D. $  
	\edfn
	\bppsn
	Let $\mathcal{H}_D^1$ denote the Hilbert space of one-forms as in Definition \ref{18thoct20191}. Suppose that the map $$ \IR \to  \mathcal{B} (\mathcal{H}) \text{ defined by }t\mapsto e^{itD}Xe^{-itD} $$ is differentiable at $t=0$ in the norm topology of $\mathcal{B} (\mathcal{H})$, so that the map $\mathcal{L}:= -d^*d$ makes sense. If we moreover assume that $\mathcal{L}(\A) \subseteq \A $, then  $$ g ( \omega \tensora \eta ) \in \A \text{ for all } \omega, \eta \in \oneform. $$
	\eppsn
	{\bf Proof:}  Thus, $d:L^2(\A,\tau) \to \mathcal{H}_D^1$ is a densely defined linear map. By applying Lemma 3.1 of \cite{goswamiqiso}, it turns out that $d$ is closable and $\A \subseteq {\rm Dom}(\mathcal{L})$. By our assumption, we have that $\mathcal{L}$ maps $\A$ into $\A$ (in general, $\mathcal{L}$ maps $\A$ into the weak closure of $\A$ in $\mathcal{B}(L^2(\A,\tau))$.
	 We claim that 
	 $$ g ( d a \tensora d b ) = - \frac{1}{2} ( \cll ( b^* a^* ) - \cll ( b^* ) a^* - b^* \cll ( a^* ) ) ~ \forall ~ a, b \in \A, $$
	where $ d a = \sqrt{ - 1 } [ D, a ] $ as above.
	
	Indeed, for all $c$ in $\cla,$ by using the self-adjointness of $\mathcal{L}$, $\mathcal{L}(x^*)=(\mathcal{L}(x))^*$ (Lemma 3.2, \cite{goswamiqiso} and Lemma 5.1 of \cite{goswamiqiso}), we have
	 \begin{eqnarray*}
	  \tau ( \left\langle \left\langle ( d a )^*, d b \right\rangle \right\rangle c ) &=& \tau ( \left\langle \left\langle ( d a )^*, d b . c \right\rangle \right\rangle )\\
		  &=& \left\langle d a^*, d b. c \right\rangle\\
				&=& \left\langle a^*, d^* ( d b. c ) \right\rangle \\
				&=& - \frac{1}{2} \left\langle a^*, ( b \cll ( c ) - \cll ( b ) c - \cll ( b c ) ) \right\rangle \\
				&=& - \frac{1}{2} \left\langle \cll ( b^* a^* ) - \cll ( b^* ) a^* - b^* \cll ( a^* ), c \right\rangle\\
				&=& - \frac{1}{2} \tau ( \left\langle \left\langle \cll ( b^* a^* ) - \cll ( b^* ) a^* - b^* \cll ( a^* ), c \right\rangle\right\rangle ).
				\end{eqnarray*}
				Thus, by the normality and faithfulness of $ \tau $ on $ \cla^{\prime \prime}, $ we conclude that 
				$$ g ( da \tensora db ) = \left\langle \left\langle ( da )^*, d b \right\rangle\right\rangle = \left\langle \left\langle da^*, d b \right\rangle\right\rangle = -(\frac{1}{2} \cll ( b^* a^* ) - \cll ( b^* ) a^* - b^* \cll ( a^* )). $$
				This proves the claim. Since $ \cll ( \A ) \subseteq \A, $ the proof of the proposition is complete. \qed
	 
\brmrk
Our $\mathcal{H}^1_D$ and $d$ are the same as the bimodule and derivation respectively constructed by Cipriani and Sauvagoet (\cite{cs}) from the Dirichlet form $$ (a,b) \mapsto -\langle \mathcal{L}(a),b \rangle; \ a,b \in {\rm Dom}((- \mathcal{L})^{\frac{1}{2}}). $$
\ermrk
		
It also follows from the definition of inner product that the map $ V_g $ is one-one. However, the invertibility of $ V_g, $ which is the nondegeneracy in our sense, has to be verified case by case.

\section{Levi-Civita connections on spectral triples} \label{mainsection}

	Let us recall the notations $\Psym$, $\sigma$ from Definition \ref{defn24thmay2018} and the map $\sigma^{{\rm can}}$ from Theorem \ref{skeide2}. The goal of this section is to prove the following theorem:
	\bthm \label{lcexistenceforbilinear}
	Suppose $(\A,\mathcal{H},D)$ is a spectral triple such that the following conditions hold:
	\begin{item}
		\item[i.] $\E := \oneform$ is a finitely generated projective right $\A$ module,
		\item[ii.] The map $u^\E:\Ecenter \otimes_{\Acenter} \A \rightarrow \E$ defined by
		$$u^\E(\sum_i e_i' \otimes_{\Acenter} a_i)=\sum_i e_i'a_i$$ is an isomorphism of vector spaces,
		\item[iii.] Suppose that there exists a right $\A$-module $\F$ such that $\E \tensora \E = {\rm Ker}(\wedge) \oplus \F$ as right $\A$-modules,
		\item[iv.] $\sigma=\sigma^{{\rm can}}$.
	\end{item}
	
	  If $g$ is any pseudo-Riemannian bilinear metric on $\E$, then there exists a unique connection on $\E$ which is torsion-less and compatible with $g$ (in the sense to be defined in Subsection \ref{metric_compatibility}). In particular, this applies to the candidate of a Riemannian bilinear map in Definition \ref{defn25thmay2018}.
	\ethm
	
	For the sake of convenience, we introduce the following definition:
\bdfn
	We call a spectral triple $(\A, \mathcal{H}, D)$ which satisfies the hypothesis of Theorem \ref{lcexistenceforbilinear} a tame spectral triple.
\edfn	


It is worthwhile to explain the significance of the equality $\sigma = \sigma^{\rm can}$. This is what we record in the following proposition:

\bppsn \label{17thoct20192}
Suppose $(\A, \mathcal{H}, D)$ is a tame spectral triple 
\begin{itemize}
\item[i.]$\E := \oneform$. Then the decomposition $\E \tensora \E = {\rm Ker}(\wedge) \oplus \F$ on simple tensors is explicitly given by
\[ \omega \tensora \eta a = \frac{1}{2}(\omega \tensora \eta a + \eta \tensora \omega a)  + \frac{1}{2}(\omega \tensora \eta a - \eta \tensora \omega a) , \]
for all $\omega$, $\eta$ in $\Ecenter$ and for all $a$ in $\A$.
\item[ii.] If $\E$ is a free right $\A$-module with a central basis $\{ e_1, e_2, ..., e_n \}$ and $g$ is a pseudo-Riemannian metric on $\E$, then the components $g_{ij} = g(e_i \tensora e_j)$ of $g$ are symmetric in $i$ and $j$.
\end{itemize}
\eppsn
{\bf Proof:}
We observe that since $\E \tensora \E$ is centered, any element of $\E \tensora \E$ is a $\IC$-linear sum of elements of the form $\omega \tensora \eta a$, where $\omega$, $\eta$ are in $\Ecenter$ and for $a$ in $\A$. Since $\sigma = \sigma^{\rm can}$,
\[ \Psym(\omega \tensora \eta a) = \frac{1}{2}(1 + \sigma^{\rm can})(\omega \tensora \eta a)= \frac{1}{2}(\omega \tensora \eta a + \eta \tensora \omega a) \]
\[ {\rm and} \quad (1 - \Psym)(\omega \tensora \eta a) = \frac{1}{2}(1 - \sigma^{\rm can})(\omega \tensora \eta a)= \frac{1}{2}(\omega \tensora \eta a - \eta \tensora \omega a). \]
Since $\Psym$ is an idempotent, this implies that $\frac{1}{2}(\omega \tensora \eta a + \eta \tensora \omega a)$ is in ${\rm Ran}(\Psym) = {\rm Ker}(\wedge)$ and $\frac{1}{2}(\omega \tensora \eta a - \eta \tensora \omega a)$ is in ${\rm Ker}(\Psym) = \F$.\\
Now we prove the second assertion. Since $g$ is a pseudo-Riemannian metric, and $\sigma = \sigma^{\rm can}$, we have
\[ g_{ij} = g(e_i \tensora e_j) = g \circ \sigma (e_i \tensora e_j) = g(e_j \tensora e_i) = g_{ji}. \] This finished the proof.
\qed

Let us make the following observation at this point:
\blmma \label{29thmay2018}
Suppose that $(\A,\mathcal{H},D)$ is a tame spectral triple. Then $\Psym$ is an $\A$-$\A$ bimodule map. In particular, a tame spectral triple is a quasi-tame spectral triple.
\elmma

{\bf Proof:} Since equation \eqref{splitting25thmay2018} is satisfied, $\Psym$ is a right $\A$-linear map by definition. But as $\sigma=\sigma^{\rm can}$ and $\sigma^{\rm can}$ is $\A$-$\A$ bilinear by Theorem \ref{skeide2}, $\sigma$ is $\A$-$\A$ bilinear. Therefore $\Psym = \frac{1+\sigma}{2}$ is also $\A$-$\A$ bilinear. \qed

We have already defined the torsion of a connection in Subsection \ref{sptripleforms}. In the next subsection, we formulate a notion of metric compatibility of a connection on the space of one-forms of a spectral triple satisfying some assumptions. We also prove (Theorem \ref{Phigiso}), a result which gives a sufficient condition for the existence and uniqueness of Levi-Civita connections. Subsection \ref{subsectionmainproof} is devoted to the proof of Theorem \ref{lcexistenceforbilinear}. 

 Throughout this section, we will work with tame spectral triples and continue to denote $\oneform$ by the symbol $\E$. By Lemma \ref{29thmay2018}, we are allowed to use all results concerning a quasi-tame spectral triple proved before and also the $\A$-$\A$ bilinearity of the map $\Psym$.

 \subsection{The metric compatibility of a connection on $\oneform$} \label{metric_compatibility}
	Throughout this subsection, $ ( \A, \mathcal{H}, D ) $ will denote a tame spectral triple. Moreover, $ g $ will denote  any pseudo-Riemannian bilinear metric ( not necessarily the canonical one ) on the bimodule $\E$ of one-forms. 

\bdfn
Let $\nabla$ be a connection on $\E$. Then we define $\Pi_g^0(\nabla):\Ecenter\tensorc\Ecenter \rightarrow \E$ by the map given by
$$\Pi_g^0(\nabla)(\omega \tensorc \eta)=(g\tensora {\rm id})\sigma_{23}(\nabla(\omega)\tensora \eta + \nabla(\eta) \tensora \omega).$$
\edfn
We have the following:

\blmma \label{pigzero}
$\Pi_g^0(\nabla)$ descends to a map from $\Ecenter \otimes_{\Acenter} \Ecenter$ to $\E$, to be denoted by the same notation. Moreover, for all $a' \in \Acenter$ and $\omega,\eta \in \Ecenter$
\be \Pi_g^0(\nabla)(\omega \otimes_{\Acenter} \eta a')= \Pi_g^0(\omega \otimes_{\Acenter} \eta)a' + g(\omega \tensora \eta)da' \label{pigone}. \ee
\elmma

{\bf Proof:} We write $\nabla(\eta)=\sum_i \eta_i^{(1)} \tensora \eta_i^{(2)}$, where $\eta_i^{(1)},\eta_i^{(2)} \in \E$ and the sum has finitely many terms. We have $\sigma_{23}(\omega\tensora da'\tensora \eta)=\omega\tensora \eta\tensora da'$, $\sigma_{23}(\nabla(\eta)a'\tensora \omega)=\sum_i \eta_i^{(1)}\tensora \omega \tensora \eta_i^{(2)}a'$. Using these, we get
\begin{eqnarray*}
&&\Pi_g^0(\nabla)(\omega a'\tensorc \eta)\\
 &=&(g\tensora {\rm id})\sigma_{23}(\nabla(\omega)a'\tensora \eta +\omega \tensora da' \tensora \eta + \nabla(\eta)\tensora \omega a')\\
&=&(g \tensora {\rm id})\sigma_{23}(\nabla(\omega) \tensora \eta a') + g(\omega \tensora \eta)da'\\
&+& \sum_i g(\eta_i^{(1)} \tensora \omega) \eta_i^{(2)}a' ~ \text{( by Proposition \ref{28thmay20183} )}\\
&=&(g \tensora {\rm id})\sigma_{23}(\nabla(\omega)\tensora \eta a')+g(\eta \tensora \omega)da'\\
&+&\sum_i g(\eta_i^{(1)} \tensora \omega)\eta_i^{(2)}a'\\
&=&(g \tensora id)\sigma_{23}(\nabla(\omega) \tensora \eta a' + \eta \tensora da' \tensora \omega + \nabla(\eta)a' \tensora \omega)\\
&=&(g \tensora {\rm id})\sigma_{23}(\nabla(\omega)\tensora \eta a' + \nabla(\eta a') \tensora \omega)\\
&=&\Pi_g^0(\nabla)(\omega \tensorc \eta a')\\
&=&\Pi_g^0(\nabla)(\omega \tensorc a'\eta).
\end{eqnarray*}

This proves the first assertion.
To prove the second assertion we make the following computation: for $a' \in \Acenter$ and $\omega,\eta \in \Ecenter$, we have:
\begin{align*}
\Pi_g^{(0)}(\nabla)(\omega \otimes_{\Acenter} \eta a')&=(g \tensora {\rm id})\sigma_{23}(\nabla(\omega)\tensora \eta a' + \nabla(\eta a') \tensora \omega )\\
&= (g \tensora {\rm id})\sigma_{23}(\nabla(\omega)\tensora \eta a' + \nabla(\eta )a' \tensora \omega + \eta\tensora da' \tensora \omega) \\
&\text{ ( since $\nabla$ is a connection )}\\
&= (g \tensora {\rm id})\sigma_{23}(\nabla(\omega)\tensora \eta + \nabla(\eta ) \tensora \omega )a' + (g \tensora {\rm id})(\eta \tensora \omega \tensora da') \\
&\text{ ( using Lemma \ref{sigmaoncenteroneside} )}  \\
&= \Pi_g^0(\omega \otimes_{\Acenter} \eta)a' + g(\omega \tensora \eta)da' \\
&\text{ ( using Proposition \ref{28thmay20183} )} 
\end{align*} 
\qed

\bdfn \label{21oct20191sm}
We define a map from $ \Ecenter \otimes_{\Acenter} \Ecenter \otimes_{\Acenter} \A $ to $ \rightarrow \E \tensora \E$ by the formula:
$$ u^{\E \tensora \E}:= (u^{\E} \tensora {\rm id}_{\E})\circ({\rm id}_{\Ecenter}\otimes_{\Acenter} u^{\E}).$$
We note that $ u^{\E \tensora \E} $ is an isomorphism since $u^\E$ is so.

 For all $ \omega,\eta \in \Ecenter $
 and $a \in \A$, define $\Pi_g(\nabla):\E \tensora \E \rightarrow \E$ by
$$\Pi_g(\nabla)\circ u^{\E \tensora \E}(\omega \otimes_{\Acenter} \eta \otimes_{\Acenter} a)= \Pi_g^0(\nabla)(\omega \otimes_{\Acenter}\eta)a + g(\omega \tensora \eta)da$$
\edfn

Therefore, for $\omega,\eta \in \Ecenter$ and $a \in \A$, we have
$$\Pi_g(\nabla)(\omega\tensora\eta a)= \Pi_g^0(\nabla)(\omega\otimes_{\Acenter}\eta)a +g(\omega\tensora\eta)da$$

\bppsn \label{proppignabla}
 The map $ \Pi_g ( \nabla ) $ defined in Definition \ref{21oct20191sm} is a well defined map from $ \E \tensora \E $ to $ \E. $ Moreover, $ \Pi_g ( \nabla ) - dg: \E \tensora \E \rightarrow \E $ is right $\A$-linear.
\eppsn
 {\bf Proof:} Since the map $u^{\E \tensora \E}$ is an isomorphism, it is enough to check that the map 
	$$ \Pi_g ( \nabla ) \circ u^{\E \tensora \E}: \Ecenter \otimes_{\Acenter} \Ecenter \otimes_{\Acenter} \A \rightarrow \E $$
	is well defined. For $ \omega, \eta \in \Ecenter, ~ a \in \Acenter, b \in \A, $ the equalities
	$$ \Pi_g ( \nabla ) \circ u^{\E \tensora \E} ( \omega a \tensorc \eta \tensorc b ) = \Pi_g ( \nabla ) \circ u^{\E \tensora \E} ( \omega \tensorc a \eta \tensorc b ) ~{\rm and}$$
	$$ \Pi_g ( \nabla ) \circ u^{\E \tensora \E} ( \omega \tensorc \eta a \tensorc b ) = \Pi_g ( \nabla ) \circ u^{\E \tensora \E} ( \omega \tensorc \eta \tensorc ab ) $$
	follow from Lemma \ref{pigzero} and \eqref{pigone} respectively.	For proving the right $\A$-linearity of the map $ \Pi_g ( \nabla ) - dg $ it is sufficient to evaluate it on $\omega\tensora \eta ab$, where $ \omega, \eta \in \Ecenter, ~ a , b \in \A $, since $u^{\E\tensora \E}$ is an isomorphism.
	\begin{align*}
	(\Pi_g ( \nabla ) - dg)(\omega\tensora \eta ab)&=\Pi_g^0(\nabla)(\omega \otimes_{\Acenter}\eta)ab +g(\omega\tensora \eta)d(ab)-dg(\omega\tensora\eta ab)\\
	&=\Pi_g^0(\nabla)(\omega \otimes_{\Acenter}\eta)ab +g(\omega \tensora\eta)(da.b+adb) -dg(\omega\tensora\eta a)b -g(\omega\tensora\eta a)db\\
	&=(\Pi_g^0(\nabla)(\omega \otimes_{\Acenter}\eta)a +g(\omega\tensora \eta)d(a)-dg(\omega\tensora\eta a))b\\
	&=((\Pi_g(\nabla)-dg)(\omega\tensora\eta a)b
	\end{align*}
		 \qed

\bdfn \label{19thaugust20191}
A connection $\nabla$ is said to be compatible with a pseudo-Riemannian metric $g$ if $\Pi_g(\nabla)=dg$.
\edfn

\bppsn
The above definition of metric compatibility coincides with that in the classical case.
\eppsn

{\bf Proof:} If $(M,g)$ is a pseudo-Riemannian manifold, then a connection $\nabla$ on $\oneformclassical$ is said to be compatible with $g$ if and only if
$$g_{13}(\nabla(\omega)\otimes_{C^{\infty}(M)} \eta + \nabla(\eta)\otimes_{C^{\infty}(M)} \omega)=dg(\omega\otimes_{C^{\infty}(M)}\eta),$$
for all $\omega,\eta \in \oneformclassical$, where $g_{13}=(g\otimes_{C^{\infty}(M)} {\rm id})\sigma_{23}^{{\rm can}}$.\\
We have $\A = \Acenter = {C^{\infty}(M)}$ and $\E = \Ecenter = \oneformclassical$. If $\nabla$ is metric compatible in our sense, then for $\omega,\eta \in C^\infty(M)$ we have
 \begin{align*}
 g_{13}(\nabla(\omega)\otimes_{C^\infty(M)}\eta +\nabla(\eta)\otimes_{C^\infty(M)}\omega) &=(g\otimes_{C^\infty(M)}{\rm id})\sigma_{23}(\nabla(\omega)\otimes_{C^\infty(M)}\eta +\nabla(\eta)\otimes_{C^\infty(M)}\omega)\\
 &=\Pi_g^0(\nabla)(\omega\otimes_{C^\infty(M)}\eta)\\
 &=\Pi_g(\nabla)(\omega\otimes_{C^\infty(M)}\eta )\\
 &=dg(\omega\otimes_{C^\infty(M)}\eta ).
 \end{align*}
Thus, our definition of metric compatibility coincides with that in the classical case. \qed

\blmma \label{remark4.29}
\begin{enumerate}
	\item[(i)] The map $\Pi_g(\nabla)-dg \in {\rm Hom}_{\A}(\E \tensora \E,\E)$ is determined by its restriction on $\E \tensorsym \E$ for any connection $\nabla$ and can be viewed as an element of ${\rm Hom}_{\A}(\E \tensorsym \E,\E)$
	\item[(ii)] For any two torsion-less connections $\nabla_1$ and $\nabla_2$, $\nabla_1 - \nabla_2 \in {\rm Hom}_{\A}(\E,\E \tensorsym \E)$
\end{enumerate}
\elmma
{\bf Proof:} By the definition of $\Pi_g^0(\nabla)$ and the equality $g\circ \sigma=g$, it follows that $\Pi_g^0(\nabla)\circ\sigma=\Pi_g^0(\nabla)$ on $\Ecenter\otimes_{\Acenter}\Ecenter$. Now for $\omega,\eta \in \Ecenter$ and $a \in \A$, we have 
\begin{align*}
(\Pi_g(\nabla)-dg)\circ\sigma(\omega\tensora\eta a) &=(\Pi_g(\nabla)-dg)(\sigma(\omega\tensora\eta)a) =(\Pi_g(\nabla)-dg)\circ\sigma(\omega\tensora\eta)a\\
&=(\Pi_g(\nabla)-dg)(\omega\tensora\eta)a =(\Pi_g(\nabla)-dg)(\omega\tensora\eta a),
\end{align*}
since $\Pi_g(\nabla)-dg$ is right $\A$-linear by Proposition \ref{proppignabla}.
Therefore, $\Pi_g(\nabla)-dg=(\Pi_g(\nabla)-dg)\circ \sigma$ on the whole of $\E\tensora\E$. This proves $(i)$.\\
Now we prove $(ii)$. If $\nabla_1$ and $\nabla_2$ are two torsion-less connections, $\wedge \circ \nabla_1 = -d = \wedge \circ \nabla_2$. Therefore, ${\rm Ran}(\nabla_1-\nabla_2)\subseteq {\rm Ker}(\wedge) =\E\tensorsym\E$. Moreover $(\nabla_1-\nabla_2)(\omega a)=\nabla_1(\omega)a-\nabla_2(\omega)a$ for $ \omega \in \E$ and for $ a \in \A$. Hence, $\nabla_1 - \nabla_2 \in {\rm Hom}_{\A}(\E,\E \tensorsym \E)$.
\qed

\bdfn We define a map
$$\Phi_g : {\rm Hom}_{\A}(\E, \E \tensorsym \E) \rightarrow {\rm Hom}_{\A}(\E \tensorsym \E, \E) ~ {\rm by} $$
$$ \Phi_g(L)=(g \tensora {\rm id}) \sigma_{23} (L \tensora {\rm id})(1+\sigma).$$
\edfn

\bppsn
$\Phi_g$ is a right $\A$-linear map.
\eppsn
{\bf Proof:}  
Let $\omega,\eta \in \Ecenter$, and $a,b \in \A$ and $L \in {\rm Hom}_\A(\E,\E \tensorsym\E)$. Then by using Proposition \ref{28thmay20183} and Lemma \ref{29thmay2018} we obtain
\begin{align*}
\Phi_g(La)(\omega\tensora\eta b)&=(g \tensora {\rm id})\sigma_{23}(La \tensora {\rm id})(1 + \sigma)(\omega \tensora\eta b)\\
&=(g \tensora {\rm id})\sigma_{23}(La \tensora {\rm id})(\omega\tensora\eta b + \eta\tensora\omega b)\\	
&=(g \tensora {\rm id})\sigma_{23}(L(a\omega)\tensora\eta b + L(a\eta) \tensora \omega b) \\
&=(g \tensora {\rm id})\sigma_{23}(L \tensora {\rm id})( a \omega \tensora\eta b + a \eta \tensora \omega b )\\
&=(g \tensora {\rm id})\sigma_{23}(L \tensora {\rm id})(1 + \sigma)(a(\omega \tensora\eta b ))\\
&=(\Phi_g(L)a)(\omega\tensora\eta b).
\end{align*}
 Hence we have that $\Phi_g(La)=\Phi_g(L)a$.
\qed

Now we are in a position to prove the following result: 
\bthm \label{Phigiso}
	If $ \Phi_g:{\rm Hom}_{\A}(\E, \E \tensorsym \E) \rightarrow {\rm Hom}_{\A}(\E \tensorsym \E, \E) $ is an isomorphism of right $\A$ modules, then there exists a unique connection on $\E$ which is torsion-less and compatible with $g.$
	\ethm
	{\bf Proof:}
	We recall the torsion-less connection $\nabla_0$ constructed in Lemma \ref{torsionless}. By $(i)$ of Lemma \ref{remark4.29}, $dg-\Pi_g(\nabla_0) \in {\rm Hom}_{\A}(\E \tensorsym \E,\E)$. Since $\Phi_g$ is an isomorphism from ${\rm Hom}_{\A}(\E,\E \tensorsym \E)$ to ${\rm Hom}_{\A}(\E \tensorsym \E,\E)$ there exists a unique element $\Phi_g^{-1}(dg - \Pi_g(\nabla_0)) \in {\rm Hom}_{\A}(\E,\E \tensorsym \E)$. Define the $\IC$-linear map $$\nabla:=\nabla_0+\Phi_g^{-1}(dg - \Pi_g(\nabla_0)).$$
	We claim that $\nabla$ is a torsion-less connection on $\E$ which is compatible with $g$. Indeed, if $\omega \in \E$ and $a \in \A$, we have
	\begin{align*}
	\nabla(\omega a)&=\nabla_0(\omega)a + \omega\tensora da + \Phi_g^{-1}(dg - \Pi_g(\nabla_0))(\omega)a\\
	&= \nabla(\omega)a + \omega \tensora da. 
	\end{align*}
	so that $\nabla$ is a connection. That $\nabla$ is a torsion-less connection is verified from the following:
	\begin{align*}
	\wedge \circ \nabla &= \wedge \circ \nabla_0 + \wedge \circ \Phi_g^{-1}(dg - \Pi_g(\nabla_0))\\
	&= \wedge \circ \nabla_0 \text{ ( since ${\rm Ran}(\Phi_g^{-1})(dg - \Pi_g(\nabla_0) \subseteq \E \tensorsym \E = {\rm Ker} (\wedge)$ ) }\\
	&= -d.
	\end{align*}
	We note that this in particular implies that $\nabla -\nabla_0 \in {\rm Hom}_{\A}(\E,\E \tensorsym \E)$ so that $\Phi_g(\nabla-\nabla_0)$ is well-defined. Moreover, for $\omega,\eta \in \Ecenter$ and $a \in \A$, we have
	\begin{align*}
	(\Pi_g(\nabla)-\Pi_g(\nabla_0))(\omega\tensora\eta a)&=\Pi_g^0(\nabla)(\omega \otimes_{\Acenter} \eta)a-\Pi_g^0(\nabla_0)(\omega \otimes_{\Acenter}\eta)a\\
	&=(g \tensora {\rm id})\sigma_{23}\Big((\nabla(\omega)\tensora\eta +\nabla(\eta)\tensora\omega)-(\nabla_0(\omega)\tensora\eta +\nabla_0(\eta)\tensora\omega)\Big)a\\
	&=(g \tensora {\rm id})\sigma_{23}((\nabla-\nabla_0)\tensora{\rm id})(1+\sigma)(\omega \tensora \eta a)\\
	&=\Phi_g(\nabla-\nabla_0)(\omega \tensora \eta a).
	\end{align*}
	Therefore, $\Phi_g(\nabla-\nabla_0)=\Pi_g(\nabla) - \Pi_g(\nabla_0)$.
	Since $\Phi_g(\nabla - \nabla_0) = dg - \Pi_g(\nabla_0)$ by the definition of $\nabla$, we have $\Pi_g(\nabla)=dg$.
	Therefore, $\nabla$ is compatible with $g$.\\
	To show uniqueness, suppose $\nabla'$ is another torsion-less connection compatible with the metric $g$. Then exactly as above, $\nabla-\nabla' \in {\rm Hom}_{\A}(\E,\E \tensorsym \E)$ and
  $$\Phi_g(\nabla-\nabla')=\Pi_g(\nabla)-\Pi_g(\nabla')=dg -dg=0,$$
	where we have used the fact that $\nabla$ and $\nabla'$ are compatible with $g$.
	Hence, $\nabla=\nabla'$, as $\Phi_g$ is an isomorphism.	
	\qed
	
\brmrk \label{19thaugust20193}
 The definition and results of this subsection go through verbatim even for a right $\A$-linear pseudo-Riemannian (i.e, not necessarily bilinear) metric. This will be used in \cite{article2}. 
\ermrk

\subsection{Existence and uniqueness of Levi-Civita connections for a class of spectral triples} \label{subsectionmainproof}

In this subsection, we prove Theorem \ref{lcexistenceforbilinear} by utilizing Theorem \ref{Phigiso}. As observed before, the isomorphism of the map $u^{\E}$ implies that $\E$ is centered. Therefore, we will freely use the fact that $\E$ is centered throughout this section, sometimes without mentioning.	The map $\zeta_{\E\tensora\E,\E}$ will be as introduced in Proposition \ref{28thmay2018}.

	\blmma \label{20thfeb20182}
Let $ L $ be an element of $ {\rm Hom}_\A ( \E, \E \tensora \E )  $ such that  $ \zeta^{-1}_{\E \tensora \E, \E} ( L ) = \xi \tensora \eta \tensora V_g ( \omega ) $ for some $ \xi, \eta, \omega \in \E. $
\begin{enumerate} 

\item Then for all $ e $ in $\E,$ we have
 $$ L ( e ) = \xi \tensora \eta g ( \omega \tensora e ). $$
\item Let us define then an element $ L^\prime \in {\rm Hom}_\A (\E,\E \tensora \E ) $ by the equation
 $$ \zeta^{-1}_{\E \tensora \E, \E} ( L^\prime ) = \eta \tensora \xi \tensora V_g ( \omega ). $$
	If $ L \in {\rm Hom}_\A ( \E, \E \tensorsym \E ) $ and $ \xi, \eta, \omega \in \Ecenter, $ then $ L = L^\prime $ as elements of $ {\rm Hom}_\A ( \E, \E \tensora \E ). $ Moreover, 
		 $$ \xi \tensora \eta \tensora V_g ( \omega ) = \eta \tensora \xi \tensora V_g ( \omega ).$$
\item The set $ \{ \zeta_{\E \tensora \E, \E} ( \xi \tensora \eta \tensora V_g ( \omega ): \xi, \eta, \omega \in \Ecenter \} $ is right $\A$-total in $ {\rm Hom}_\A ( \E, \E \tensora \E ). $		
\end{enumerate}	
\elmma

{\bf Proof:} Let $ e $ denote an element of $ \E. $ By the definition of $ \zeta_{\E \tensora \E, \E}, $ it follows that 
$$ L ( e ) = \xi \tensora \eta V_g ( \omega ) ( e ) = \xi \tensora \eta g ( \omega \tensora e ). $$
Now we prove part 2. By part 1., we have
 $$ \Psym L ( e ) = \frac{1}{2} ( \xi \tensora \eta + \eta \tensora \xi ) g ( \omega \tensora e ). $$ 
 Since $ L ( e ) \in \E \tensorsym \E, $ we have $ \Psym L ( e ) = L ( e ). $ 
	Therefore, $ \frac{1}{2} ( \xi \tensora \eta + \eta \tensora \xi ) g ( \omega \tensora e ) = \xi \tensora \eta g ( \omega \tensora e ) $ which implies that $ \xi \tensora \eta g ( \omega \tensora e ) = \eta \tensora \xi g ( \omega \tensora e ). $ This proves that $ L ( e ) = L^{\prime} ( e ). $
	Hence, 
	$$ \xi \tensora \eta \tensora V_g ( \omega ) = \zeta^{-1}_{\E \tensora \E, \E} ( L ) = \zeta^{-1}_{\E \tensora \E, \E} ( L^\prime ) = \eta \tensora \xi \tensora V_g ( \omega ). $$
 Finally, for part 3., we note that since $ g $ is bilinear, the set $ S = \{ \xi \tensora \eta \tensora V_g ( \omega ): \xi, \eta, \omega \in \Ecenter \} $ is right $\A$-total in $ \E \tensora \E \tensora \E^{*} $ and therefore $ \zeta_{\E \tensora \E, \E} ( S ) $ is right $ \A $-total in $ {\rm Hom}_\A ( \E, \E \tensora \E ). $ 
\qed

\blmma \label{lemma3}
For all $\omega,\eta \in \E,~ \Vgtwo \sigma ( \omega \tensora \eta ) = \Vgtwo ( \omega \tensora \eta ) \sigma. $ In particular, $ \Psym = \Psym^\ast. $
\elmma
 {\bf Proof:} By Lemma \ref{29thmay2018} $\Psym$ is bilinear. Therefore, the map $\Psym^{*}$ (defined as in Lemma \ref{lemma4}) makes sense. It is enough to prove that for all $ \omega, \eta, \omega^\prime, \eta^\prime \in \Ecenter, $
  $$ \Vgtwo ( \sigma ( \omega \tensora \eta ) ) ( \omega^\prime \tensora \eta^\prime ) = \Vgtwo ( \omega \tensora \eta ) \sigma ( \omega^\prime \tensora \eta^\prime ). $$
	We compute
	 \begin{eqnarray*}
		 \Vgtwo ( \sigma ( \omega \tensora \eta ) ) ( \omega^\prime \tensora \eta^\prime )
		 &=& g^{(2)} ( ( \eta \tensora \omega ) \tensora ( \omega^\prime \tensora \eta^\prime ) )\\
			&=& g ( \eta \tensora \eta^\prime ) g ( \omega \tensora \omega^\prime )\\
			&=& g ( \omega \tensora \omega^\prime ) g ( \eta \tensora \eta^\prime ) ~ {\rm (} ~ {\rm by}~ {\rm Proposition}~ \ref{17thoct20191}  ~{\rm )} \\
			&=& \Vgtwo ( \omega \tensora \eta ) \sigma ( \omega^\prime \tensora \eta^\prime ).
		 \end{eqnarray*}
		This finishes the proof.
	\qed

\blmma \label{lemma1}
 Let $ L \in {\rm Hom}_\A ( \E, \E \tensora \E ) $ be such that $ \zeta^{-1}_{\E \tensora \E, \E} ( L ) = \xi \tensora \eta \tensora V_g ( \omega ) $ for some $ \xi, \eta, \omega ~ \in ~ \Ecenter. $ Then 
 \be \label{lemma1eqn} \Phi_g ( L ) = \zeta_{\E, \E \tensora \E} ( \eta \tensora \Vgtwo ( \xi \tensora \omega + \omega \tensora \xi ) ). \ee
\elmma
{\bf Proof:} Let us observe that it is enough to prove that for all $ \omega^\prime, \eta^\prime \in \Ecenter, $ 
 $$ \Phi_g ( L ) ( \omega^\prime \tensora \eta^\prime ) = \eta \Vgtwo ( \xi \tensora \omega + \omega \tensora \xi ) ( \omega^\prime \tensora \eta^\prime ).$$
By using part 1. of Lemma \ref{20thfeb20182}, we compute
 \begin{eqnarray*} 
 \Phi_g ( L ) ( \omega^\prime \tensora \eta^\prime )
  &=& ( g \tensora {\rm id} ) \sigma_{23} ( L ( \omega^\prime ) \tensora \eta^\prime + L ( \eta^\prime ) \tensora \omega^\prime )\\
	 &=& ( g \tensora {\rm id} ) \sigma_{23} ( \xi \tensora \eta g ( \omega \tensora \omega^\prime ) \tensora \eta^\prime + \xi \tensora \eta g ( \omega \tensora \eta^\prime ) \tensora \omega^\prime )\\
	&=& ( g \tensora {\rm id} ) ( \xi \tensora \eta^\prime \tensora \eta g ( \omega \tensora \omega^\prime ) + \xi \tensora \omega^\prime \tensora \eta g ( \omega \tensora \eta^\prime ) )\\
	&=& g ( \xi \tensora \eta^\prime ) \eta g ( \omega \tensora \omega^\prime ) + g ( \xi \tensora \omega^\prime ) \eta g ( \omega \tensora \eta^\prime )\\
	&=& \eta g ( \xi \tensora \eta^\prime ) g ( \omega \tensora \omega^\prime ) + \eta g ( g ( \xi \tensora \omega^\prime ) \omega \tensora \eta^\prime ) ~ {\rm ( } ~ {\rm since} ~ g ~ {\rm is} ~ {\rm bilinear} ~ {\rm )} \\
	&=& \eta \Vgtwo ( \xi \tensora \omega + \omega \tensora \xi ) ( \omega^\prime \tensora \eta^\prime ).	
 \end{eqnarray*} \qed

	\bppsn 
	 Let $ L \in {\rm Hom}_\A ( \E, \E \tensorsym \E ). $ Then 
		\be \label{21stfeb2018} \frac{1}{2} \Phi_g ( L ) = \zeta_{\E, \E \tensora \E} ( {\rm id} \tensora \Vgtwo ) ( P_{\rm sym} )_{23} ( {\rm id} \tensora V^{-1}_g ) \zeta^{-1}_{\E \tensora \E, \E} ( L ) . \ee
	\eppsn 
	{\bf Proof:} Let $ L \in {\rm Hom}_\A ( \E, \E \tensorsym \E ) $ be such that $ \zeta^{-1}_{\E \tensora \E, \E} ( L ) = \xi \tensora \eta \tensora V_g ( \omega ) $ for some $ \xi,\eta, \omega \in \Ecenter. $
	  Then by part 2. of Lemma \ref{20thfeb20182}, we have $ \xi \tensora \eta \tensora V_g ( \omega ) = \eta \tensora \xi \tensora V_g ( \omega ). $
		Therefore, 
		\begin{eqnarray*}
		 && \zeta_{\E, \E \tensora \E} ( ( {\rm id} \tensora \Vgtwo ) ( P_{\rm sym} )_{23} ( {\rm id} \tensora V^{-1}_g ) \zeta^{-1}_{\E \tensora \E, \E} ( L ) )\\
		 &=& \zeta_{\E, \E \tensora \E} ( ( {\rm id} \tensora \Vgtwo ) ( P_{\rm sym} )_{23} ( {\rm id} \tensora V^{-1}_g ) ( \xi \tensora \eta \tensora V_g ( \omega ) ) )\\
			&=& \zeta_{\E, \E \tensora \E} ( ( {\rm id} \tensora \Vgtwo ) ( P_{\rm sym} )_{23} ( \eta \tensora \xi \tensora \omega ) )\\
			&=& \frac{1}{2} \eta \tensora \Vgtwo ( \xi \tensora \omega + \omega \tensora \xi ) ~ {\rm ( } ~ {\rm since} ~ \xi, \omega ~ \in ~ \Ecenter ~ {\rm )}\\
		&=& \frac{1}{2} \Phi_g ( L ) ~ {\rm ( } ~ {\rm by} ~ {\rm Lemma} ~ \ref{lemma1} ~ {\rm )}.\\
			\end{eqnarray*}
			Thus, we have proved \eqref{21stfeb2018} for all $ L $ of the above form. But since the maps $ \zeta_{\E \tensora \E, \E}, ~ \Phi_g, \Vgtwo $ and $ \Psym $ are all right $ \A $-linear, we can conclude that \eqref{21stfeb2018} holds for all $ L $ in $ {\rm Hom}_\A ( \E, \E \tensorsym \E ) $ by appealing to part 3. of Lemma \ref{20thfeb20182}. 
			\qed
			
\blmma \label{lemma2}
$ \Vgtwo $ is nondegenerate as a map from $ \E \tensorsym \E $ to $ ( \E \tensorsym \E )^*. $ 
\elmma
{\bf Proof:} Let us start by claiming that $ ( \E \tensorsym \E )^\ast $ can be identified with the bimodule $ \{ \phi \in ( \E \tensora \E )^\ast: \phi \circ ( 1 - \Psym ) = 0 \} .$ Indeed, if $ \psi \in ( \E \tensorsym \E )^*, $ then $ \psi $ can be uniquely extended to an element $ \phi \in ( \E \tensora \E )^* $ by using the decomposition $ \E \tensora \E = {\rm Ran} ( \Psym ) \oplus {\rm Ran} ( 1 - \Psym ). $ Clearly, $ \psi = \phi \circ \Psym. $ Conversely, if $ \phi \in ( \E \tensora \E  )^* $ then $ \phi \circ \Psym $ defines an element of $ ( \E \tensorsym \E  )^*. $ This proves our claim.

Now we use our claim to prove that $ \Vgtwo $ is one-one and onto as a map from $\E \tensorsym \E$ to $(\E \tensorsym \E)^{*}$. Let $ \phi \in ( \E \tensora \E )^* $ be such that $ \phi \circ ( 1 - \Psym ) = 0. $ Since $ \Vgtwo : \E \tensora \E \rightarrow (\E \tensora \E)^* $ is non-degenerate by Proposition \ref{vg2nondegenerate}, there exists $ \psi \in \E \tensora \E $ such that $ \Vgtwo ( \psi ) = \phi. $ We claim that $ \Psym \psi = \psi. $ Indeed, 
\begin{eqnarray*}
	\Vgtwo ( \Psym \psi ) &=& \Vgtwo ( \psi ) \circ \Psym = \phi \circ \Psym \\
	&=& \phi \circ \Psym + \phi \circ ( 1 - \Psym ) = \phi \\
	&=& \Vgtwo ( \psi ),
\end{eqnarray*}
where we have used Lemma \ref{lemma3}. By using Proposition \ref{vg2nondegenerate}, we conclude that $ \Psym \psi = \psi. $ This proves that $ \Vgtwo $ maps onto $ ( \E \tensorsym \E )^* .$

To prove that $ \Vgtwo $ is one-one as a map from $ \E \tensorsym \E $ to  $ ( \E \tensorsym \E )^*,$ let $ \psi \in  \E \tensorsym \E  $ be such that $ \Vgtwo ( \psi ) \circ \Psym = 0. $ Therefore, by Lemma \ref{lemma3}, we have
\be \Vgtwo ( \psi ) = \Vgtwo \Psym ( \psi ) = \Vgtwo ( \psi ) \circ \Psym = 0, \ee
so that by Proposition \ref{vg2nondegenerate}, we have $ \psi = 0. $
\qed
			
\vspace{4mm}

{\bf Proof of Theorem \ref{lcexistenceforbilinear}:} We need to prove that the map $ \Phi_g $ is an isomorphism from $ {\rm Hom}_\A ( \E, \E \tensorsym \E ) $ to $ {\rm Hom}_\A ( ( \E \tensorsym \E ) \tensora \E ). $ By Lemma \ref{module_symm_iso}, the map
					$$ ( \Psym )_{23}: ( \E \tensorsym \E ) \tensora \E \rightarrow \E \tensora (  \E \tensorsym \E ) $$ 
					is an isomorphism of right $\A$ modules. Since $ ( {\rm id} \tensora V^{-1}_g ) \zeta^{-1}_{\E \tensora \E, \E} $ is an isomorphism from $ {\rm Hom}_\A ( \E, \E \tensorsym \E ) $ to $ ( \E \tensorsym \E  ) \tensora \E $ and $ \Vgtwo $ is an isomorphism from $ \E \tensorsym \E $ to $ ( \E \tensorsym \E )^*, $ by Lemma \ref{lemma2} we see that 
				$	\zeta_{\E, \E \tensora \E} ( {\rm id} \tensora \Vgtwo ) ( P_{\rm sym} )_{23} ( {\rm id} \tensora V^{-1}_g ) \zeta^{-1}_{\E \tensora \E, \E} $ is an isomorphism from $ {\rm Hom}_\A ( \E, \E \tensorsym \E ) $ to $ {\rm Hom}_\A ( ( \E \tensorsym \E ) \tensora \E ). $ Finally, the equation \eqref{21stfeb2018} implies that $ \Phi_g $ is an isomorphism. \qed

\subsection{A remark on the isomorphism of the map $u^{\E}$} \label{isomorphismremark}

Before going into the examples, it is worthwhile to derive a sufficient condition which ensures the isomorphism of the map $u^{\E}$, The following result will be crucially used in Section \ref{prelimdef}, where we prove the existence of the Levi-Civita connection on the Connes-Landi isospectral deformation of classical spectral triples.

\bppsn \label{isomorphismofEprime}
	Suppose $(\A,\mathcal{H},D)$ is a spectral triple. Suppose that there exists a unital subalgebra $\A'$ of $\Acenter$ and an $\A'$-submodule $\E'$ of $\Ecenter$ such that $\E'$ is projective and finitely generated over $\A'$. If the map
	$$u^{\E}_{\E'}:\E'\otimes_{\A'}\A \rightarrow \E,$$
	defined by
	$$u^\E_{\E^\prime} \left( \sum_{i} e^{\prime}_i  \otimes_{\A^\prime} a_i \right) =\sum_{i} e^{\prime}_i a_i $$
	is an isomorphism of vector spaces, then $u^{\E}:\Ecenter \otimes_{\Acenter} A \rightarrow \E$ is an isomorphism. Moreover, if $\Ecenter$ is a finitely generated projective module over $\Acenter$, then $u^{\E}$ is an isomorphism if and only if there exists $\E'$ and $\A'$ such that $u^{\E}_{\E'}$ is an isomorphism.
\eppsn	
	{\bf Proof:} If $u^{\E}_{\E'}$ is an isomorphism, we claim that $\Ecenter \cong \E' \otimes_{\A'} \Acenter$. If our claim is true, then we have
$$ \clz ( \E ) \otimes_{\clz ( \A )} \A \cong \E^\prime \otimes_{\A^\prime} \clz ( \A ) \otimes_{\clz ( \A )} \A = \E^\prime \otimes_{\A^\prime} \A \cong \E, $$
so that $u^{\E}$ is an isomorphism. Thus, it is enough to prove our claim.\\
By a verbatim adaptation of the proof of Proposition \ref{assumption3implycentered}, we have that $ \E \cong \E^\prime \otimes_{\A^\prime} \A $ as bimodules where the bimodule structure of $ \E^\prime \otimes_{\A^\prime} \cla $ is defined by $ z_1 ( e^\prime \otimes_{\A^\prime} a ) z_2 = e^\prime \otimes_{\Aprime} z_1 a z_2. $ Thus, $ \{ \sum_i e^\prime_i \otimes_{\A^\prime} a_i: e^\prime_i \in \E^\prime, ~ a_i \in \clz ( \cla ) \} \subseteq \clz ( \E^\prime \otimes_{\A^\prime} \cla  ). $

For the reverse inclusion, let us suppose that there exists a free $\A^\prime$ module $\clf$ and an idempotent $P$ on $ \clf $ such that $ P ( \clf ) = \E^\prime. $ Let $ m_1, m_2, \cdots m_n $ be a basis of $ \clf. $
Therefore,
 $$ \E \cong \E^\prime \otimes_{\A^\prime} \A = P ( \clf ) \otimes_{\A^\prime} \A = ( P \otimes_{\A^\prime} {\rm id}_\A ) ( \clf \otimes_{\A^\prime} \A ). $$
Clearly, $ P \otimes_{\A^\prime} {\rm id}_\A $ is an idempotent on $ \clf \otimes_{\A^\prime} \A $ and thus for all $ y \in \E^\prime \otimes_{\A^\prime} \A \subseteq \clf \otimes_{\A^\prime} \A, $ we have
  \be \label{19thfeb2017} ( P \otimes_{\A^\prime} {\rm id}_\A ) ( y ) = y. \ee 
		On the other hand, $ \clz ( \E^\prime \otimes_{\A^\prime} \cla ) $ is also a submodule of $ \clf \otimes_{\A^\prime} \cla $ and if $ x $ is an element of $ \clz ( \E^\prime \otimes_{\A^\prime} \cla ), $ there exists unique elements $ a_i \in \A $ such that
		 $ x = \sum_i m_i \otimes_{\A^\prime} a_i. $
		Since $ x b = b x $ for all $ b \in \cla, $ we see that $ a_i \in \clz ( \A ) $ for all $ i. $
		Hence,
		$$ ( P \otimes_{\A^\prime} {\rm id}_\A ) ( x ) = \sum_i ( P \otimes_{\A^\prime} {\rm id}_\A ) ( m_i \otimes_{\A^\prime} a_i ) = \sum_i P ( m_i ) \otimes_{\A^\prime} a_i \in \E^\prime \otimes_{\A^\prime} \clz ( \cla ) .$$
		But by \eqref{19thfeb2017}, $ ( P \otimes_{\A^\prime} {\rm id}_\A ) ( x ) = x $ so that $ x \in \E^\prime \otimes_{\A^\prime} \clz ( \cla ).$ Since $ x $ is an arbitrary element of $ \clz ( \E^\prime \otimes_{\A^\prime} \cla ) \cong \clz ( \E ), $ this completes the proof. 
 \qed

\section{Levi-Civita connection for fuzzy 3-spheres} \label{fuzzysection}
			
			Let $G$ denote the compact Lie group $SU(2)$ and $V_j,$ $j\in \frac{1}{2}\mathbb{N}_0 $, denote the $(2j+1)$ dimensional irreducible representation of $SU(2)$.  Let $ k $ be a positive integer and $\mathcal{H}_0:=\bigoplus_{j=0,\frac{1}{2},...,\frac{k}{2}}V^{*}_j\tensorc V_j$ and $\mathcal{A}:=\mathcal{B}(\mathcal{H}_0)$. Let $W$ be the carrier vector space of the irreducible representation of the Clifford algebra generated by the vector space $T_eG$ with respect to the Killing form on $G$ as defined by equations (3.4) and (3.5) of \cite{frolich}. There exists a spectral triple $(\mathcal{A},\mathcal{H},D)$, where $\mathcal{H}:=\mathcal{H}_0\tensorc W$, called the ``fuzzy" or non-commutative 3-sphere. We refer to \cite{frolich} for the details. It turned out (Corollary 3.8 of \cite{frolich}) that there exists a unique real unitary connection (in the sense of Definition 2.24 of \cite{frolich}) on  $\E:=\oneform$ which is also torsion-less. In this section, we prove that with our definition of metric compatibility of a connection, there exists a unique Levi-Civita connection and this connection coincides with the real unitary and torsion-less connection in Corollary 3.8 of \cite{frolich}.\\
			In what follows, we will denote the elements $1 \tensorc \psi_i$ in the center of $\oneform$ as in \cite{frolich} by the symbol $e_i$, so that 
			$$ e_j \wedge e_k=-e_k \wedge e_j $$
			 and $\{e_i \wedge e_j:i \leq j\}$ is linearly independent. Let $\E:=\oneform$. We have the following result.

			\bthm \label{30thapril20180} (Equation (3.19) and Theorem 3.2, \cite{frolich})
			 The space of forms for the spectral triple $ ( \A, \clh, D  ) $ has the following description:\\
			 1. The module $ \E $ is isomorphic to $ {\rm Span} \{ e_i a_i : i = 1,2,3  \} $ and thus is a free right $\A$ module of rank three.\\
			2. The module $\twoform \cong ~ {\rm Span} \{  e_i \wedge e_j a_{ij} : a_{ij} = - a_{ji} \}$ is a free right $\A$ module of rank three.
	\ethm
	
	Moreover, it was also proven in \cite{frolich} that the space of three-forms is a free rank one module and all the higher forms are zero. The bimodule structure for $\E:=\oneform$ (and similarly, for the higher forms) is given by
	 $$ a ( b \tensorc \psi_i ) c = a b c \tensorc \psi_i = e_i abc. $$
	Thus, we can identify $ \E \tensora \E $ with $ {\rm Span} \{  e_i \tensora e_j a: i,j = 1,2,3 \}. $
	
	\blmma \label{30thapril20182}
	 ${\rm Ker} ( \wedge ) $ is generated (as a right $\A$ module) by the set $ \{ e_i \tensora e_i, e_i \tensora e_j + e_j \tensora e_i: i = 1,2,3   \}. $
	\elmma 
		{\bf Proof:} Throughout this proof, we will be using the fact that the elements $e_i$ are in $\Ecenter$. Let $ \omega = \sum_j e_j a_j, \eta = \sum_k e_k b_k $ be elements of $ \E. $ If $ \epsilon_{ijk} $ denotes the Levi-Civita tensor, i.e, 
		$$\epsilon_{ijk} = 
		\begin{cases}
			0 \text{, if any two indices are repeated}\\
			1 \text{, if $(ijk)$ is an even permutation}\\
			-1 \text{, if $(ijk)$ is an odd permutation},
		\end{cases}
		$$
		then by equation ( 3.29 ) of \cite{frolich}, we have 
		\begin{eqnarray*}
		 \omega \wedge \eta  &=& \sum_{i,j,k} ( \epsilon_{ijk} )^2 e_j a_j \wedge e_k  b_k = \sum_{i,j,k} ( \epsilon_{ijk} )^2 e_j \wedge e_k  a_j  b_k\\
		  &=& \sum_{j,k = 2,3} ( \epsilon_{1jk} )^2  e_j \wedge e_k a_j b_k + \sum_{j,k = 1,3} ( \epsilon_{2jk} )^2  e_j \wedge e_k a_j b_k + \sum_{j,k = 1,2} ( \epsilon_{3jk} )^2 e_j \wedge e_k a_j b_k\\
			&=& \sum_{j \neq k }  e_j \wedge e_k a_j b_k\\
			&=& \sum_{j < k }  e_j \wedge e_k ( a_j b_k - a_k b_j ).
		\end{eqnarray*}
	Therefore, we have 
		 $$ e_i \wedge e_i = 0 = e_i \wedge e_j + e_j \wedge e_i.$$
		Hence, $ \{ e_i \tensora e_i, e_i \tensora e_j + e_j \tensora e_i: 1 \leq i \leq j \leq 3 \} \subseteq {\rm Ker} ( \wedge ). $
		
		Conversely, if $ a_{ij} \in \A $ is such that $ \wedge ( \sum_{i,j} e_i \tensora e_j ) = 0, $ then by the above computation, we can conclude that $ \sum_{i < j} e_i \wedge e_j ( a_{ij} - a_{ji} ) = 0. $
		Since $ \{ e_i \wedge e_j: i < j \} $ is linearly independent, we have $ a_{ij} = a_{ji}. $
		Therefore, $ {\rm Ker} ( \wedge ) \subseteq \{ e_i \tensora e_i, e_i \tensora e_j + e_j \tensora e_i: i, j = 1,2,3 \}. $ This finishes the proof. 
		\qed	\\
		
		From the description of $ {\rm Ker} ( \wedge ) $ in Lemma \ref{30thapril20182} and the isomorphism $ \twoform \cong {\rm Span}\{ e_i \wedge e_j a_{ij} : a_{ij} = - a_{ji} \} $ (2. of Theorem \ref{30thapril20180}), it is clear that we have a right $\A$-linear splitting: $\E \tensora \E = {\rm Ker} ( \wedge ) \oplus \F$ where $\F = {\rm Span}\{ e_i \tensora e_j a_{ij}: a_{ij} = - a_{ji} \}$ is satisfied. Moreover, it is easy to verify that for all $\omega, \eta \in \Ecenter, $ the map 
		$$ \omega \tensora \eta \mapsto \frac{1}{2} ( \omega \tensora \eta + \eta \tensora \omega  ) $$
		extends to  a bilinear idempotent map on $ \E \tensora \E $ with range equal to $ {\rm Ker} ( \wedge ). $ Thus, for all $ \omega, \eta \in \Ecenter, $ we have
		$$ \Psym ( \omega \tensora \eta ) = \frac{1}{2} ( \omega \tensora \eta + \eta \tensora \omega ) $$
		and therefore, $\sigma = \sigma^{\rm can}$.
		
		\bppsn \label{26thjune2018}
		The bilinear form $g$ constructed in Subsection \ref{assumptionsonsptriple}, given by $g(e_i \tensora e_j) = \delta_{ij} 1_\A$ in the case of the fuzzy 3-spheres, is a Riemannian bilinear metric.
		\eppsn
		{\bf Proof:} From Equation (3.49) of \cite{frolich}, we see that $g:\E \tensora \E \to \A$ is defined by 
		$$g(\omega \tensora \eta)=\sum_{i=1,2,3} \omega_i \eta_i,$$
		where $\omega= \sum_{i=1,2,3} e_i \omega_i $, $\eta = \sum_{i=1,2,3} e_i \eta_i$.\\
		We need to check the conditions of Definition \ref{defn25thmay2018}. From the definition of $g$, it is clear that $g$ is an $\A$-valued map. Next, we check that the map $V_g$ is nondegenerate. Let $\omega \in \oneform$ be such that $V_g(\omega)(\eta)=0$ for all $\eta$. In particular, $g(\omega \tensora e_j)=0$ for all $j=1,2,3$. If $\omega =\sum_{i=1,2,3} e_i \omega_i$, we conclude that $\omega_i = 0$ for all $i$. Therefore, $\omega= 0$, proving that $V_g$ is one-one.\\
		Now we prove that $V_g$ is onto. Let us define $\phi_{\omega} \in \oneform^*$ by 
		$$\phi_{\omega}( e_i b)=\omega_i b$$ 
		where $\omega= \sum_{i=1,2,3} \omega_i \tensorc e_i$. Any $\phi \in {\oneform}^* $ is of the form $ \phi_\omega $ for some $\omega.$  Since $V_g(\sum_{i=1,2,3} e_i \omega_i) = \phi_{\omega}$, $V_g$ is onto.\\
		Now we prove that $g$ satisfies the equation $g \circ \sigma = g$. We have
		\begin{equation*}
		g\circ \sigma ( e_i \tensora e_j )=g( e_j \tensora e_i)
		=\delta_{ij}1_\A
		=g( e_i \tensora e_j ).
		\end{equation*}
		Since $\Ecenter = {\rm Span} \{ e_i: i=1,2,3\}$, $g \circ \sigma (\omega \tensora \eta)=g(\omega \tensora \eta)$ for all $\omega, \eta \in \oneform$.
		\qed
		
		\bthm \label{11thjuly2018}
		For every pseudo-Riemannian bilinear metric $g$ on $\E$ there exists a unique torsion-less connection which is compatible with $g$.
		\ethm
			{\bf Proof:} We need to show that the hypotheses of Theorem \ref{lcexistenceforbilinear} are satisfied. The only hypothesis we are left to verify is that $u^{\E}: \Ecenter \otimes_{\Acenter} \A \rightarrow \E $ is an isomorphism. But this is clear, since $ \Acenter = \IC.1 $ and $ \Ecenter $ is the $\IC$-linear span of $ e_1, e_2, e_3. $.
			 \qed

			The authors of \cite{frolich} investigated the existence of torsion-less and unitary connections on $\E$. While the definition of torsion of a connection discussed in their paper is the same as that in ours, the definitions of ``metric compatibility" of a connection are different, since the paper \cite{frolich} views a Riemannian metric as a sesquilinear form as opposed to a bilinear form as in our paper. In Proposition 3.7 of \cite{frolich}, it is proven that there exists a nontrivial family of torsion-less connections which are also unitary. However, once the additional condition of the connection to be real is imposed, then Corollary 3.8 of \cite{frolich} proves that such a connection is unique. We have the following result:
			
			\bthm \label{11thjuly20182}
			The Levi-Civita connection of Theorem \ref{11thjuly2018} coincides with the unique real unitary and torsion-less connection in Corollary 3.8 of \cite{frolich}.
			\ethm
			{\bf Proof:} We take basis elements $e_i$ $\oneform$, and use the fact that $e_i$ are elements of  $\Ecenter$. We denote by $\Gamma^i_{jk}$ the Christoffel symbols given by $\nabla(e_i)=\sum_{j,k}e_j \tensora e_k \Gamma^i_{jk}$. Then, we explicitly compute the metric compatibility criterion for the fuzzy 3-sphere by our definition:
			
			\begin{align*}
			0&=d(\delta_{ij})=d(g(e_i \tensora e_j))\\
			 &=(g \tensora {\rm id})({\rm id}\tensora \sigma)(\nabla(e_i)\tensora e_j + \nabla(e_j)\tensora e_i)\\
			 &=(g \tensora {\rm id})({\rm id}\tensora \sigma)(\sum_{k,l}e_k\tensora e_l\tensora e_j \Gamma^i_{kl}+\sum_{k,l}e_k\tensora e_l \tensora e_i \Gamma^j_{kl})\\
			 &=(g \tensora {\rm id})(\sum_{k,l}e_k \tensora e_j \tensora e_l \Gamma^i_{kl} + \sum_{k,l}e_k\tensora e_i \tensora e_l \Gamma^j_{kl})\\
			 &=\sum_l e_l(\Gamma^i_{jl}+\Gamma^j_{il}) {\rm, \ for \ all} \ l {\rm  \ and \ for \ all} \ i \neq j.
			\end{align*}
			Hence, the metric compatibility criterion gives us that $\Gamma^i_{jl}=-\Gamma^j_{il}$. In \cite{frolich}, combining the necessary and sufficient condition for a connection to be unitary (Equation (3.51) of \cite{frolich}) and to be a real connection, i.e. the connection coefficients must be anti-Hermitian, we get that the connection coefficients must satisfy $\Gamma^i_{jk}=-\Gamma^j_{ik}$. We see that this is the same condition that we arrive at for a metric compatible connection in our sense.
			
			The torsion-less criterion gives us that for all basis elements $e_i \in \oneform$,
			\begin{equation*}
			0=(\wedge \circ \nabla + d)(e_i)
			 =\sum_{j,k} e_j \wedge e_k\Gamma^i_{jk} - \sqrt{-1}\sum_{j,k} \epsilon^{ijk} e_j \wedge e_k,
			\end{equation*}
			where we obtain the expression for $d(e_i)$ from Equation (3.31) of \cite{frolich}.
			From Proposition 6.6 and Proposition 3.7 of \cite{frolich}, we know that this is equivalent to the criterion $\Gamma^i_{jk}-\Gamma^i_{kj}=\sqrt{-1}\epsilon^{ijk}$
			We see that the solution $\Gamma^i_{jk} = \frac{\sqrt{-1}}{2} \epsilon^{ijk}$ satisfies both the metric compatibility as well as the torsion-less criteria. Hence these are the Christoffel symbols to our unique Levi-Civita connection.
	
			 Hence, the unique real unitary and torsion-less connection in Corollary 3.8 of \cite{frolich} and the unique Levi-Civita connection for the fuzzy 3-sphere obtained by Theorem \ref{11thjuly2018} coincide. \qed
			
\section{Levi-Civita connection for quantum Heisenberg manifold} \label{heisenberguniqueness}
						
			In this section, we consider the examples of the quantum Heisenberg manifolds introduced in \cite{rieffel_heisenberg}.	In \cite{chak_sinha}, a family of spectral triples and the corresponding space of forms were studied. However, it turned out that with a particular choice of a metric and the definition of the metric compatibility of the connection 	in the sense of \cite{frolich}, there exists no connection on the space of one-forms which is both torsion-less and compatible with the metric. We will see that with our definition of metric compatibility of a connection, every pseudo-Riemannian bilinear metric on this noncommutative manifold admits a unique Levi-Civita connection. Before we begin to describe the quantum Heisenberg manifold, let us refer to the papers \cite{kang1,kang2} for the work on compatible connections for Hermitian metrics and Yang-Mills theory on the quantum Heisenberg manifolds. 
			
			 The description of the Dirac operator and the space of one-forms require the Pauli spin matrices denoted by $ \sigma_1, \sigma_2, \sigma_3 $ in \cite{chak_sinha}. In particular, the $\sigma_i$'s satisfy the following relations: 
			 \be \label{paulispin} \sigma^2_j = 1, \sigma_j \sigma_k = - \sigma_k \sigma_j, ~ \sigma_1 \sigma_2 = \sqrt{- 1} \sigma_3, \sigma_2 \sigma_3 = \sqrt{- 1} \sigma_1, \sigma_1 \sigma_3 = \sqrt{- 1} \sigma_2.\ee
			 Moreover, we are going to work with right connections instead of left connections as done in \cite{chak_sinha}.\\

			The description of the quantum Heisenberg manifold in \cite{rieffel_heisenberg} is as follows.
			\begin{dfn}
				For any positive integer $c$, let $S^c$ denote the space of infinitely differentiable functions $\Phi: \IR \times \IT \times \IZ \to \IC$ that satisfy the following two conditions:
				\begin{itemize}
					\item[(a)] $\Phi(x+k,y,p)=e^{2\pi i ckpy}\Phi(x,y,p)$ for all $k \in \IZ$,
					\item[(b)] for every partial differential operator $\tilde{X}=\frac{\partial^{m+n}}{\partial x^m \partial y^n}$ on $\IR \times \IT$ and every polynomial
				\end{itemize} $p$ on $\IZ$, the function $P(p)(\tilde{X}\Phi)(x,y,p)$ is bounded on $K \times \IZ$ for any compact set $k$ of $\IR \times \IT$.\\
				For each $\hbar,\mu,\nu \in \IR$ with $\mu^2 + \nu^2 \neq 0$, let $\A^\infty_\hbar$ denote the space $S^c$ equipped with product and involution defined, respectively, by
				\begin{itemize}
					\item[] $(\Phi \star \Psi)(x,y,p)$
							$=\sum_q \Phi(x - \hbar(q-p)\mu,y-\hbar(q-p)\nu,q)\Psi(x-\hbar q \mu, y - \hbar q ,\nu, p-q)$,
					\item[] $\Phi^*(x,y,p)=\Phi(x,y,-p)$.
				\end{itemize}
				Let $\pi$ be the representation of $\A^\infty_\hbar$ on $L^2(\IR \times \IT \times \IZ)$ given by
				$$ (\pi(\Phi)\xi)(x,y,p) = \sum_q \Phi(x- \hbar (q-2p)\mu, y - \hbar (q-2p)\nu, q)\xi(x, y, p-q). $$
				Then $\pi$ gives a faithful representation of the involutive algebra $\A^\infty_\hbar$. The norm closure of $\pi(\A^\infty_\hbar)$, denoted by $\A^{c,\hbar}_{\mu,\nu}$ is called the quantum Heisenberg manifold.
			\end{dfn}
			 For the rest of this section, we will denote the involutive algebra $\A^\infty_\hbar$ by $\A$. The algebra $ \cla $ admits an action of the Heisenberg group. The symbol $ \tau $ will denote a certain state on $ \cla $ invariant under the action of the Heisenberg group. Let $ X_1, X_2, X_3 $ denote the canonical basis of the Lie algebra of the Heisenberg group so that we have associated self-adjoint operators $ d_{X_i} $ on $ L^2 ( \cla, \tau ) $ in the natural way. Then the triple $ ( \cla, L^2 ( \cla, \tau ) \tensorc \IC^2, D ) $ defines a spectral triple on $ \cla $ where $ \cla $ is represented on $ L^2 ( \cla, \tau ) \tensorc \IC^2 $ diagonally and the Dirac operator $ D $ is defined as 
			$$ D = \sum_j d_{X_j} \tensorc \sigma_j, $$
			where $ \{ \sigma_j: j = 1,2,3\} $ are the self-adjoint Pauli spin matrices satisfying \eqref{paulispin}.
			
			Let us denote the operator $1 \tensorc \sigma_i$ by the symbol $e_i$. Then, the following lemma is a direct consequence of the proof of Proposition 9 of \cite{chak_sinha}.

\blmma \label{27thjan1}
For all $ a $ in $ \cla, $
$$ d ( a ) = \sum^3_{j = 1} e_j \partial_j ( a ), $$
 $$ {\rm where} ~ \partial_1 ( a ) = \frac{\partial a}{\partial x}, ~ \partial_2 ( a ) = - 2 \pi \sqrt{- 1} c p x a + \frac{\partial a}{\partial y}, ~ \partial_3 ( a ) = - 2 \pi \sqrt{- 1} c p \alpha a $$
for some $ \alpha $ greater than $ 1. $
The derivations $ \partial_1, \partial_2, \partial_3 $ satisfy the following relation:
			 \be \label{delirelation} [ \partial_1, \partial_3 ] = [ \partial_2, \partial_3 ] = 0, ~ [ \partial_1, \partial_2 ] = \partial_3. \ee
\elmma

The space of one-forms and two-forms for the spectral triple $(\A,L^2(\A,\tau)\tensorc \IC^2,D)$ are as follows:

\bppsn \label{oneandtwoform}
The module of one-forms $ \E:= \oneform $ is a free module generated by $ e_1, e_2, e_3. $ Moreover, $ e_1, e_2, e_3 $ are central elements.
As a subset of $ B ( L^2(\A,\tau) \tensorc \IC^2 ),$ $ \E $ can be described as follows:
   $$ \E = \{ \sum_i  a_i \tensorc \sigma_i : a_i \in \A  \}. $$
The set of  junk forms ( cf. Subsection \ref{sptripleforms} ) is equal to $ \{ a \tensorc 1 : a \in \A \}, $ and therefore is isomorphic to $\A.$ 
Finally, the space of two forms $ \twoform $ is isomorphic to $  \cla \oplus \cla \oplus \cla. $ More precisely, 
 \[\twoform = \{  a_1 \tensorc \sigma_1 \sigma_2 + a_2 \tensorc \sigma_2 \sigma_3 + a_3 \tensorc \sigma_1 \sigma_3 : a_1, a_2, a_3 \in \A \} \subseteq B ( L^2(\A,\tau) \tensorc \IC^2  ).\]
\eppsn
{\bf Proof:} The space of one-forms is described in Proposition 21 of \cite{chak_sinha}. The fact that $ e_1, e_2, e_3 $ are central can be easily seen from the definition of the representation of $ \cla $ on $ L^2 ( \cla, \tau ) \tensorc \IC^2. $ The statement about the two forms follow from Proposition 22 of the same paper. \qed

\bppsn \label{16thoctober20192sm}
The bilinear form $g$ constructed in Subsection \ref{assumptionsonsptriple} satisfies the conditions of Definition \ref{defn25thmay2018}, i.e, it is the canonical Riemannian bilinear metric for the spectral triple.
\eppsn
 {\bf Proof:} We need to check the conditions of Definition \ref{defn25thmay2018}. This essentially follows from the results of \cite{chak_sinha}. We will use Proposition \ref{oneandtwoform} to identify $\E$ with $\A \tensorc \IC^3, $ the bimodule structure being defined as:
$$ a ( e_i b ) c = e_i a b c . $$
We will let $\tau$ denote the functional on $ \mathcal{B} ( \mathcal{H} ) $ as in Subsection \ref{assumptionsonsptriple}. Let $\psi: \A \rightarrow \IC $ be the faithful normal tracial state on $\A^{\prime \prime}$ as in Section 2 of \cite{chak_sinha} ( denoted by $\tau$ in \cite{chak_sinha} ). By Proposition 14 of \cite{chak_sinha},
$$ \tau ( X ) = ( \frac{1}{2} \psi \tensorc {\rm Tr} ) ( X ) ~{\rm for} ~ {\rm all} ~ X ~ \in ~ \oneform. $$ 
since $ \psi $ is faithful on $ \A^{\prime \prime}, $ we can conclude that $ \tau $ is faithful on the $ \ast $-algebra generated by $\A$ and $ \{ [ D, a ]: a \in \A \}. $ Moreover, by identifying $ \A \subseteq \oneform = \A \tensorc \IC^3 $ via $ a \mapsto a \otimes I_2, $ $ \tau ( a ) = \psi ( a ) $ for all $ a $ in $\A. $

If $ \omega = \sum^3_{i = 1} e_i a_i  $ and $ \eta = \sum^3_{i = 1}e_i b_i  $ are two one-forms, then 
$$ \frac{1}{2} ( I \tensorc {\rm Tr} ) ( \omega \eta ) = \sum^3_{i = 1} a_i b_i. $$
Therefore, for all $ c $ in $\A,$ the formula $ g ( \omega \tensora \eta ) = \left\langle \left\langle \omega^*, ~ \eta \right\rangle\right\rangle $ ( Lemma \ref{22ndjune2018}  ) implies that
\begin{align*}
 \tau ( g ( \omega \tensora \eta ) c ) &= \tau ( \left\langle \left\langle \omega^*, ~ \eta \right\rangle\right\rangle c  )\\
                    &= \tau ( \omega \eta c ) \\
				    &= ( \frac{1}{2} \psi \tensorc {\rm Tr} )  ) ( \omega \eta c )\\
					&= \sum_i \psi ( a_i b_i c )\\
					&= \tau ( (  \sum^3_{i=1} a_i b_i  )  c ).
\end{align*}
Therefore, $ g ( \omega \tensora \eta ) = \sum^3_{i = 1} a_i b_i \in \A. $
The nondegeneracy of the map $ V_g $ follows just as in Proposition \ref{26thjune2018}. 
 \qed

\bppsn \label{16thoctober20191sm}
Let $\F$ denote the right $\A$-linear span of the set $ \{ e_i \tensora e_j - e_j \tensora e_i: 1 \leq i < j \leq 3  \} $. Then, the bimodule $\E \tensora \E$ admits a decomposition $\E \tensora \E = {\rm Ker}(\wedge) \oplus \F$ as right $\A$-modules. Moreover, the map $\sigma = 2 \Psym - 1$ is equal to the map $\sigma^{\rm can}$ as in Theorem \ref{skeide2}, i.e. for all $e$, $f$ in $\Ecenter$, and $a$ in $\A$,
\[ \sigma(e \tensora f a) = f \tensora e a. \]
\eppsn
{\bf Proof:} We will use the fact that $ e_i $ are central elements throughout the proof. Moreover, let $ \wedge, m_0, \mathcal{J}$, be as in Subsection \ref{sptripleforms} while $\Psym$ will be as in Definition \ref{defn24thmay2018}. By the description of $\mathcal{J}$ and that of $\twoform$ in Proposition \ref{oneandtwoform}, it is easy to see that $ {\rm Ker} ( \wedge ) $ is spanned by $ \{ e_i \tensora e_j + e_j \tensora e_i : 1 \leq i \leq j \leq 3 \} $ and $ \clf = Span \{ e_i \tensora e_j - e_j \tensora e_i: 1 \leq i < j \leq 3  \}. $ Clearly, $ \E \tensora \E = {\rm Ker} ( \wedge ) \oplus \clf $ as right $\A$ modules.

Since $e_1,e_2,e_3 \in \Ecenter$, it can be easily checked that $u^{\E}$ is an isomorphism. In particular, $\E$ is centered. Moreover, by the description of ${\rm Ker}(\wedge) $ as above, we have
$$ P_{{\rm sym}} ( e_i \tensora e_j - e_j \tensora e_i ) = 0, ~ P_{{\rm sym}} ( e_i \tensora e_j + e_j \tensora e_i ) = e_i \tensora e_j + e_j \tensora e_i $$
$$ {\rm and} ~ {\rm thus} ~ 2 P_{{\rm sym}} ( e_i \tensora e_j ) = e_i \tensora e_j + e_j \tensora e_i. $$
$${\rm Therefore,} ~ \sigma ( e_i \tensora e_j ) = ( 2 P_{{\rm sym}} - 1 ) ( e_i \tensora e_j ) = e_j \tensora e_i. $$
Therefore, $\sigma=\sigma^{\rm can}$. 
\qed

\bthm \label{existenceheisenberg}
 For any pseudo-Riemannian bilinear metric $g$ on $ \E, $ there exists a unique Levi-Civita connection on the module $ \E $ compatible with $g$.
\ethm
{\bf Proof:} In Proposition \ref{oneandtwoform} and Proposition \ref{16thoctober20191sm}, we have verified all the assumptions of Theorem \ref{lcexistenceforbilinear}. Hence we have our result. \qed

\brmrk
In \cite{article2},  we have given the explicit computation of the Levi-Civita connection and associated curvature for the Levi-Civita connection with respect to the canonical Riemannian metric $g$, where it has been shown that the scalar curvature is a negative constant multiple of the identity element of $\A.$ 
\ermrk

\section{Levi-Civita connection for Connes-Landi deformed spectral triples}
  \label{prelimdef}
		
		Suppose $M$ is a compact Riemannian manifold such that the maximal torus of the isometry group of $M$ has rank greater than or equal to 2. Then the action of the maximal torus on $C^{\infty}(M)$ allows us to define a deformed algebra $ C^{\infty}(M)_{\theta}$ (\cite{rieffel}, \cite{Connes-dubois}). Moreover, the torus equivariant spectral triple on $M$ can be deformed to a new spectral triple on $ C^{\infty}(M)_{\theta}$(\cite{connes_landi}).	
		The goal of this section is to prove the following theorem:
		\bthm \label{existencerieffeldeformation}
		Suppose $M$ is a compact Riemannian manifold equipped with a free isometric action of $\IT^n$. Let $\E:=\oneformclassical$ denote the space of one-forms of the spectral triple $(C^\infty(M),\bigoplus_kL^2(\Omega^k(M),d+d^{*})$ and let $\E_\theta$ be the deformation of $\E$ as in Subsection \ref{rieffelgen}. Then for any Riemannian bilinear metric $g$ on $\E_\theta$ there exists a unique Levi-Civita connection on the bimodule $\E_\theta$.
		\ethm
		
		 In the first subsection, we prove some preparatory results on the fixed point algebra under the action of a compact abelian Lie group. In Subsection \ref{rieffelgen} we prove some results on generalities of Rieffel deformations. In Subsection \ref{bimetricexistencedef} we prove that there exists a Riemannian bilinear metric on $\E_\theta$ and that it is the deformation of the canonical metric on $\E$. In Subsection \ref{rieffelsubsection}, we prove that under our assumptions, the deformed module of one-forms on the Rieffel deformed manifold satisfies the conditions of Theorem \ref{lcexistenceforbilinear}. 
			 
		We recall that for any action $  \beta $ of $\mathbb{T}^n$ on a module $\mathcal{G}$ (or an algebra $\mathcal{D}$), the spectral subspace corresponding to
a character $\underline{m}\equiv(m_{1},...,m_{n})\in\widehat{\mathbb{T}^n}\cong \mathbb{Z}^n$, denoted by $\mathcal{G}_{\underline{m}}$ ( respectively $ \mathcal{D}_{\underline{m}} $ ),
consists of all $\xi$ such that $  \beta_{t}(\xi)=\chi_{\underline{m}}(t)\xi$ for all $t=(t_1,\ldots, t_n) \in\mathbb{T}^n$, where $\chi_{\underline{m}}(t):=
t_1^{m_1} \ldots t_n^{m_n}.$
It is easily seen that $\mathcal{D}_{\underline {m}} {\mathcal D}_{\underline{n}} \subseteq {\mathcal D}_{\underline{m}+\underline{n}}.$

 Suppose that $\mathcal{G}$ is a $\mathcal{D}$ bimodule. Moreover, let us assume that both $\mathcal{D}$ and $\mathcal{G}$ are equipped with actions of $ \IT^n $ in such a way that $ \mathcal{G}$ becomes an equivariant $\mathcal{D}$ bimodule. Then we have $\mathcal{G}_{\underline {m}} {\mathcal D}_{\underline{n}} \subseteq \G_{\underline{m}+\underline{n}}$ and $ \D_{\underline {n}} {\mathcal G}_{\underline{m}} \subseteq \G_{\underline{m}+\underline{n}} .$ 
The subspace  $ {\rm Span} \{ \D_{\underline{m}}: \underline{m} \in \IZ^n \} $ comprise the so-called `spectral subalgebra' for the action. Similarly, $ {\rm Span} \{\G_{\underline{m}}: \underline{m} \in \IZ^n \} $ is called the spectral submodule of the action.

		Let $ G $ be a group. Let us recall that a spectral triple $(\mathcal{A},\mathcal{H},D)$ is called $ G $ equivariant if there exists a unitary representation $ \beta $ of $ G $ on $ \clh $ such that $ \beta_g D = D \beta_g. $ Moreover, 	we recall	 the following well known fact ( see \cite{Connes-dubois} for the details ).
			
			\bppsn \label{equivarianttriple}
			Suppose that $ M $ is a compact Riemannian manifold with an isometric action of the $ n $-torus $ \IT^n $ on $ M. $ Consider the spectral triple $ ( C^\infty ( M ), \mathcal{H}, d + d^* ) $ where $ \clh $ is the Hilbert space of forms and $ d $ is the de-Rham differential on $ \mathcal{H}. $ The $ \IT^n $action on smooth forms extends to a unitary representation $\beta$ on $ \clh $ and the spectral triple is equivariant w.r.t this representation of $ \IT^n. $  In particular, if $ \alpha $ denotes the action of $ \IT^n $ on $ C^\infty ( M ) $ and $ \delta ( \cdot ) = \sqrt{-1}[ d + d^*, \cdot ], $ then 
			$$ \beta_t ( f \delta ( g ) ) = \alpha_t ( f ) \beta_t ( \delta ( g ) ) = \alpha_t ( f ) \delta ( \alpha_t ( g ) ) ~ \forall ~ t ~ \in \IT^n. $$
			\eppsn

		In this set up, it is easy to see the following result:
		
		\blmma \label{mequivariant}
		If $ \D $ is a subalgebra of $ C^\infty ( M ) $ kept invariant by the action of a compact group $ G $ acting by isometries on $ M, $ then the map $ \wedge: \Omega^1 ( \D ) \otimes_{\D} \Omega^1 ( \D ) \rightarrow \Omega^2 ( \D ) $ is $ G $ equivariant. 
		\elmma
		
		As an immediate corollary, we have
		
		\bcrlre \label{Finvariant}
		With the notations of Lemma \ref{mequivariant}, $ {\rm Ker} ( \wedge ) $ is invariant under the action of $ G .$ Moreover, if $ \Omega^1 ( \D ) \otimes_{D} \Omega^1 ( \D ) = {\rm Ker} ( \wedge ) \oplus \mathcal{G} $ gives the decomposition as in Subsection \ref{splittingsubsection}, then $ \mathcal{G} $ is also kept invariant by $ G. $
		\ecrlre
		{\bf Proof:} The $ G $ invariance of $ {\rm Ker} ( \wedge ) $ follows from the $ G $ invariance of $ \wedge. $ Moreover, we have $ \mathcal{G}  = {\rm Ker} ( 1 - \sigma^{{\rm can}} ) .$ Since $ \sigma^{{\rm can}} $ is $ G $ equivariant, $ \G $ is $ G $ invariant. \qed\\

			\subsection{Some results on the fixed point algebra} \label{classicalcase}
			
		Let us consider a compact Riemannian manifold $ M $ with the $ \IT^n $ equivariant spectral triple $(C^\infty(M),\mathcal{H},d+d^{*})$ as in Proposition \ref{equivarianttriple}. Throughout this section, we will follow the notations introduced in the following definition.
		
		\bdfn \label{defn22ndfeb2017}
		  Let $ \E := \oneformclassical$ and $\A:=C^\infty(M)$. $\F$ will denote the $ \IT^n $ equivariant spectral submodule of $ \E $. The symbol $ \F_{\underline{k}} $ will denote the $ \underline{k} $-th spectral subspace of $ \F. $ Thus, $ \F = {\rm Span} \{ \F_{\underline{k}}: \underline{k} ~ \in ~ \IZ^n \}. $ Similarly, we define $ \clc $ to be the spectral subalgebra $ {\rm Span} \{ \clc_{\underline{k}}: \underline{k} \in \IZ^n \}$ of $\A$ where $ \clc_{\underline{k}} $ is the $ \underline{k} $-th spectral subspace of $ \C. $ In particular, $\C_{\underline{0}}$ and $\F_{\underline{0}}$ denote the $\IT^n$ invariant spectral subalgebra and the $\IT^n$ invariant spectral submodule respectively.			 
 		\edfn
		
		\brmrk
		It is clear from the definition of spectral subspaces of algebras and modules that if $\A_{\underline{k}}$ and $\E_{\underline{k}}$ denote the spectral subspaces of $\A$ and $\E$ respectively, then $\A_{\underline{k}} = \C_{\underline{k}}$ and $\E_{\underline{k}} = \F_{\underline{k}}$. We will from now on use this fact, often without mentioning.
		\ermrk		
		
		\brmrk \label{prop616}
		Since the representation $\beta$ as in Proposition \ref{equivarianttriple} commutes with $ d + d^*, $ it is easy to see that $ \beta_t ( \F ) \subseteq \F $ for all $ t \in \IT^n. $ Moreover, it is easy to see that the space of one-forms for the spectral triple $ ( \clc, \clh, d + d^* ) $ is precisely $ \F. $
		\ermrk

		The aim of this subsection is to prove that if the action of $\IT^n$ on $M$ is free, then the spectral subalgebra $\C_{\underline{0}}$ and the spectral submodule $\F_{\underline{0}}$ satisfy the hypotheses of Proposition \ref{isomorphismofEprime}.

\blmma
\label{free_case}
Suppose that the $T \equiv \IT^n$ action on $M$ is free. Then $ \F_{\underline{0}} $ is a finitely generated projective right module over $\C_{\underline{0}}.$ 
\elmma
{\it Proof :} For a module $ \G $ equipped with an action of $ \IT^n, $ let us denote the $ \IT^n $ invariant submodule of $ \G $ by the symbol $ \G^{\IT^n}. $ 
Since the $\IT^n$-action on $ M $ is free, $M/\IT^n$ is a smooth compact manifold and $M$ is a principal 
$\IT^n$-bundle over $M/\IT^n$. Let $\pi$ denote the projection map from $ M$ onto $M/\IT^n$. 
Given any point in $M$, we can find a $\IT^n$ invariant open neighborhood $U$ which is $\IT^n$ equivariantly diffeomorphic 
with $U/\IT^n \times \IT^n.$ Moreover, we can choose $U$ in such a way that $U/\IT^n$ is the domain of a local coordinate chart for the manifold $M/\IT^n$, say $U=\pi^{-1}(V)$,
where $V$ is the domain of some local chart for $M/\IT^n$. This gives the following isomorphism: 
$$\Omega^1(U)^{\IT^n} \cong \Omega^1(U/\IT^n) \otimes_\IC \Omega^1(\IT^n)^{\IT^n} \cong \Omega^1(U/T) \otimes_\IC \cll,$$ 
$\cll$ being the complexified Lie algebra 
of $\IT^n$ which is nothing but $\IC^n$. As $U/\IT^n$ is the domain of a local coordinate chart, the module of one-forms is a free $C^\infty(U/\IT^n)$ module, 
say $C^\infty(U/\IT^n) \otimes_{\IC}
\IC^k$, hence 
$\Omega^1(U/\IT^n)$ is isomorphic with $C^\infty(U/\IT^n) \tensorc \IC^{n+k}$, i.e. $\Omega^1(U/\IT^n)$ is free. By covering $M$ with finitely many such neighborhoods, we can complete the proof of 
$C^\infty(M/\IT^n)$ projectivity and finite generation of $\Omega^1(M)^{\IT^n}$. \qed\\

		Now, we will make use of the notation $ u^\F_{\F_{\underline{0}}} : \F_{\underline{0}} \otimes_{\clc_{\underline{0}}} \clc \rightarrow \F $ introduced in Subsection \ref{isomorphismremark}.

\blmma
 \label{lemma_for_assump_iv}
 If for each $\underline{m} \in {\mathbb Z}^n$, we can find $a_1, \ldots, a_k \in {\mathcal C}_{\underline{m}}$ and $ b_1, \ldots b_k \in {\mathcal C}_{ - \underline{m}}$
 ($k$ depends on $\underline{m}$ ) such that $\sum_i b_i a_i=1,$  then 
  the map $u^\F_{\F_{\underline{0}}}$ is an isomorphism.
\elmma
{\bf Proof:} We need to prove that under the above assumption, the map $ u^\F_{\F_{\underline{0}}} $ has a right $ \clc $-linear inverse. However, since $ u^\F_{\F_{\underline{0}}} $ is right $ \clc $-linear to start with, it suffices to prove that $ u^\F_{\F_{\underline{0}}} $ defines an isomorphism of vector spaces. Hence, it is sufficient to prove that for all  $ {\underline{m}}$, the restriction $ p^\F_{\underline{m}} $ of $u^\F_{\F_{\underline{0}}}$ to ${\mathcal F}_{\underline{0}} \otimes_{{\mathcal C}_{\underline{0}}} {\mathcal C}_{\underline{m}}$
 is a vector space isomorphism onto its image $ {\mathcal F}_{\underline{m}}$.

Then the map 
$$q^\F_{\underline{m}} : {\mathcal F}_{\underline{m}} \rightarrow {\mathcal F}_{\underline{0}} \otimes_{{\mathcal C}_{\underline{0}}} {\mathcal C}_{\underline{m}} ~ {\rm defined} ~ {\rm by} ~ q^\F_{\underline{m}}(e):=\sum_i e b_i \otimes_{{\mathcal C}_{\underline{0}}} a_i$$
satisfies $p^\F_{\underline{m}} \circ q^\F_{\underline{m}}={\rm id}.$ On the other hand,
 as $b_i a \in {\mathcal C}_{\underline{0}}$ if $a \in {\mathcal C}_{\underline{m}}$, we have 
	$$q^\F_{\underline{m}} \circ p^\F_{\underline{m}} (e \otimes_{{\mathcal C}_{\underline{0}}} a)
 =\sum_i ea b_i \otimes_{{\mathcal C}_{\underline{0}}} a_i =e \otimes_{{\mathcal C}_{\underline{0}}} \sum_i a b_i a_i =e \otimes_{{\mathcal C}_{\underline{0}}} a.$$
	This finishes the proof of the lemma.\qed

Now we shall identify the bimodule $ \F_{\underline{m}} $ with the bimodule of sections of a certain vector bundle over $M$.

\blmma
\label{sec}
 Let $ M $ be a smooth compact Riemannian manifold equipped with a smooth and free right action of a compact connected abelian Lie group $ K. $
	Let $ M \times_{\chi_{ - \underline{m}}} {\mathbb C} \rightarrow M/K $ denote the associated vector bundle ( of $ M \rightarrow M /K $ ) corresponding to the character $ \chi_{-\underline{m}}. $
			
			Then the space of all smooth functions $ f $ on $ M $ satisfying $ f ( x.t ) = \chi_{\underline{m}} ( t ) f ( x ) $ is in one to one correspondence with the set of all smooth sections of the vector bundle $ M \times_{\chi_{ - \underline{m}}} {\mathbb C} \rightarrow M / K. $  
\elmma

{\bf Proof:}
The elements of the total space of the associated vector bundle $M \times_{\chi_{ - \underline{m}}} {\mathbb C}$ are given by the equivalence class $[y, \lambda]$ of $(y,\lambda)\in M\times \mathbb{C}$ such that
$(y,\lambda) \sim (y.t,\chi_{ - \underline{m}}(t^{-1})\lambda)$ for all $t\in K.$ Now, for $ f \in \E_{\underline{m}}, $ 
we can define a section of the above vector bundle $s_f$ by
 $$s_f([x])= [x,f(x)],$$
 where $[x]$ denotes the class of the point $x$ in $M/K. $
 We need to show that this is well defined. But for any $t \in K,$ $s_f([x.t])=[x.t,f(x.t)]= [x.t,\chi_{\underline{m}}(t)f(x)]=[x.t,\chi_{ - \underline{m}}(t^{-1})f(x)]
 =[x,f(x)].$ This proves that $s_f$ is well defined. Similarly, given a section $s$ of the above vector bundle we can define a
 function $f_s$ on $M$ by $f_s(x)=\lambda_x$ where $\lambda_x \in {\mathbb C}$ is such that $s([x])=[x,\lambda_x]$. Clearly, $\lambda_x$ is uniquely
 determined, because the $K$ action is free. Moreover, $[x,\lambda_x]=[x.t,\chi_{ - \underline{m}}(t)^{-1}\lambda_x] $implies $\lambda_{xt}=\chi_{ - \underline{m}}(t)\lambda_x$, i.e. $f_s \in \E_{\underline{m}}$.

  Finally, it is easy to verify that the maps $f \mapsto s_f$ and $s \mapsto f_s$ are inverses of one another, completing the proof. \qed
	
	The following lemma is well-known. However, we give a proof for it since we could not find any appropriate references.
	
\blmma
\label{class}
For a complex smooth Hermitian vector bundle over a compact manifold $M$ there are finitely many smooth
sections $s_{i}$'s such that $\sum_i \left\langle \left\langle s_{i},~ s_{i} \right\rangle\right\rangle = 1$ where $ \left\langle \left\langle \cdot , \cdot \right\rangle \right\rangle$ denotes the $C^{\infty}(M)$- valued inner product coming from the Hermitian structure.
\elmma

{\bf Proof:} Corresponding to a finite open cover $\{ U_{i},~i=1, \ldots, l \}$ choose finitely many smooth sections $\gamma_{i}$ which are non zero on $U_{i}$. Then choosing a smooth
partition of unity $ \psi_i, i=1, \ldots, l$, we can construct $t_{i}=\psi_i \gamma_i$'s so that $t=\sum \langle \langle t_{i},t_{i} \rangle \rangle $ is nowhere zero.
The sections $s_i=\frac{t_i}{t^{\frac{1}{2}}}$ satisfy the conditions of the lemma.\qed\\
This gives us the following:

\blmma
 \label{classical_iso}
 Suppose $ M $ is a compact Riemannian manifold equipped with a free and isometric action of $ \IT^n $. Then the map $ u^\F_{\F_{\underline{0}}} $ is an isomorphism. 
\elmma

{\bf Proof :} Without loss of generality, we can assume $M$ to be connected. In general, if $M$ has $k$ connected components $M_1, M_2, \cdots M_k,$ the module $\F$ decomposes as $\F_1 \oplus \cdots \F_k,$ where $\F_i$ is the linear span of spectral subspaces of $\Omega^1 ( M_i ), $ and it is suffices to prove that for all $i,$ $u^{\F}_{\F_{\underline{0}}}$ is an isomorphism from $ ( \F_i )_{\underline{0}} \otimes_{( \clc_i )_{\underline{0}}} \clc_i $ onto $ \F_i. $ Since the action of $ \IT^n $ on $ M $ is free, $ M \rightarrow M/ \IT^n $ is a principal $ \IT^n $-bundle. Therefore, for the sections $ s_i $ as in Lemma \ref{class}, we have functions $ f_{s_i} $ in $ \clc_{\underline{m}} $ by Lemma \ref{sec}. By the definition of $ f_{s_i} $ and the relation $\sum \langle \langle s_{i},s_{i} \rangle \rangle =1,$ it follows that 
$$ \sum_i \overline{f_{s_i}}f_{s_i}=1.$$
Since $ f_{s_i} $ belongs to $ \clc_{\underline{m}} ,$ the function $ \overline{f_{s_i}} $ belongs to $ \clc_{ - \underline{m}}. $ Thus, we can apply Lemma \ref{lemma_for_assump_iv} to deduce the conclusion of the theorem. \qed \\

Next we prove that with the hypothesis of Lemma \ref{classical_iso}, the map $u^\E_{\E_{\underline{0}}}$ is also an isomorphism.

\blmma \label{isomorphismofEzero}
The map $u^{\E}_{\mathcal{E}_{\underline{0}}}:\E_{\underline{0}}\otimes_{\A_{\underline{0}}}\A \rightarrow \E$ is an isomorphism.
\elmma
{\bf Proof:}
Let us start by proving that the map is one-to-one. Let $u^{\E}_{\mathcal{E}_{\underline{0}}}(\sum_ie'_i\otimes_{\A_{\underline{0}}}f_i)=0$, i.e. $\sum_ie'_if_i=0$. Then each spectral projection $\mathcal{P}_{\underline{m}}(\sum_ie'_if_i)=0 $ i.e, $\sum_ie'_i\mathcal{P}_{\underline{m}}(f_i)=0$. So, $u^{\E}_{\mathcal{E}_{\underline{0}}}(\sum_ie'_i\otimes_{\mathcal{A}_{\underline{0}}}\mathcal{P}_{\underline{m}}(f_i))=0$ for all $m$. But $u^{\E}_{\mathcal{E}_{\underline{0}}}$ restricted to each spectral subspace is as isomorphism as already proved in Lemma \ref{classical_iso}. Hence, $\sum_ie'_i\otimes_{\mathcal{A}_{\underline{0}}}\mathcal{P}_{\underline{m}}(f_i)=0$ for all $m$. Thus $\sum_ie'_i\otimes_{\A_{\underline{0}}}f_i=\underset{N}{\rm lim}(\sum_{i,|m|\leq N}e'_i\otimes_{\mathcal{A}_{\underline{0}}}\mathcal{P}_{\underline{m}}(f_i))=0$, where ${\rm lim}$ denotes the limit in the Fre\'chet topology. Therefore, the map is one-to-one.

Now we show that the map is onto. Since the map $u^\E_{\E_{\underline{0}}}$ is right $\A$-linear, it suffices to check that for all $f\in \A$, $df$ has a pre-image in $\mathcal{E}_{\underline{0}}\otimes_{\A_{\underline{0}}}\A$. Consider the principal $T=\mathbb{T}^{n}$ bundle $\pi:M\rightarrow M/T$. Since $M/T$ is compact, we can take a finite atlas $(U_i,\phi_i)$ on it such that the bundle $\pi^{-1}(U_i)\rightarrow U_i$ is $T$ equivariantly diffeomorphic with the canonical bundle $U_i\times T\rightarrow U_i$. Let $\{\psi_i\}_i$ be a partition of unity on $M$ subordinate to $(U_i,\phi_i)$. Then $f=\sum_i f\psi_i$ and $df=\sum_i d(f\psi_i)$. Thus in particular we can assume $f$ is supported in $\pi^{-1}(U_i)$ or equivalently in $U\times T$. Let $\{dx_i\}$ be a basis for differential forms along the direction of $U$ i.e. the horizontal direction of the bundle $U\times T \rightarrow U$ and $\{\omega_j\}$ be a basis of right invariant 1-forms in the vertical direction corresponding to the basis $\{\chi_j\}$ of right invariant vector fields along the direction of $T$. Then $df=\sum_i dx_i.\frac{\partial}{\partial x_i}(f)+\sum_j\omega_j.\chi_j(f)$. The right action of $T$ on $U\times T$ acts trivially in the direction of U, hence $dx_i\in \mathcal{E}_{\underline{0}}$. Since $\omega_j$ is invariant under the action induced by the right action of $T$ on $U\times T$, so $\omega_j \in \mathcal{E}_{\underline{0}}$. Hence $df$ has a pre-image $\sum_i dx_i\otimes_{\mathcal{A}_{\underline{0}}}\frac{\partial}{\partial x_i}(f)+\sum_j\omega_j\otimes_{\mathcal{A}_{\underline{0}}}\chi_j(f)\in \mathcal{E}_{\underline{0}}\otimes_{\mathcal{A}_{\underline{0}}}\A$. Therefore, we have that $u^\E_{\E_{\underline{0}}}$ is an onto map.
This completes the proof. \qed

\subsection{Some generalities on Rieffel-deformation}	\label{rieffelgen}

Our main reference for Rieffel deformation of a $ C^* $ algebra endowed with a strongly continuous action by $ \mathbb{T}^n $ is \cite{rieffel}. However, we will also need to use equivalent descriptions of this deformation given in
		 \cite{Connes-dubois}, \cite{connes_landi}, \cite{scalar_5} and \cite{hanfeng}.

We next define the deformation of an algebra $\D$. We refer to \cite{scalar_5} for details.

\bdfn[Definition 2.3 \cite{scalar_5}]
Let $\mathcal D$ be a $\IT^n$ smooth algebra as in Definition 2.2 of \cite{scalar_5}. For a skew symmetric $n\times n$ matrix $\theta$, consider the  bi-character $\chi_\theta$ defined by
$$\chi_\theta(\underline{k},\underline{l})=e^{\pi i \langle k,\theta l \rangle}, \ \underline{k},\underline{l}\in \IZ^n,$$
where the pairing $\langle .,. \rangle$ is the usual dot product in $\mathbb{R}^n$. The deformation of $\D$ is the algebra $\D_\theta$ whose underlying vector space is equal to $\D$ while the multiplication $\times_\theta$ is deformed as follows:
\begin{equation}
a \times_\theta b = \sum_{\underline{k},\underline{l}\in \IZ^n} \chi_\theta(\underline{k},\underline{l}) a_{\underline{k}} b_{\underline{l}}, \ \forall \ a,b \in \D, \label{deformationformula}
\end{equation}
where $a=\sum_{\underline{k}}a_{\underline{k}},\ b=\sum_{\underline{l}}b_{\underline{k}} $ are the isotypical decompositions.
\edfn

\brmrk \label{convergenceisotypical}
By Proposition 2.1  of \cite{scalar_5}, the isotypical decompositions converge absolutely to the element.
\ermrk

$\D_\theta$ turns out to be a $\IT^n$ smooth algebra and the deformed product is associative.

Similarly, one can deform $\IT^n$ smooth $\D$ bimodules (refer to Definition 2.2 of \cite{scalar_5}) as follows:

\bdfn
Let $\G$ be a $\IT^n$ smooth $\D$-$\D$ bimodule. Then the deformed bimodule $\G_\theta$ is a $\D_\theta$-$\D_\theta$ bimodule whose underlying vector space is equal to $\G$ while the deformed right module action is as follows
\begin{equation}
e\times_\theta a= \sum_{\underline{k},\underline{l}\in \IZ^n} \chi_\theta(\underline{k},\underline{l}) e_{\underline{k}} a_{\underline{l}}, \ \forall \ e \in \G, \ \forall \ a \in \D, \label{moduledeformationformula}
\end{equation}
where $e = \sum_{\underline{k}}e_{\underline{k}}$ and $a = \sum_{\underline{l}}a_{\underline{l}}$ are the isotypical decompositions, and the deformed left module action is defined similarly.
\edfn

Using the fact that $\G_\theta$ is isomorphic as a vector space to $\G$, for $e \in \G$, we will denote its image under this isomorphism in $\G_\theta$ by $e_\theta$ from now on.

As in the case of deformed algebras, $\G_\theta$ turns out to be a $\IT^n$ smooth bimodule. In particular, a deformed bimodule admits a $\IT^n$ action $ \beta^\theta$ induced from the action $ \beta$ on the bimodule prior to deformation as follows:
\begin{equation}\label{actionondeformed}
 \beta^\theta_t(e_\theta)=\sum_{\underline{k}}\chi_{\underline{k}}(t)e_{\underline{k}} \ \forall \ t \in \IT^n. 
\end{equation}

\brmrk \label{4thjune2018remark}
It is easy to see that $(\D_\theta)_{\underline{0}} \cong \D_{\underline{0}}$ as algebras. Moreover, $(\G_\theta)_{\underline{0}} \cong \G_{\underline{0}}$ as $\D_{\underline{0}}$ bimodules, since $\D_{\underline{0}}$ and $\G_{\underline{0}}$ are the fixed point subalgebra and the fixed point submodule under the $\IT^n$ action respectively. We also note that by \eqref{moduledeformationformula},  $(\G_\theta)_{\underline{0}} \subseteq \mathcal{Z}(\G_\theta)$ and in particular $(\D_{\theta})_{\underline{0}} \subseteq \mathcal{Z}(\D_\theta)$.
\ermrk

We have the following easy consequence of the definitions above: 

\blmma \label{equivarince}
Let ${\mathcal D}$ be a $\IT^n$ smooth algebra and
$\mathcal{G}_1, \mathcal{G}_2$ be $\IT^n$ smooth $\mathcal D$ bimodules, in the sense discussed above. Let $L : {\mathcal G_1} \rightarrow {\mathcal G_2}$ be a $\mathbb{T}^{n}$ equivariant continuous $\mathcal{D}$ bimodule map.
Then the underlying vector space map $L$ from $\G_1$ to $\G_2$ becomes a $\mathbb{T}^{n}$ equivariant continuous $\mathcal{D}_{\theta}$ bimodule map, denoted by $$L_{\theta}: ({\mathcal G_1})_\theta \rightarrow ({\mathcal G_2})_\theta$$
 defined by the equation
\begin{equation}
L_\theta(e_\theta)=(L(e))_\theta \ \forall \ e\in \G_1
\end{equation}
If $ L $ is a $ \mathcal{D} $-bimodule isomorphism, then $ L_\theta $
will be a $ \D_\theta $-bimodule
isomorphism. If $ \mathcal{G}_1 $ 
and $ \G_2 $ are algebras in particular, then $ L_\theta $ is an algebra homomorphism.  

Now suppose that $ {\rm Ker} ( L ) $ is complemented as a $ \D $-bimodule in $ {\mathcal G}_1, $ i.e, there exists a bimodule $ M \subseteq {\mathcal G}_1 $ such that $ {\mathcal G}_1 \cong {\rm Ker} ( L ) \oplus M. $ Then

1. $ {\rm Ker} ( L ) $ is invariant under the action of $ \mathbb{T}^n. $

2. $ M \cong {\rm Im} ( L ). $

3. If $ M $ is $ \mathbb{T}^n $ invariant, then $ { (\G_1)}_\theta = {\rm Ker} ( L_\theta ) \oplus M_\theta $ and $ M_\theta \cong {\rm Im} ( L_\theta ). $

4. If $ \G_2 = {\mathcal G_1} $ and $ L $ is an idempotent, then $ L_\theta $ is also idempotent.   

\elmma

The following lemma will also be of use to us.

\blmma
\label{12345}
Let ${\mathcal D}$ be an algebra equipped with ${\mathbb T}^n$-action and
$\mathcal{G}_1, \mathcal{G}_2$ be equivariant $\mathcal{D}$ bimodules, in the sense discussed above. Then
$ ( \mathcal{G}_1 \otimes_{\mathcal{D}} \mathcal{G}_2 )_{\theta} \cong (\mathcal{G}_1)_{\theta} \otimes_{{\mathcal D}_\theta} (\mathcal{G}_2)_{\theta} $ as $ {\mathcal D}_\theta $ bimodules.
\elmma	

Now we recall the Connes-Landi deformation (\cite{connes_landi}) of a spectral triple and its associated space of forms. We will work in the set up of Proposition \ref{equivarianttriple}. In particular, $ \A = C^\infty(M)$ and $ \E=\oneform $. By Sections 2 and 3 of \cite{scalar_5}, $\A$ and $\E$ can be deformed to the algebra $\A_\theta$ and the bimodule $\E_\theta$ respectively. Moreover, we have the following:

		\bthm \label{25thfebtheorem}
		With the algebra structure of $ \A_\theta $ as in \eqref{deformationformula}, we define $ \pi^\theta : \E \rightarrow \clb ( \clh ) $ by 
		\be \label{23feb2018} \pi^\theta ( e ) ( h ) = \sum_{\underline{m}, \underline{n} \in \mathbb{Z}^n} \chi_\theta ( \underline{m}, \underline{n} ) e_{\underline{m}} ( h_{\underline{n}} ),\ee
	 where  $e = \sum_{\underline{m}} e_{\underline{m}}, ~ h = \sum_{\underline{n}} h_{\underline{n}}.$
	\be \label{thirdpicture} {\rm Then} ~	 \pi^\theta ( \E ) \cong \E_\theta.\ee 
	We note that an analogous formula defines a representation of $ \A_\theta $ on $ \clh,$ to be denoted by $ \pi_\theta $ again. Also, 
	 $(\A_\theta, {\mathcal H}, d + d^*)$ defines a spectral triple. 
		
		Moreover, $ {\Omega^1} ( \A_\theta ) $ and ${\Omega^2} ( \A_\theta )$ are canonically isomorphic as $\A_\theta$ bimodules
	  with $ \E_\theta $ and $ ( {\Omega^2} ( \A ) )_\theta$ respectively. If $\delta: \A \raro \A$
	 denotes the map which sends $a$ to $[d + d^*,a],$ then we have a deformed map $ {\delta}_\theta $ from $ \A_\theta $ to $\E_\theta.$ 
	\ethm
		
		{\bf Proof:} For the proof that $(\A_\theta, {\mathcal H}, d + d^*)$ is a spectral triple, we refer to \cite{Connes-dubois}. For the isomorphism $ \pi^\theta ( \E) \cong \E_\theta, $ we refer to Proposition 2.8 of \cite{scalar_5}.			
			We refer to \cite{scalar_5} for the fact that $\E=\oneform$ and $\twoform$ can be deformed. Then the isomorphism follows by using the identification of \eqref{thirdpicture}. The last assertion follows by observing that the map $ \delta $ is $ \mathbb{T}^n $ equivariant. \qed

		Henceforth we will make the identifications
	  $\E_\theta \cong {\Omega^1} ( \A_\theta ), ~ {\Omega^2} ( \A_\theta ) \cong ( {\Omega^2} ( \A ) )_\theta$ without explicitly mentioning.

\subsection{The canonical Riemannian  bilinear metric on $\E_\theta$} \label{bimetricexistencedef}

In this subsection, we prove that the prescription of Subsection \ref{assumptionsonsptriple} is indeed a Riemannian bilinear metric on $\E_\theta.$ We prove this in two steps. In the first step, we deform the $\A$-bilinear map $g$ to an $\A_\theta$-bilinear map $ g_\theta $ and show that $ g_\theta $ is a Riemannian bilinear metric. In the second step, we show that the $\A$-bilinear map obtained from Lemma \ref{22ndjune2018} ( for the spectral triple $ ( \A_\theta, \mathcal{H}, D  ) $ ) coincides with $g_\theta.$ 

Let us recall the following definition:

\bdfn \label{actiononhom}
Let $ \G_1 $ and $ \G_2 $ be two $\mathcal{D}$ bimodules admitting actions by a group $ \IT^n $ and denoted by $  \beta_1 $ and $  \beta_2 $ respectively. Then ${\rm Hom}_{\mathcal {A}}(\mathcal{G}_1,\mathcal{G}_2)$
admits a natural $ \IT^n $ action $ \gamma $ defined by 
$$(\gamma_t.T)(e)= ( \beta_2)_t.(T({( \beta_1)_t}^{-1}.e)). $$
Here, $ t, T $ and $e $ belong to $ \IT^n, {\rm Hom}_{\mathcal{D}} ( \G_1, \G_2 ) $ and $ \G $ respectively.
\edfn 

\blmma  
In the set up of Definition \ref{actiononhom}, assume furthermore that $ \D $ admits an action $ \alpha $ of $ \IT^n $ and $  \beta_1,  \beta_2 $ are both $ \alpha $ equivariant. Then, for $ a $ in $ \D $, $ \omega $ in $ \G_1 $ and $T\in {\rm Hom}_{\mathcal{D}} ( \G_1, \G_2 ) $, we have 
$$ \gamma_t ( T a ) ( \omega ) =  \gamma_t ( T ) ( \alpha_t ( a )  ) ( \omega ) ~~ {\rm and } ~~ \gamma_t ( a T ) ( \omega ) = \alpha_t ( a ) (\gamma_t ( T ) ( \omega )) $$
\elmma
{\bf Proof:} We compute
\begin{align*}
 \gamma_t ( T a ) ( \omega ) &= ( \beta_2)_t ( ( T a ) ( ( \beta_1)^{- 1}_t ( \omega ) ) ) = ( \beta_2)_t T ( a ( \beta_1)^{- 1}_t ( \omega ) )\\
 &= ( \beta_2)_t T ( ( \beta_1)_{z^{- 1}} ( ( \alpha_t a ) \omega )  ) = \gamma_t ( T ) ( ( \alpha_t a ) \omega )\\
 &=  \gamma_t ( T ) ( \alpha_t ( a )  ) ( \omega ).
\end{align*}

The other equality follows similarly. \qed

\bppsn \label{dualdeformation}
$ ( \E^* )_{\theta} \cong ( \E_\theta )^* $ as $ \A_\theta$ bimodules. 
\eppsn
{\bf Proof:} We define a map $ T^\E_\theta: ( \E^* )_\theta \rightarrow ( \E_\theta )^* $ by 
$$ (T^\E_\theta (\phi_\theta)) ( e_\theta ) = (\phi(e))_\theta, $$
where we have used the fact that any arbitrary $\D$-$\D$ bimodule $\G$ is isomorphic as a vector space with the twisted bimodule $\G_\theta$, and have denoted the image of an element $f$ in $\G$ under the isomorphism by $f_\theta$. Let us denote the right $\A_\theta$ module structure of $\E_\theta$ by $\times_\theta$.

Since $ (T^\E_\theta (\phi_\theta))(e_\theta \times_\theta a_\theta)=(T^\E_\theta (\phi_\theta))((ea)_\theta)=(\phi(ea))_\theta=(\phi(e))_\theta\times_\theta a_\theta= (T^\E_\theta (\phi_\theta))(e_\theta)\times_\theta a_\theta$, we have that $(T^\E_\theta (\phi_\theta)) \in (\E_\theta)^*$.
That $T^\E_\theta$ is right $\A_\theta$-linear can be seen from the following.
\begin{align*}
(T^\E_\theta(\phi_\theta\times_\theta a_\theta))(e_\theta)&=(T^\E_\theta(\phi(a))_\theta)(e_\theta)\\
&=((\phi(a)(e))_\theta\\
&=(\phi(ae))_\theta\\
&=(T^\E_\theta(\phi_\theta))((ae)_\theta)\\
&=(T^\E(\phi_\theta))(a_\theta\times_\theta e_\theta)
\end{align*}

Let $ \gamma $ denote the action of $ \mathbb{T}^n $ on $ \E^* := {\rm Hom}_\A(\E,\A)$ defined by Definition \ref{actiononhom}. The $ \mathbb{T}^n $ actions on $ \E_\theta $ and $ ( \E^* )_\theta $ will be denoted by $ \beta^\theta $ and $ \gamma^\theta $ respectively as in \eqref{actionondeformed}. Moreover, the $ \mathbb{T}^n $ action on $ ( \E_\theta )^* := {\rm Hom}_{\A_\theta}(\E_\theta,\A_\theta)$ as obtained from Definition \ref{actiononhom} will be denoted by $ \gamma^\prime. $ We claim that the map $ T^\E_\theta $ is equivariant w.r.t the $ \IT^n $-action on $ ( \E^* )_\theta $ and $ ( \E_\theta )^*, $ i.e,
\be \label{24thjune} T^\E_\theta ( \gamma^\theta_t ( \phi_\theta ) ) = \gamma^\prime_t ( T^\E_\theta (\phi_\theta) ).\ee 
Indeed, we have: 
\begin{eqnarray*}
	T^\E_\theta ( \gamma^\theta_t ( \phi_\theta ) ) ( e_\theta ) &=& T^\E_\theta ( (\gamma_t ( \phi ))_\theta ) ( e_\theta )\\
	&=& ((\gamma_t (\phi ))(e))_\theta\\
	&=& ( \alpha_t (\phi ( \beta_{z^{- 1}} (e) ) ))_\theta \\
	&=& \alpha^\theta_t (\phi( \beta_{z^{- 1}} (e )))_\theta\\
	&=& \alpha^\theta_t ( T^\E_\theta ( \phi_\theta ) ( e_\theta ) )\\
	&=& (\gamma^\prime_t ( T^\E_\theta ( \phi_\theta ))) ( e_\theta ) .
\end{eqnarray*}
This proves \eqref{24thjune}. 

Thus, we have a well defined equivariant morphism 
$$ T^{\E_\theta}_{- \theta} : (  (\E_\theta)^* )_{- \theta} \rightarrow ( ( \E_\theta )_{- \theta} )^* \cong \E^*, $$
and subsequently, a morphism
$$ ( T^{\E_\theta}_{- \theta} )_{\theta} : ( \E_\theta )^* \cong ( (  (\E^*)_\theta )_{- \theta} )_{\theta} \rightarrow  ( \E^* )_\theta. $$
Finally, it is easy to check that the maps $ T^\E_\theta $and $ ( T^{\E_\theta}_{- \theta} )_{\theta} $ are inverses of one another. This finishes the proof. \qed

Let $\beta$ denote the action of $ \IT^n $ on $ C^\infty ( M ). $ Since $\IT^n$ acts on $M$ by isometries, the Riemannian metric $g$ is equivariant under the $\IT^n$ action i.e, for all $ \omega, \eta $ in $ \E , $ we have 
\be \label{gisinvariant}
g ( \beta_t ( \omega ) \otimes_{\A} \beta_t ( \eta ) ) = \alpha_t ( g ( \omega \otimes_{\A} \eta ) ). 
\ee

By using the $ \IT^n $ invariance of $ g, $ it is easy to see that the map $ V_g : \E \rightarrow {\E}^{\ast} $ is $ \IT^n $ equivariant. Thus, we have a $ \clc_\theta $-bimodule isomorphism $ ( V_g )_\theta $ from $ \E_\theta $ to $ ( \E^{\ast} )_{\theta}. $  

We view $ g $ as an element of ${\rm Hom}_{\A} (\E \otimes_{\A} \E, \A).$ The bimodule $ \E \otimes_{\A} \E $ is equipped with the natural diagonal action $ \beta \times \beta $ of $ \IT^n.$ Therefore, by Definition \ref{actiononhom}, we have an action of $ \IT^n $ on ${\rm Hom}_{\A} (\E \otimes_{\A} \E, \A).$ Since by \eqref{gisinvariant} $ g $ is equivariant, we have a deformed map 
$ g_{\theta} \in {\rm Hom}_{\A_\theta} ( {( \E \otimes_{\A} \E)}_{\theta} , \A_\theta ). $ However, by Lemma \ref{12345}, $ {( \E \otimes_{\A} \E)}_{\theta} \cong \E_{\theta} \otimes_{\A_\theta} \E_{\theta}. $ Thus, we have a map $g_\theta \in {\rm Hom}_{\A_\theta}(\E_\theta \otimes_{\A_\theta} \E_\theta,\A_\theta)$ which is the candidate for the Riemannian metric on $ \E_\theta = {\Omega^1} ( \clc_\theta ) $.

Our next result connects $(V_g)_\theta$ with $V_{g_\theta}$.

\bppsn \label{vgthetaisom}
\be \label{vgoneo} T^\E_\theta \circ (V_{g} )_{\theta} = V_{g_{\theta}} \ee 
and hence the map $ V_{g_{\theta}} : \E_{\theta} \rightarrow ( \E_{\theta} )^* $ is an isomorphism. 
\eppsn

{\bf Proof:} 
By the equivariance of $ g, $ it easily follows that the map $ V_g $ is equivariant and hence the map $ ( V_g )_{\theta} $ is an element of $ {\rm Hom}_{\Atheta} ( \E_\theta, ( \E^* )_\theta ) .$ By Proposition \ref{dualdeformation}, for all $e_\theta \in \E_\theta$, $ (T^\E_\theta \circ  (V_g )_\theta) (e_\theta) = (V_g(e))_\theta $ . For any element $f_\theta \in \E_\theta$, $(V_g(e))_\theta(f_\theta)=(g(e\tensora f))_\theta = g_\theta(e_\theta \otimes_{\A_\theta} f_\theta)=(V_{g_\theta}(e_\theta))(f_\theta)$. Thus, \eqref{vgoneo} holds. Moreover, since $ V_g $ is an isomorphism from $ \E $ to $ \E^* $,  Lemma \ref{equivarince} implies that  $ ( V_{g} )_\theta $ is an isomorphism from $ \E_\theta $ to $(\E^*)_{\theta} $. As $T^\E_\theta$ is an isomorphism from $(\E^*)_\theta$ to $(\E_\theta)^*$, the isomorphism of $ V_{g_\theta} $ follows from \eqref{vgoneo}. \qed

\bppsn
$ g_{\theta} $ is a noncommutative Riemannian bilinear metric on $ \E_\theta.$
\eppsn
{\bf Proof}: Clearly, $\sigma$ ($=2P_{\rm sym}-1$) is $ \IT^n $ equivariant, and as $g \circ \sigma =g$, we have $g_\theta \circ \sigma_\theta=g_\theta$ too,
i.e. $g_\theta$ is symmetric. It is also clear that $g_\theta$ is a bimodule map. Finally, by Proposition \ref{vgthetaisom}, $ V_{g_{\theta}} $ is nondegenerate. \qed

\bppsn
 Let $ g^\prime_\theta: \Etheta \otimes_{\Atheta} \Etheta \rightarrow \Atheta^{\prime \prime} $ be the $\Atheta$-bilinear map from Lemma \ref{22ndjune2018}. Then $ g^\prime_\theta = g_\theta $ and hence $ g^\prime_\theta $ is a Riemannian bilinear metric on $\Etheta.$ 
\eppsn
{\bf Proof:} Let $\omega = [D,a_1]a_2$ and $\eta = [D,b_1]b_2$ be elements in $\E$ to be viewed as elements of $\mathcal{B}(\mathcal{H})$. Let us denote the images of $\omega$ and $\eta$ in $\E_\theta$ by $\omega_\theta$ and $\eta_\theta$ respectively. Similarly, the representation of $\A_\theta$ in $\mathcal{B}(\mathcal{H})$ will be denoted by $\pi_\theta$. Finally, $\tau_\theta$ will denote the state on $\mathcal{B}(\mathcal{H})$ (for the spectral triple $(\A_\theta,\mathcal{H},D)$ as in Subsection \ref{assumptionsonsptriple}).
Then, if $p$ is the dimension of the manifold $M$, we compute
\begin{align*}
	\tau_\theta(\omega_\theta \eta_\theta \times{}_\theta a_\theta)&={\rm Lim}_{\omega} \frac{{\rm Tr}([D,\pi_\theta(a_1)]\pi_\theta(a_2)[D,\pi_\theta(b_1)]\pi_\theta(b_2)\pi_\theta(a)|D|^{-p})}{{\rm Tr}(|D|^{-p})}\\
	&={\rm Lim}_{\omega}\frac{{\rm Tr}([D,a_1]a_2[D,b_1]b_2a|D|^{-p})}{{\rm Tr}(|D|^{-p})} \text{ (by Proposition 4.4.2 of \cite{thesis})}\\
	&= \tau(\omega \eta a)\\
	&=\tau(g(\omega \tensora \eta)a)\\
	&=\tau_\theta(g(\omega \tensora \eta_\theta)\pi_\theta(a)) \text{ (by Proposition 4.4.2 of \cite{thesis})}\\
	&=\tau_\theta(g_\theta(\omega_\theta \otimes_{\A_\theta} \eta_\theta)\times{}_\theta a_\theta)	
\end{align*} 
This proves that the bilinear form of Lemma \ref{22ndjune2018} for the spectral triple $(\A_\theta,\mathcal{H},D)$ is equal to $g_\theta$ and hence it satisfies all the conditions of Definition \ref{defn25thmay2018}. \qed

\subsection{Existence and uniqueness of Levi-Civita connections for the Connes-Landi deformed spectral triple} \label{rieffelsubsection}

  We will continue to use the notations introduced in Definition \ref{defn22ndfeb2017}. The goal of this subsection is to apply the results deduced in the last two subsections for proving Theorem \ref{existencerieffeldeformation}.

\blmma \label{starprime}
$\E_\theta$ is a finitely generated projective right module over $\A_\theta$.
\elmma
{\bf Proof:} By Lemma \ref{free_case}, $\E_{\underline{0}}$ is a finitely generated projective right $\A_{\underline{0}}$ module. Then $\E_{\underline{0}}\otimes_{\A_{\underline{0}}}\A$ is a finitely generated projective right $ \A$ module. Since the isomorphism $u^\E_{\E_{\underline{0}}}:\E_{\underline{0}}\otimes_{\A_{\underline{0}}}\A \rightarrow \E$ as given by Lemma \ref{isomorphismofEzero} (and similarly the isomorphism  $ \A_{\underline{0}}\otimes_{\A_{\underline{0}}}\A \rightarrow \A$) is $\IT^n$ equivariant, $\E_\theta \cong (\E_{\underline{0}} \otimes_{\A_{\underline{0}}} \A)_\theta \cong \E_{\underline{0}} \otimes_{\A_{\underline{0}}}\A_\theta$ is finitely generated as a right $\A_\theta$ module. Here, we have used the fact that $(\E_{\underline{0}})_\theta \cong \E_{\underline{0}}$ as right $\A_{\underline{0}}$ modules since $\E_{\underline{0}}$ is the fixed point submodule for the action of $\IT^n$. \qed

\blmma \label{4thjune20182}
The map $u^{\E_\theta}_{\E_{\underline{0}}} = (u^\E_{\E_{\underline{0}}})_\theta : \E_{\underline{0}} \otimes_{\A_{\underline{0}}} \A_\theta \rightarrow \E_\theta$ is an isomorphism. Moreover, the map $u^{\E_\theta}:\mathcal{Z}(\E_\theta)\otimes_{\mathcal{Z}(\A_\theta)} \A_\theta \rightarrow \E_\theta$ is an isomorphism.
\elmma
{\bf Proof:} By Lemma \ref{isomorphismofEzero}, the $\IT^n$ equivariant map $u^\E_{\E_{\underline{0}}}:\E_{\underline{0}}\otimes_{\A_{\underline{0}}}\A \rightarrow \E$ is an isomorphism. Hence, by Lemma \ref{equivarince} and Lemma \ref{12345} the map $(u^\E_{\E_{\underline{0}}})_\theta:\E_{\underline{0}}\otimes_{\A_{\underline{0}}}\A_\theta \rightarrow \E_\theta $ is an isomorphism.

For the second assertion, we note that by Lemma \ref{free_case}, $\F_{\underline{0}}=\E_{\underline{0}}=(\E_\theta)_{\underline{0}}$ is finitely generated projective over $\C_{\underline{0}}=\A_{\underline{0}}=(\A_\theta)_{\underline{0}}$. By Remark \ref{4thjune2018remark}, $(\A_\theta)_{\underline{0}} \subseteq \mathcal{Z}(\A_\theta)$  and $(\E_\theta)_{\underline{0}} \subseteq \mathcal{Z}(\E_\theta)$. Therefore, by Proposition \ref{isomorphismofEprime}, we conclude that the map $u^{\E_\theta}:\mathcal{Z}(\E_\theta)\otimes_{\mathcal{Z}(\A_\theta)} \A_\theta \rightarrow \E_\theta$ is an isomorphism.
 \qed

\blmma \label{kermsplit}
The bimodule $ \E_\theta \otimes_{\A_\theta} \E_\theta$ admits a decomposition $\E_\theta \otimes_{\A_\theta} \E_\theta = {\rm Ker} ( \wedge_\theta ) \oplus \F_\theta$ of right $\A_\theta$ modules, where $\F_\theta \cong \Omega^2(\A_\theta)$ is satisfied.
\elmma
{\bf Proof:} This follows by applying Lemma \ref{equivarince}, Lemma \ref{12345} and Corollary \ref{Finvariant} applied to the $ \IT^n $ equivariant map $\wedge. $ \qed

\blmma \label{deformedsigmacan}
We have  $ \sigma_\theta=\sigma^{\rm can}$ where $ \sigma^{\rm can} $   is the map given by Theorem \ref{skeide2}, i.e. for all $e_\theta$, $f_\theta$ in $\mathcal{Z}(\E_{\theta})$ and $a$ in $\A_\theta$,
\[ \sigma_\theta( e_\theta \otimes_{\A_\theta} f_\theta a) = f_\theta \otimes_{\A_\theta} e_\theta a. \]
\elmma
{\bf Proof:} We begin by observing that the map $ \sigma^{\rm can} $ makes sense as $ \E_\theta $ is a centered bimodule since $u^{\E_\theta} $ is an isomorphism by Lemma \ref{4thjune20182}.  For the purpose of the proof of this lemma, we use an equivalent description of Rieffel deformation as described in \cite{Connes-dubois}. 
Let $C^{\infty}(\mathbb{T}^{n}_{\theta})$ be the Fre\'chet   algebra corresponding to the noncommutative $n$-torus, where the deformation parameter is given by a real,
 $n \times n$ skew symmetric matrix $\theta=(( \theta_{kl} ))$. We denote the canonical ${\mathbb{T}}^n$-action on this algebra by $\nu$ and let $\A({\mathbb{T}}^n_\theta)$ denote the canonical
 `polynomial subalgebra', i.e. the dense unital $\ast$-subalgebra of $ C^{\infty}(\mathbb{T}^{n}_{\theta}))$ generated by polynomials in $U_i$ and their inverses, where $U_i, i=1, \ldots,n$
  are the canonical unitary `coordinates', i.e. elements satisfying 
	$ U_k U_l={\rm exp}(2 \pi i \theta_{kl}) U_l U_k.$
	For $ \underline{m} = ( m_1, m_2, \cdots m_n ) \in \IZ^n, $ we write $ U^\theta_{\underline{m}} = U^{m_1} \cdots U^{m_n}. $ Hence, $ U^\theta_{\underline{m}} $ are   the canonical generators of $ ( \A ( \IT^n_\theta ) )_{\underline{m}}. $
	
In this picture, $ \E_\theta $ is identified with a suitable Fre\'chet  completion of the fixed point submodule 	
 $({\mathcal E} {\otimes}_{\mathbb C}\A(\mathbb{T}^{n}_{\theta}))^{\beta\times\nu^{-1}}$
 of the module ${\mathcal E} {\otimes}_{\mathbb C} \A(\mathbb{T}^{n}_{\theta})$ with respect
to the action $\beta \times \nu^{-1}$. 
We can identify an element $ e_\theta \in ( \E_\theta )_{\underline{m}} $ with the element $ e \tensorc U^\theta_{\underline{m}} \in ( \E \tensorc \A({\mathbb{T}}^n_\theta) )^{\beta \times \nu^{-1}}.$
The action of the group $ \IT^n $  on an element $ e_\theta \in ( \E_\theta )_{\underline{m}} $ is given by the formula: 
$$ \beta^\theta_t ( e_\theta  ) = \beta^\theta_t (  e \tensorc U^\theta_{ \underline{m}} ) := \beta_t ( e ) \tensorc U^\theta_{ \underline{m}}. $$
Therefore, if $ e_\theta \in ( \E_\theta )_{\underline{m}} $ and $ f_\theta \in ( \E_\theta )_{\underline{n}}, $ we have
$$ \sigma_\theta ( e_\theta \otimes_{\Atheta} f_\theta  ) = \sigma ( e \tensora f ) \tensorc U^\theta_{\underline{m}} U^\theta_{\underline{n}}.$$
Thus, if $ e, f \in ( \E_\theta )_{\underline{0}}, $ we get
$$ \sigma_\theta ( e_\theta \otimes_{\Atheta} f_\theta ) = \sigma ( e \tensora f ) \tensorc 1 = ( f \tensora e ) \tensorc 1 = \sigma^{\rm can} ( e_\theta \otimes_{\Atheta} f_\theta ), $$
since $ (  \E_\theta  )_{\underline{0}} \subseteq \mathcal{Z} ( \Etheta ) $ by Remark \ref{4thjune2018remark}.

Now, by Lemma \ref{4thjune20182}, $ \E_{\underline{0}} $ is right $ \A_\theta $-total in $ \E_\theta $ and hence $\{e_\theta \otimes_{\A_\theta} f_\theta : e_\theta, f_\theta \in (\E_\theta)_{\underline{0}} \}$  is right $\A_\theta$-total in $\E_\theta\otimes_{\A_\theta}\E_\theta$. Thus, our claim follows from Proposition \ref{28thmay2018}. \qed

\vspace{4mm}

 Collecting the above results, we get the following:

\vspace{4mm}
 
{\bf Proof of Theorem \ref{existencerieffeldeformation}}
We start by recalling that we have already proved (Lemma \ref{starprime}) that $\E_\theta$ is a finitely generated projective right module over $\A_\theta$. By Lemma \ref{4thjune20182}, the map $u^{\E_\theta}:\mathcal{Z}(\E_\theta)\otimes_{\mathcal{Z}(\A_\theta)} \A_\theta \rightarrow \E_\theta$ is an isomorphism.
Next,  ${\rm Ker} ( \wedge_\theta ) $ is complemented in $ \E_\theta \otimes_{\A_\theta} \E_\theta $ by Lemma \ref{kermsplit}. 
Lastly, the equality $\sigma_\theta=\sigma^{\rm can}$ follows from Lemma \ref{deformedsigmacan}. 
Thus we have shown that all the hypotheses of Theorem \ref{lcexistenceforbilinear} are satisfied and this completes the proof. \qed

\vspace{4mm}
 
 \brmrk 
 For the deformed spectral submodule $\F_\theta$, analogues of the results Lemma \ref{starprime}, Lemma \ref{4thjune20182}, Lemma \ref{kermsplit} and Lemma \ref{deformedsigmacan} are proved the same way. Hence the analogous result of Theorem \ref{existencerieffeldeformation} also holds for the deformed module $\F_\theta.$
 \ermrk
 
 \bcrlre
 Under the assumptions of Theorem \ref{existencerieffeldeformation}, the Levi-Civita connection $\nabla$ on the bimodule $\E$ deforms to the Levi-Civita connection $\nabla_\theta$ on $\E_\theta.$
 \ecrlre
 
 {\bf Proof:} Let us observe that given a $\IC$-linear map $T$ between $\IT^n$ smooth bimodules $\G_1$ and $\G_2$, there exists a $\IC$ linear deformation $$T_\theta: (\G_1)_\theta \to (\G_2)_\theta$$ defined by $T_\theta(e_\theta)=(T(e))_\theta$. Here, $e \in \G_1$ and $e_\theta$ is the image of $e$ in $(\G_1)_\theta$. Moreover, if $T$ is $\IT^n$ equivariant, then $T_\theta$ is $\IT^n$ equivariant.\\
 Since the spectral triple $(\A,\mathcal{H},D)$ is $\IT^n$ equivariant, the maps $d: \A \to \E$ and $d: \E \to \twoform$ are $\IT^n$ equivariant. It is easy to see that the map $\A_\theta \to \E_\theta$ given by $a_\theta \mapsto [D,\pi_\theta(a_\theta)]$ is nothing but the deformation of the map $d: \A \to \E$. Moreover, the maps $d_\theta:\A_\theta \to \E_\theta$ and $d_\theta: \E_\theta \to \twoformdeformed$ are $\IT^n$ equivariant.\\
 Since the map $\nabla$ is the Levi-Civita connection, $\nabla$ is $\IT^n$ equivariant. Thus, we have a $\IC$ linear map $\nabla_\theta:\E_\theta \to (\E \tensora \E)_\theta \cong \E_\theta \otimes_{\A_\theta} \E_\theta$ and it can be easily checked that $\nabla_\theta$ is a connection.\\
 By Lemma \ref{mequivariant}, $\wedge: \E \tensora \E \to \twoform$ is a $\IT^n$ equivariant $\A$ bimodule map. Hence, $\wedge_\theta: \E_\theta \otimes_{\A_\theta} \E_\theta \to \twoformdeformed$ is defined, and $\wedge_\theta \circ \nabla_\theta = (\wedge \circ \nabla)_\theta = -d_\theta$. Therefore, $\nabla_\theta$ is a torsionless connection.\\
 Lastly we show that $\nabla_\theta$ is compatible with the metric $g_\theta$. We need to show that $\Pi_{g_\theta}(\nabla_\theta)=d_\theta g_\theta$. However, by Lemma \ref{proppignabla}, the map $\Pi_{g_\theta}(\nabla_\theta) -d_\theta g_\theta$ is right $\A_\theta$ linear. Since $\{ \omega_\theta \otimes_{\A_\theta} \eta_\theta : \omega_\theta, \eta_\theta \in \mathcal{Z}(\E_\theta) \}$ is right $\A_\theta$ total in $\Etheta \otimes_{\A_\theta} \Etheta$, it is enough to show that for all $\omega_\theta, \eta_\theta \in \mathcal{Z}(\E_\theta)$, we have $(\Pi_{g_\theta}(\nabla_\theta)-d_\theta g_\theta)(\omega_\theta \otimes_{\A_\theta} \eta_\theta)=0$ for all $\omega_\theta, \eta_\theta \in \mathcal{Z}(\E_\theta)$. Let $\omega_\theta, \eta_\theta \in \mathcal{Z}(\E_\theta)$. Then,
 \begin{align*}
 (\Pi_{g_\theta}(\nabla_\theta))(\omega_\theta \otimes_{\A_\theta} \eta_\theta)&=(g_\theta \otimes_{\A_\theta} {\rm id}_{\E_\theta})(\sigma_\theta)_{23}(\nabla_\theta(\omega_\theta)\otimes_{\A_\theta} \eta_\theta + \nabla_\theta(\eta_\theta)\otimes_{\A_\theta} \omega_\theta)\\
 &=((g \tensora {\rm id}_{\E})\circ \sigma_{23})_\theta(\nabla(\omega)\otimes_{\A} \eta + \nabla(\eta)\otimes_{\A} \omega)_\theta\\
 &=\big(((g \tensora {\rm id}_{\E})\circ \sigma_{23})(\nabla(\omega)\otimes_{\A} \eta + \nabla(\eta)\otimes_{\A} \omega)\big)_\theta\\
 &=(\Pi_g(\nabla)(\omega \tensora \eta))_\theta\\
 &=(-dg(\omega \tensora \eta))_\theta\\
 &=-d_\theta g_\theta(\omega_\theta \otimes_{\A_\theta} \eta_\theta).
 \end{align*}
Therefore, $\nabla_\theta$ is compatible with the metric $g_\theta$. \qed

\vspace{4mm}

{\bf Acknowledgement:} D.G. is partially supported by the J.C. Bose National Fellowship.

\end{document}